\documentclass[12pt,a4paper,twoside,final,notitlepage, leqno]{article}
\usepackage[english]{babel}
\usepackage[T1]{fontenc}
\usepackage{epsfig, graphicx, amssymb}
\usepackage{amsmath,amsthm,epsfig,amsfonts}
\usepackage{float}
\usepackage{color}
\newcommand{\monitem}{ \smallskip \noindent $\bullet$ \quad  }
\newcommand{\moneq}{\vspace*{-7pt} \begin{equation} \displaystyle }
\newcommand{\moneqstar}{\vspace*{-6pt} \begin{equation*} \displaystyle }
\newcommand{\monendstar}{\vspace*{-6pt} \end{equation*}   }
\newcommand{\monend}{\vspace*{-7pt} \end{equation}   }
\newcommand{\moneqarraystar}{ \begin{eqnarray*} \displaystyle }
\newcommand{\monendarraystar}{ \end{eqnarray*}   }
\newcommand{\dd}{{\rm d}}

\newcommand{\RR}[0]{\mathbb{R}}

\newcommand{\ZZ}[0]{\mathbb{Z}}
\newcommand{\NN}[0]{\mathbb{N}}

\definecolor{vertfonce}{rgb}{0.0, 0.5, 0.0}
\def\section*#1{}
\usepackage{fancyhdr}
\renewcommand{\headrulewidth}{0pt}
\begin{document}

\fancypagestyle{plain}{ \fancyfoot{} \renewcommand{\footrulewidth}{0pt}}
\fancypagestyle{plain}{ \fancyhead{} \renewcommand{\headrulewidth}{0pt}}
%%   \bibliographystyle{alpha}

%%%%%%%%%%%%%%%%%%%%%%%%%%%%%%%%%%%%%%%%%%%%%%%%%%%%%%%%%%%%%%%%%%%%%%%%%%%%%%%
~

  \vskip 2.1 cm

\centerline {\bf \LARGE On single distribution lattice Boltzmann schemes}

\bigskip

\centerline {\bf \LARGE for the approximation of Navier Stokes equations}

 \bigskip  \bigskip \bigskip

\centerline { \large    Fran\c{c}ois Dubois$^{abc}$ and Pierre Lallemand$^{d}$}

\smallskip  \bigskip

\centerline { \it  \small
  $^a$   Laboratoire de Math\'ematiques d'Orsay, Facult\'e des Sciences d'Orsay,}

\centerline { \it  \small   Universit\'e Paris-Saclay, France.}

\centerline { \it  \small
$^b$    Conservatoire National des Arts et M\'etiers, LMSSC laboratory,  Paris, France.}

\centerline { \it  \small
  $^c$    International Research Laboratory 3457, Centre National de la Recherche Scientiﬁque,}

\centerline { \it  \small
Centre de Recherches Math\'ematiques, Universit\'e de Montr\'eal, Montr\'eal, QC, Canada.}

\centerline { \it  \small
 $^d$ Beijing Computational Science Research Center, Haidian District, Beijing 100094,  China.}

%%%   \bigskip

\bigskip  \bigskip

\centerline {October 2023$\,$\footnote {\rm  \small  $\,$ 
{\it Communications in Computational Physics}, volume 34, pages 613-671, 2023. 
A preliminary version of this contribution was first  presented %%%at the Onera seminar, 19 June 2020.
at the International Conference of Mesoscopic Methods in Engineering and Science, Newark (Delaware, USA) the 
11~July 2018. Edition October 2024.}}
%%%    {\footnote

 \bigskip \bigskip
 {\bf Keywords}: partial differential equations, asymptotic analysis

 {\bf AMS classification}:
%%    65Q05,   %% Difference and functional equations, recurrence relations,
 76N15,  %%%  Gas dynamics, general
 82C20.   %%%  Dynamic lattice systems (kinetic Ising, etc.) and systems on graphs

 {\bf PACS numbers}:
02.70.Ns, % meth particulaires
%%% 05.20.Dd, % theorie cinetique
47.10.+g  % generalites en meca flu
%%%%%%%%%%%%%  47.11.+j. % cfd

\bigskip  \bigskip
\noindent {\bf \large Abstract}

\noindent
In this contribution we study the formal ability of a
multi-resolution-times  lattice Boltzmann scheme to approximate 
isothermal and thermal compressible Navier Stokes equations
with a single particle distribution.
More precisely, we consider a total of 12 classical square lattice Boltzmann schemes
with  prescribed sets of conserved and nonconserved moments.
The question is to determine the algebraic expressions of the equilibrium functions for the nonconserved moments
and  the relaxation parameters associated to each scheme.
We compare  the fluid equations  
and the result  of the  Taylor expansion method at second order accuracy
for bidimensional examples with a maximum of 17 velocities and
three-dimensional schemes with at most 33 velocities.  
In some cases, it is not possible to fit exactly the physical model.
For several  examples, we adjust  the Navier Stokes equations
and propose nontrivial expressions for the equilibria.

\newpage 
%%%%%%%%%%%%%%%%%%%%%%%%%%%%%%%%%%%%%%%%%%%%%%%%%%%%%%%%%%%%%%%%%%%%%%%%%%%%%%%  section 1
\bigskip \bigskip    \noindent {\bf \large    1) \quad  Introduction} %%%  : etat de l art NS et LBM   }
%%%%%%%%%%%%%%%%%%%%%%%%%%%%%%%%%%%%%%%%%%%%%%%%%%%%%%%%%%%%%%%%%%%%%%%%%%%%%%%%%%%%%%%%%%

%%%%%%%%%%%%%%%%%%%%%%%%%%%%%%%%%%       31 decembre 2018       %%%%%%%%%%%%%%%%%%%%%%%%%
\fancyhead[EC]{\sc{Fran\c{c}ois Dubois and Pierre Lallemand}}
\fancyhead[OC]{\sc{Single lattice Boltzmann distribution for Navier Stokes equations}}
%%%%%%%%%%%%%%%%%%%%%%%%%%%%%%%%%%       31 decembre 2018       %%%%%%%%%%%%%%%%%%%%%%%%%
%%%%%%%%%%%%%%%%%%%%%%%%%%%%%%%%%%%%%%%%%%%%%  jolie numerotation des pages
\fancyfoot[C]{\oldstylenums{\thepage}}
%%%%%%%%%%%%%%%%%%%%%%%%%%%%%%%%%%%%%%%%%%%%%  fin jolie numerotation des pages

\noindent
The study of fluid mechanics is a natural problem set by the pioneers
of the lattice Boltzmann schemes in their modern  form
%%% with the contributions of d'Humi\`eres, Lallemand, Succi
(see \cite{HSB91, QHL92, DDH92, LL00} and many others). 
An underlying  lattice Boltzmann equation is discretised on a cartesian grid
with a finite set of velocities chosen in such a way that during one time step, an exact transport is done
between two  vertices of the mesh.
A lattice Boltzmann scheme is composed with two steps:
a nonlinear  local relaxation step,   followed by a linear advection scheme coupling
a given vertex with a given family of neighbours.
The relaxation step follows in general the approximation 
introduced by Bhatnagar, Gross and Krook \cite{BGK54}.

\noindent
From a completely defined lattice Boltzmann model (elementary velocity set, equilibrium
and relaxation rates) it is possible to determine the macroscopic behaviour of
the model through equivalent partial differential equations. One  tries to match
these equivalent equations to those governing a given physical situation.
Instead of matching equivalent partial differential equations and physical  partial differential equations,
term by term, it is common place to
compute the ``hydrodynamic'' modes of the two approaches and try to match them.

\noindent 
%This numerical method allows the simulation of an important number of physical phenomena as
This process allows the design of new simulation techniques 
for an important number of physical phenomena as
isothermal flows, compressible flows with heat transfer,  non-ideal fluids, 
multiphase and multi-component flows,
microscale gas flows, soft-matter flows,...  up to quantum mechanics.
The lattice Boltzmann method is inspired from a mesoscopic Boltzmann model
but is not able  in general to solve the Boltzmann equation or associated kinetic models.
It is admitted that conservative macroscopic models can be  approximated 
with the lattice Boltzmann schemes.
For the present status of the method and the various applications,
we refer {\it e.g.} to the books of Guo and Shu  \cite{GS13} 
or Kr\"uger {\it et al.} \cite{K3S2V17}.

%%%   \monitem Isothermal Navier Stokes
\noindent 
The simulation of incompressible flows, isothermal flows or thermal flows
with moderate compressible effects %%  with the Boussinesq  approximation
is very classical and we refer among others to the contributions   \cite{CCL19,HL97,GP18,GSW00,KT04,LL00,LLKY20,OF20,Su01}
%%%  Coreixas, Chopard, and Latt \cite{CCL19}
%%%  qui a propose le D3Q33 pour la premiere fois ?  Geier et Andrea Pasqualia \cite{GP18} l'utilisent
%%%   Otte and Frank \cite{OF20} observe that the D3Q33 admits a linear stability structure for BGK-type collision operators 
%%%   Kataoka and Tsutahara \cite{KT04} used D2Q9 and D3Q15 lattice Boltzmann models for simulating  compressible Euler equations. 
and to some  operational softwares like
OpenLB \cite{Kr10}, 
Palabos  \cite{LC07}
Powerflow \cite{CTM98}, 
LaBS-ProLB \cite{TRL13}
or pylbm \cite {GG17}. 
%%%  
In his  prospective article \cite{Su15},  Succi points the fact that the simulation of compressible flows
including the presence of two thermodynamic variables and eventually high Mach numbers is
one of the main open questions related to lattice Boltzmann schemes. 

%%%  \monitem 2 f
\noindent
Following an idea initially proposed in \cite{KP94,SC93},
a popular approach consists in adding a second particle distribution
to treat the conservation of energy. Following this framework, 
Guo {\it et al.} \cite{GZSZ07} use two particle distributions on the standard D2Q9 lattice
  to simulate compressible thermal flows at low Mach number, 
Nie {\it et al.} \cite{NSC08} propose a double discrete distribution
for thermal lattice Boltzmann model
using a three-dimensional %%% ninth-order accurate quadrature 
scheme employing 121 velocities, 
Latt {\it et al.} \cite{LCBP20} use a 
double-distribution-function based on the D3Q39 scheme 
for the simulation of polyatomic gases in the supersonic regime, 
Frapolli {\it et al.} \cite{FCK20} present a huge variety of three-dimensional lattice Boltzmann models
with two particle distributions to simulate compressible flows with the entropic lattice Boltzmann method. 
We proposed in~\cite{DGL16}
to recover  the full Navier Stokes equations in one space dimension using two lattice Boltzmann schemes and treating an entropy
equation with a second particle distribution.

%%%%%%%%%%%%%%%%%%%%%%%%   \monitem entre 1 f et 2 f    puis 1f 
\noindent
The approximation of thermal Navier Stokes equations with
lattice Boltzmann schemes involving only  one particle distribution
is not popular. Nevertheless, 
Sun and Hsu \cite{SH03} 
studied a three-dimensional compressible lattice Boltzmann model with large distribution velocity sets
using up to 96 velocities, 
Shan {\it al.}  \cite{SYC06} use  an  Hermite expansion
  and consider isothermal as well as thermal flows, 
Prasianakis and Karlin \cite{PK07} introduce two corrections to standard lattice Boltzmann schemes
to simulate thermal flows for two space dimensions,   
Yudistiawan {\it et al.} \cite{YKPA10} use variants of D3Q27 scheme 
and add some  degrees of freedom to the D3Q27 scheme for simulating
%%% isothermal and thermal
fluid flows with a single distribution, and 
Gillissen \cite{Gi16}, after Ricot {\it et al.}   \cite{RMSB09},  uses spatial filtering for stabilizing thermal Navier Stokes models
when using the D3Q33 scheme.

In  \cite{LL03},  Lallemand  and Luo put in evidence the stability difficulties
due to a merging of  viscous and thermal modes for moderate wave numbers. 
We did in \cite{LD15}
a linear analysis of D2Q13 lattice Boltzmann scheme  for advective acoustics and
tune  the parameters of the D2Q13  scheme for higher moments to improve stability. 
With this approach, it has been possible to simulate  the De Vahl Davis  \cite{DVD83}  thermal test case for natural convection 
at Rayleigh number = $ \, 10^{5} \, $ and  Prandtl number $ = 0.71 $.

\noindent
In this contribution, we consider classical schemes for two and three space dimensions,  
as the D2Q9 scheme (studied in detail in \cite{LL00}) and  
the D2Q13 scheme  as proposed in  \cite{Qi93}.
We are also interested with not so common schemes as the 
D2Q17  of Qian and Zhou \cite{QZ98}, 
D2V17 of Philippi {\it et al.}  \cite{PH06} and 
D2W17   of one of us  \cite{La10}.
For three space dimensions, the popular 
D3Q19 proposed by  d'Humi\`eres {\it et al.} \cite{DGKLL02} %% ok dans leur article de 2002
has got our interest. It is also the case for the  D3Q27 scheme~\cite{DL11,SKTC15}. 
It has been also necessary to consider the D3Q33 scheme,
intensively used in \cite{GP18}.
%%We are happy to put in evidence the interest af a simpler scheme
We put in evidence the interest of a simpler scheme
with 27 velocities, named in this contribution ``D3Q27-2'',
initially proposed for the simulation of viscoelastic behaviour~\cite{LHLR03}. 
%%%%   included in the simulation of
%%%%   viscoelastic behaviour.

\noindent
%
%%%{\color{red}
More precisely, we work with the paradigm of 
multi-relaxation-times lattice Boltzmann scheme proposed by d'Humi\`eres \cite{DDH92}.
For all these schemes, the moments are supposed to be given previously to our study; 
in particular the set of conserved moments that determine the conservation laws of fluid models, 
and also the microscopic moments that are not directly constrained by the physics. 
Then two sets of parameters determine completely the  lattice Boltzmann scheme:
the value of the microscopic moments at equilibrium and the relaxation parameters.
In this contribution, we use the general asymptotic  Taylor expansion method of analysis developed in  \cite{Du08, Du2022}
to establish equivalent partial differential equations at second order.
We constrain the study with the hypothesis that the lattice Boltzmann scheme has 
only a single particle distribution.
The result of this expansion is a very large algebraic expression of the derivatives of the
equilibrium functions and of the relaxation parameters. Our objective
in the present contribution is to match these
 equivalent partial differential equations first with the isothermal and secondly to the thermal
 Navier Stokes equations of a perfect gas for two and three space dimensions. 
%%% }
%
%%   \noindent
%%   in the particular case of a single particle distribution. 
%%   The equilibrium value of nonequilibrium moments is {\it a priori} completely undetermined. 
%%   With the Taylor expansion method \cite{Du08, Du2022}, we derive the equivalent partial
%%   differential equation of the lattice Boltzmann scheme.
%%   Then we confront this general differential model
%%   with the isothermal and thermal Navier Stokes equations for two and three space dimensions.
%%   \noindent
In some cases, it is not possible to fit exactly the physical model and we recover classical
equilibria that are correct up to third order relative to velocity.
In a certain number of examples, we fit exactly the Navier Stokes equations
and propose nontrivial equilibria.

\noindent
In Section~2, we recall our notations and hypotheses for the compressible Navier Stokes equations,
in the isothermal case or with thermal effects, in two or three space dimensions.
We describe in section~3 our approach for lattice Boltzmann schemes with multiple
relaxation times and in particular the asymptotic analysis.
Various schemes are proposed for the isothermal Navier Stokes equations
for two space dimensions in Section~4.
This  isothermal model is developed for the dimension 3 in Section~5.
Thermal Navier Stokes hypotheses are taken into account in Section~6
in the two-dimensional case.
Finally we present in Section~7 two lattice Boltzmann schemes that
are exactly second order accurate for thermal Navier Stokes in three dimensions. 

\noindent
A preliminary  version of this work including three-dimensional results
for thermal Navier Stokes equations  was proposed in march 2021 at the 
``Colloque des sciences math\'ematiques du Qu\'ebec'', Montr\'eal (Qu\'ebec, Canada)
during the delegation period of one of the authors in Montr\'eal~\cite{Du2021-montreal}.

%%%%%%%%%%%%%%%%%%%%%%%%%%%%%%%%%%%%%%%%%%%%%%%%%%%%%%%%%%%%%%%%%%%%%%%%%%%%%%%  section 2  
\bigskip \bigskip     \noindent {\bf \large    2) \quad  Compressible Navier Stokes equations}
%%%%%%%%%%%%%%%%%%%%%%%%%%%%%%%%%%%%%%%%%%%%%%%%%%%%%%%%%%%%%%%%%%%%%%%%%%%%%%%%%%%%%%%%%%

\noindent
We study the compressible Navier Stokes equations for gas dynamics in a very basic form.
The physical hypotheses are very classical and we refer {\it e.g.} to the treatises \cite{An81,GHP01,LL87,LR57}. 
To fix the ideas, we first consider  the algebraic formulas in one space dimension. %%  in order to simplify the presentation.
The conserved variables  can be stated as 
density $ \, \rho $, momentum $ \, j \equiv \rho \, u \, $ with the velocity $ \, u $, 
and  total energy 
\moneq \label{energie-totale}
E = {1\over2} \, \rho \, | {\bf u} |^2 + \rho \, e \,, 
\monend 
with $ \, e \, $ the internal energy of the fluid. 
We suppose a  polytropic perfect gas equation of state: 
$ \, p = (\gamma - 1) \, \rho \, e $,
with a  linear relation between internal energy and temperature $ \, T $,
and a constant specific heat $ \, c_v $:  $ \, e = c_v T $.
Moreover, the ratio  $ \, \gamma = {{c_p}\over{c_v}} \, $
between the specific heats is also supposed to be constant.
The viscosity $ \, \mu \, $ can be a function of the thermodynamic variables $ \, \rho \, $
and $ \, e \, $ and the Prandtl number is defined
from the  thermal conductivity  $ \, \kappa \, $  
according to  the relation $\, Pr = {{ \mu \, c_p}\over{\kappa}} $.
Then the conservations of mass, momentum and total energy take the form
\moneq \label{NS1d}    \left \{ \begin {array}{l}
 \partial_t \rho + \partial_x (\rho \, u) = 0 \\ 
\partial_t (\rho \, u) + \partial_x (\rho \, u^2 + p ) - \partial_x ( \mu \, \partial_x u) = 0 \\ 
\partial_t E   + \partial_x (E \, u  + p \, u ) - \partial_x ( \mu \, u \, \partial_x u)
- \, {{\gamma}\over{Pr}} \,  \partial_x  ( \mu \, \partial_x e ) = 0 
\end {array} \right. \monend 
The Fourier law of heat dissipation takes the form
$ \, Q =  - {{\gamma}\over{Pr}} \,  \partial_x ( \mu \, \partial_x e) \, $
and the viscous work is equal to  $ \,\, \partial_x ( \mu \, u \, \partial_x u) $.
In the isothermal case, there is only one thermodynamic variable, 
the third equation of (\ref{NS1d}) is not considered and the pressure is related to
density through the relation
\moneq \label{pression-isotherme} 
p = c_s^2 \, \rho \, . 
\monend 

%%%%%%%%%%%%%%%%%%%%%%%%%%%%%%%%%%%%%%%%%%%%%%%%%%%   navier stokes isotherme 
\smallskip \monitem 
In the bidimensional isothermal case, we have a vector
$ \, W \equiv (\rho ,\, j_x \equiv \rho\, u ,\, j_y \equiv \rho \, v )^{\rm t} \, $ 
of three conserved variables that define the two components $ \, u \, $ and $ \, v \, $
of the velocity.
We denote by $ \, V \equiv (\rho ,\,  u ,\,  v )^{\rm t} \, $
the primitive variables.
Introduce the divergence $ \, {\rm div}{\bf u} \equiv  \partial_x u +    \partial_y v \, $
of the velocity field, the shear viscosity $ \, \mu \, $
and the bulk viscosity $ \, \zeta $.
The symmetric viscous tensor $ \, \tau \, $ is defined according to
\moneq \label{tau-2d} 
\tau_{xx} = 2 \, \mu \, \partial_x u + ( \zeta -   \mu ) \,  {\textrm {div}} \, {\bf u}  \,,\,\,
\tau_{xy} = \mu \, (\partial_x v +  \partial_y u)  \,,\,\, 
\tau_{yy} = 2 \, \mu \, \partial_y v +( \zeta -   \mu ) \,  {\textrm {div}}\, {\bf u} \, . 
\monend
The conservation of mass and momentum   takes the form
\moneq \label{NS2d-isotherme}  
\partial_t W + \partial_x F^E_x(W) + \partial_y F^E_y(W) + 
\partial_x F^V_x(W,\, \nabla V) + \partial_y F^V_y(W,\, \nabla V) = 0 
\monend
with 
\moneqstar
\left\{ \begin {array}{l}
F^E_x(W) = ( \rho \, u \,,\, \rho\, u^2 + p \,,\, \rho \, u \, v )^{\rm t} \,,\,\, 
F^E_y(W) = ( \rho \, v \,,\,  \rho \, u \, v \,,\, \rho \, v^2 + p   )^{\rm t} \\
F^V_x(W) = -( 0  \,,\,  \tau_{xx} \,,\,  \tau_{xy} )^{\rm t} \,,\,\, 
F^V_y(W) = -( 0  \,,\,  \tau_{xy} \,,\,  \tau_{yy} )^{\rm t} \,. 
%%%   F^E_x(W) = \begin{pmatrix} \rho \, u \\ \rho\, u^2 + p \\ \rho \, u \, v \end{pmatrix} \,,\,\, 
%%%   F^E_y(W) = \begin{pmatrix} \rho \, v \\ \rho\, u \, v \\ \rho \, v^2 + p  \end{pmatrix} \,,\,\, 
%%%   F^V_x(W) = -\begin{pmatrix} 0  \\ \tau_{xx} \\ \tau_{xy} \end{pmatrix} \,,\,\, 
%%%   F^V_y(W) = -\begin{pmatrix} 0  \\ \tau_{xy} \\ \tau_{yy} \end{pmatrix} \,.
 \end{array} \right. \monendstar 
We can also write  the opposite of the divergence $ \,- \big(\partial_x F^V_x(W,\, \nabla V) + \partial_y F^V_y(W,\, \nabla V) \big)\,$
of the viscous fluxes under the form 
\moneq \label{flux-visqueux-2d} \left( \begin{array} {c}
  0 \\ \partial_j \tau_{xj} \\ \partial_j \tau_{yj}   \end{array} \right) \equiv
 \left( \begin{array} {c} 0 \\ 
 \partial_x ( 2 \, \mu \,   \partial_x u
+ (\zeta-\mu)  ( \partial_x u  +   \partial_y v ) )
+ \partial_y (\mu ( \partial_x v  +   \partial_y u ) )  \\
 \partial_x (\mu (   \partial_x v +    \partial_y u ) )
+ \partial_y (  (\zeta-\mu) (   \partial_x u  +   \partial_y v ) + 2 \, \mu \, \partial_y v )) \, . 
 \end{array} \right)  \monend

\smallskip \noindent
For three space dimensions, the isothermal Navier Stokes equations follow the same structure with
one more conservative variable:
$ \, W \equiv (\rho ,\, j_x \equiv \rho\, u ,\, j_y \equiv \rho \, v ,\, j_z \equiv \rho \, w )^{\rm t} $,
and one more primitive variable with
 $ \, V \equiv (\rho ,\,  u ,\,  v ,\,  w)^{\rm t} $. We have 
\moneq \label{NS3d-isotherme}  \left\{ \begin {array}{l}
\partial_t W + \partial_x F^E_x(W) + \partial_y F^E_y(W) +  \partial_z F^E_z(W) \\ 
\qquad + \, \partial_x F^V_x(W,\, \nabla V) + \partial_y F^V_y(W,\, \nabla V) + \partial_z F^V_z(W,\, \nabla V) = 0 
 \end{array} \right. \monend
with the fluxes defined according to 
\moneqstar
\left\{ \begin {array}{l}
F^E_x(W) = ( \rho \, u \,,\, \rho\, u^2 + p \,,\, \rho \, u \, v \,,\, \rho \, u \, w )^{\rm t}  \,,\,\, 
F^E_y(W) = ( \rho \, v \,,\,  \rho \, u \, v \,,\, \rho \, v^2 + p  \,,\, \rho \, v \, w   )^{\rm t} \\
F^E_z(W) = ( \rho \, w \,,\,  \rho \, u \, w \,,\,  \rho \, v \, w   \,,\, \rho \, w^2 + p  )^{\rm t} \,,\,\, 
F^V_x(W) = -( 0  \,,\,  \tau_{xx} \,,\,  \tau_{xy}  \,,\,  \tau_{xz})^{\rm t} \\ 
F^V_y(W) = -( 0  \,,\,  \tau_{xy} \,,\,  \tau_{yy}  \,,\,  \tau_{yz})^{\rm t} \,,\,\,   
F^V_z(W) = -( 0  \,,\,  \tau_{xz} \,,\,  \tau_{yz}  \,,\,  \tau_{zz})^{\rm t}  
\end{array} \right. \monendstar
and the viscous tensor $ \, \tau \, $ satisfying the relations 
\moneq \label{tau-3d}   \left \{ \begin {array}{l}
  \tau_{xx} = 2 \, \mu \, \partial_x u + \big( \zeta -  {2\over3} \, \mu \big) \,  {\textrm {div}} \, {\bf u}  \,,\,\,
  \tau_{yy} = 2 \, \mu \, \partial_y v + \big( \zeta -  {2\over3} \, \mu \big) \,  {\textrm {div}} \, {\bf u}  \\ 
  \tau_{zz} = 2 \, \mu \, \partial_z w + \big( \zeta -  {2\over3} \, \mu \big) \,  {\textrm {div}} \, {\bf u} \,,\,\,
  \tau_{xy} = \mu \, ( \partial_x v +  \partial_y u)  \\
  \tau_{yz} = \mu \, ( \partial_y w +  \partial_z v)  \,,\,\, \tau_{zx} = \mu \, ( \partial_z u +  \partial_x w ) \, 
\end {array} \right. \monend
with $ \, {\rm div}{\bf u} \equiv  \partial_x u + \partial_y v + \partial_z w $. 

%%%%%%%%%%%%%%%%%%%%%%%%%%%%%%%%%%%%%%%%%%%%%%%%%%%   navier stokes thermique 
\smallskip \monitem
When we add the conservation of energy $\, E \, $ defined by the relation (\ref{energie-totale})
with $ \, | {\bf u} | \, $ the modulus of velocity, all the thermodynamical description is analogous
to the one dimensional case.
The conservations of mass and momentum take the form (\ref{NS2d-isotherme}) or 
(\ref{NS3d-isotherme}) adapted to three  space dimensions. 
The shear viscosity $ \, \mu \, $ and the bulk  viscosity $ \, \zeta \, $
are now functions of two thermodynamical variables.

\smallskip \noindent
For two space dimensions, the conservation of energy is added to the relations  (\ref{NS2d-isotherme})
and it can be written
\moneq \label{conservation-energie-totale-2d}   \left \{ \begin {array}{l}
\partial_t E + \partial_x (u \, E + p \, u)  + \partial_y (v\, E  + p \, v) \\ 
\qquad - \partial_x( u\, \tau_{xx} + v \, \tau_{xy}) - \partial_y ( u\, \tau_{xy} + v \, \tau_{yy}) 
- {{\gamma}\over{Pr}} \, \big[  \partial_x ( \mu \, \partial_x e) + \partial_y ( \mu \, \partial_y e) \big] = 0 \, . 
\end {array} \right. \monend
Recall that we suppose a constant Prandtl number $ \, Pr $.
The opposite of the viscous fluxes for momentum and energy can be written as 
\moneq \label{flux-visqueus-ns-thermique-2d}
\Phi_{\rm NS}^{2D} = \left( \begin {array}{l}
  \partial_j \tau_{xj} \equiv \partial_x ( 2 \, \mu \,   \partial_x u
+ (\zeta-\mu)  ( \partial_x u  +   \partial_y v ) )
+ \partial_y (\mu ( \partial_x v  +   \partial_y u ) )  \\
\partial_j \tau_{yj}  \equiv \partial_x (\mu (   \partial_x v +    \partial_y u ) )
+ \partial_y (  (\zeta-\mu) (   \partial_x u  +   \partial_y v ) + 2 \, \mu \,   \partial_y v )) \\
%% 26 \, [
\partial_j ( u_i \, \sigma_{ij} ) + {{\gamma}\over{Pr}} \, \big( \partial_x ( \mu \, \partial_x e  ) + \partial_y ( \mu \,  \partial_y e  ) \big) %% ]
\end{array} \right) \, .  \monend

\smallskip \noindent 
For three  space dimensions, we have the same kind of relation for the
conservation of total energy, with more terms:
\moneq \label{conservation-energie-totale-3d}    \left \{ \begin {array}{l}
\partial_t E + \partial_x (E \, u + p \, u)  + \partial_y (E \, v + p \, v)  + \partial_z (E \, w + p \, w) \\
\quad - \partial_x( u\, \tau_{xx} + v \, \tau_{xy} + w \, \tau_{xz} )
- \partial_y( u\, \tau_{xy} + v \, \tau_{yy} + w \, \tau_{yz} )
- \partial_z( u\, \tau_{xz} + v \, \tau_{yz} + w \, \tau_{zz} ) \!\!\!\!   \\ 
\qquad - {{\gamma}\over{Pr}} \, \big[  \partial_x ( \mu \, \partial_x e) + \partial_y ( \mu \, \partial_y e)  
+ \partial_z ( \mu \, \partial_z e)\big] = 0 \, . 
\end {array} \right. \monend
And we have a relation of the type 
\moneq \label{flux-visqueus-ns-thermique-3d}
\Phi_{\rm NS}^{3D} \! = \!\! \left( \!\! \begin {array}{l}
  \partial_j \sigma_{xj} \equiv \partial_x \big( 2 \, \mu \,   \partial_x u \!+\! (\zeta-{2\over3}\, \mu) \,  {\rm div}\,{\bf u} \big) 
\!+\! \partial_y (\mu ( \partial_x v  \!+\!   \partial_y u ) ) \!+\! \partial_z (\mu ( \partial_x w  \!+\!   \partial_z u ) )  \\
\partial_j \sigma_{yj} \equiv \partial_x (\mu ( \partial_x v  \!+\!   \partial_y u ) ) 
\!+\! \partial_y  \big( 2 \, \mu \,   \partial_y v \!+\! (\zeta-{2\over3}\, \mu) \,  {\rm div}\,{\bf u} \big) 
 \!+\! \partial_z (\mu ( \partial_y w  \!+\!   \partial_z v ) )  \\
 \partial_j \sigma_{zj} \equiv \partial_x (\mu ( \partial_x w  \!+\!   \partial_z u ) )
 \!+\! \partial_y (\mu ( \partial_y w  \!+\!   \partial_z v ) )  
\!+\! \partial_z \big( 2 \, \mu \,   \partial_z w \!+\! (\zeta-{2\over3}\, \mu) \,  {\rm div}\,{\bf u} \big)  \\
%% 26 \, [
\partial_j ( u_i \, \sigma_{ij} ) + {{\gamma}\over{Pr}} \, \big( \partial_x ( \mu \, \partial_x e  ) + \partial_y ( \mu \,  \partial_y e  )
+ \partial_z ( \mu \,  \partial_z e  )\big) %% ]
\end{array} \!\! \right)  \monend
for the  opposite of the viscous fluxes for momentum and energy.

%%%%%%%%%%%%%%%%%%%%%%%%%%%%%%%%%%%%%%%%%%%%%%%%%%%%%%%%%%%%%%%%%%%%%%%%%%%%%%%  section 3
\bigskip \bigskip     \noindent {\bf \large  3) \quad  Lattice Boltzmann schemes  with multiple relaxation times}
%%%%%%%%%%%%%%%%%%%%%%%%%%%%%%%%%%%%%%%%%%%%%%%%%%%%%%%%%%%%%%%%%%%%%%%%%%%%%%%%%%%%%%%%%%

\noindent
In the  space $ \,  \RR^d \, $ of dimension $\, d$, we consider a finite set
of $ \, q \, $ discrete velocities $ \, v_j  \in  {\cal V} \, $
with components $\,  v_j^\alpha \, $ for $ \, 1 \leq \alpha \leq d $.
In most cases, the set ${\cal V}$ satisfies the symmetries of the square or cube, respectively for
2 or 3 dimensions of space.
The unknowns of the  lattice Boltzmann schemes are the particle densities $ \, f_j $. They are
functions of discrete space $\, x  $, discrete time $\, t \, $ and discrete velocities $\, v_j \, $:
%%%  with~$ \, 0 \leq j < q $:
%
\moneqstar
f_j = f_j(x,\, t), \quad x \in \Delta x \, \ZZ^d ,\, t = n \, \Delta t ,\, n \in \NN ,\,  0 \leq j < q \, . %%% v_j \in  {\cal V} \, .
\monendstar
The vector  $ \, f (x,\, t) \in \RR^q \, $ is constructed with the scalars  $\, f_j(x,\, t) \, $
for $ \,  0 \leq j < q $. When this vector is known at discrete time $ \, t $,  the lattice Boltzmann scheme
computes the distribution  $ \, f (x,\, t+\Delta t) \, $ at the new time step.

%%%%%%%%%%%%%%%%%%%%%%%%%%%%%%%%%%%%%%%%%%%%%%%%%%%%%%%%%%%%%    les moments et le schema 
\smallskip \monitem
Our framework concerns multi relaxation times:  %%% correction typo 04 et 06 juin 2020
%%%   with  d'Humi\`eres  \cite{DDH92},
we introduce a constant invertible matrix~$ \, M \, $
called ``d'Humi\`eres matrix''  \cite{DDH92} or moment matrix in this contribution.
This matrix defines the vector of moments $ \, m \in \RR^q \, $ by a simple product:
\moneq \label{moments}
 m_k \equiv \sum_{0 \leq j < q}  \, M_{k j} \,\, f_j \,.
 \monend
%
%%% {\color{red} 
Our method to construct the matrix $ \, M \, $ proposed in this work is composed by two steps.
First we choose a family of homogenesous polynomials to define the moments. Secondly, we use a Gram Schmidt
algorithm to obtain an orthogonal matrix. The first point can be a difficult task when the number
of velocities is increasing. For example with the D2W17 schemes, the set of velocities differ from
the previous D2Q17 scheme. The moment number~15 named~$ \, xx_{xy} \, $  has to be introduced instead
of the moment $ \, h_4 \, $ used previously that leads to a noninvertible matrix
(see the Tables \ref{d2q17-ns-polynomes} %% 21
and  \ref{d2w17-ns-polynomes} %%  23
at the end of the contribution).
The Gram-Schmidt algorithm we have chosen corresponds to a very simple 
scalar product in the space of velocities: $ \, (f,\, g) = \sum_\ell f_\ell \, g_\ell  $.

\noindent 
The precise polynomials used for the construction of the matrix $ \, M \, $
are  presented in the Appendix at the end of the contribution.
%%%  }
%
Nevertheless,
symmetry properties are a constant guide for this construction,
as pointed by Rubinstein and Luo \cite{RL08}. 

\noindent 
We introduce a new parameter: the number of conservation laws $ \, N \, $ ($  1 \leq N < q $).
The $ \, N \, $ first moments are ``conserved'' by the scheme.
The last $ \, q-N \, $  moments are not conserved; they are sometimes denominated as  ``microscopic moments''. 
Then it is natural to divide the vector of moments  into two families:
\moneq \label{moments-2-familles}
m \equiv \begin{pmatrix} W \\ Y \end  {pmatrix} \, .
\monend
The conserved moments or macroscopic  variables $ \, W \, $ constitute
a linear space of dimension~$ \, N $.
Observe that the nonconserved moments or microscopic  variables $ \, Y \, $
generate  a linear space of dimension $ \, q - N $.
%%%
%%% \smallskip \noindent
We introduce now equilibrium states $\, f^{\rm eq} $.  
They are  characterized in terms of conserved moments
with the help of a regular nonlinear vector field $ \, \Phi \,: \RR^N \longrightarrow \RR^{q-N}  \, $
such that
\moneq \label{f-equilibre}
f^{\rm eq} = M^{-1} \,  \begin{pmatrix} W \\ \Phi(W) \end  {pmatrix} \, .
\monend
In other words,  the vector field $ \, W \longmapsto Y^{\rm eq} \equiv \Phi(W) \, $ defines
the set of equilibrium states.
The numerical scheme is composed by a succession of two steps: a relaxation step and an advection step.
The relaxation step is local in space and modifies the vector $ \, f(x,\, t) \, $ into a new vector $ \, f^*(x,\, t) \, $
in the following way. We introduce a diagonal relaxation  matrix $ \, S \, $ of order $ \, (q-N) $:
\moneqstar 
S = {\rm diag} ( s_1 ,\, s_2 ,\, {\rm ...} ,\, s_{q-N} ) \, .
\monendstar
This matrix is proportional to the identity matrix  for the Bhatnagar-Gross-Krook \cite{BGK54}
variant of lattice Boltzmann schemes.
It can be an arbitrary strictly positive diagonal matrix in the
variant of lattice Boltzmann schemes with 
multiple relaxation times  \cite{DDH92}, provided $s_i < 2$ and invariance in the interchange of
coordinates.  
Then the  moments~$\, m^* \, $ after relaxation are defined with the relations 
\moneq \label{moments-relaxation}
m^* \equiv \begin{pmatrix} W^* \\ Y^*  \end  {pmatrix} =
\begin{pmatrix} W \\ Y + S \, (\Phi(W)-Y)  \end  {pmatrix} \, .
\monend
The particle distribution  $ \, f^*(x,\, t) \, $
after relaxation satisfies $ \, f^* =  M^{-1} \, m^* $. 
When this new distribution is computed, the advection step propagates the informations
to the neighbours:
\moneq \label{f-advection}
f_j(x,\, t + \Delta t) = f_j^*(x-v_j\, \Delta t ,\, t) \,,\quad x \in \Delta x \, \ZZ^d ,\, t = n \, \Delta t ,\, n \in \NN ,\,  0 \leq j < q \, .
\monend

%%%%%%%%%%%%%%%%%%%%%%%%%%%%%%%%%%%%%%%%%%%%%%%%%%%%%%%%%%%%%    analyse asymptotique  
\smallskip \monitem
We have developed in \cite{Du08,Du2022} a formal asymptotic analysis and the result is an explicitation of
equivalent partial differential equations satisfied by the conserved variables. %%  $ \, W $. 
This method is based on Taylor expansion and we call it also ``ABCD'' method. 
The hypotheses of the  formal expansion are precise. 
First, we adopt the acoustic scaling: the ratio $ \, \lambda \equiv {{\Delta x}\over{\Delta t}} \, $
is supposed to be constant in all this work.
Secondly the relaxation matrix $ \, S \, $ is fixed and invertible;
it is also the case for the matrix $ \, S^{-1} $. 

\smallskip \noindent
A  unintuitive result has been discovered by  H\'enon \cite{He87}.
In particular, the lattice Boltzmann scheme
put in evidence what we call the  H\'enon  matrix $ \, \Sigma \, $ in this contribution.
It is defined by
\moneq \label{Henon}
 \Sigma \equiv  S^{-1} - {1\over2} \, {\rm I} \, .
\monend
This matrix emerges from the very classic second order analysis.
For applications to fluid dynamics,
this matrix is closely related to viscosities. 
In general, some  values of the H\'enon matrix $ \, \Sigma \, $ are chosen as small as possible in order
to simulate flows with high Reynolds number.
In consequence, over-relaxation is a mandatory practice for lattice Boltzmann schemes
applied to  high Reynolds number flows.
Recall that this matrix remains fixed in this contribution.

%%%%%%%%%%%%%%%%%%%%%%%%%%%%%%%%%%%%%%%%%%%%%%%%%%%%%%%%    momentum-velocity operator
\smallskip \monitem
For a given lattice Boltzmann scheme of dimension $ \, d $,
we introduce the momentum-velocity operator matrix $ \, \Lambda \, $ defined by the relation
\moneq \label{Lambda}
 \Lambda  =  M \,\, {\rm diag} \,\Big(  \sum_{1  \leq \alpha \leq d}  v^\alpha\, \partial_\alpha \Big) \,\, M^{-1} \, .
\monend
It is a $ \, q \times q \, $ operator matrix
composed by first-order space differential operators.
It is obtained by conjugation of the first order advection operator $\, v . \nabla \, $
by the d'Humi\`eres matrix~$ \, M $.
The operator matrix $ \, \Lambda \, $ is nothing else than the advection operator seen in the basis of moments.
In particular the eigenvalues of this momentum-velocity operator are simply the advections  
$ \,  \sum_{1  \leq \alpha \leq d}  v^\alpha\, \partial_\alpha \, $ 
associated with the set of  discrete velocities of the lattice. 
When the invertible matrix $ \, M \, $ is changed,
the eigenvalues of this momentum-velocity operator do not change.

\smallskip \noindent 
We introduce a block decomposition  of the  momentum-velocity operator matrix associated
to the decomposition (\ref{moments-2-familles}) of the moments.  We define
a $ \, N \times N \, $ operator matrix $ \, A $,
a $ \, N \times (q-N) \, $ operator matrix $ \, B $,
a $ \,  (q-N) \times N \, $ operator matrix $ \, C \, $ and
a $ \,  (q-N) \times (q-N) \, $ operator matrix $ \, D \, $ according to
\moneq \label{bloc-Lambda}
\Lambda \equiv  \begin{pmatrix} A &   B \\ C  &  D  \end  {pmatrix} \, .
\monend
Remember that in the following, the matrices $ \, A $, $ \, B $, $ \, C \, $ and $ \, D \, $
are matrices composed with first order space  operators.
Then it is possible to explicit equivalent partial differential equations at an arbitrary
order. In \cite{Du2022} we have given the general relations up to fourth order accuracy.
In this contribution, the second order is sufficient. The equivalent partial differential equations
of the lattice Boltzmann scheme can be written
\moneq \label{edp-ordre-2}
\partial_t W + \Gamma_1(W) + \Delta t \, \Gamma_2(W) = {\rm O}(\Delta t^2) \, . 
\monend
The dynamic vectors $ \, \Gamma_1(W) \, $ and $ \, \Gamma_2(W) \, $ contain nonlinear terms
and partial differential operators of order 1 and 2 respectively.
Their compact  expression is simple with the  momentum-velocity operator matrix: 
\moneq \label{Gamma-1-et-Gamma-2}
\Gamma_1(W) = A \, W + B \, \Phi(W) \,,\,\, 
\Gamma_2(W) = B \, \Sigma \, \Psi_1 
\monend
with 
\moneq \label{Psi_1} 
\Psi_1 = \dd \Phi(W) .  \Gamma_1  - (C \, W + D \, \Phi(W))  \, .
\monend 
Starting from the lattice Boltzmann scheme, it is a good exercice to explicit the
two dynamic vectors of (\ref{Gamma-1-et-Gamma-2}). Observe that this algorithm is operational in all
generality in the pylbm software \cite{GG17}.

\noindent 
%%%%  {\color{red}
It is possible to transcribe all the algebraic relations for the non conserved moments
in terms of the equilibrium particle distribution. 
We have just to apply the relation (\ref{f-equilibre}). 
%% $ \, f^{eq} = M^{-1} \, \begin{pmatrix} W \\  \Phi ( W ) \end{pmatrix}  $.
Moreover,  changing the transfer matrix $ \, M \, $ is changing also the relaxation process
and in consequence all the lattice Boltzmann scheme.
We have to keep in mind that our results concern only the set of matrices~$ \, M \, $ we have
explicitly used. They have to be revisited if this  matrix is replaced by an other one. 
%%%  }

\noindent
The question is now to know if it is possible or not
to identify the first order vector $ \, \Gamma_1\, $ with the Euler equations
of gas dynamics and the  second order vector $ \, \Gamma_2 \, $
with the dissipative terms of the Navier Stokes equations.
This question for various physical models and various lattice Boltzmann schemes.

%%%%%%%%%%%%%%%%%%%%%%%%%%%%%%%%%%%%%%%%%%%%%%%%%%%%%%%%%%%%%%%%%%%%%%%%%%%%%%%  section 4 
\bigskip \bigskip    \noindent {\bf \large    4) \quad  Two-dimensional isothermal Navier Stokes}
%%%%%%%%%%%%%%%%%%%%%%%%%%%%%%%%%%%%%%%%%%%%%%%%%%%%%%%%%%%%%%%%%%%%%%%%%%%%%%%%%%%%%%%%%%%%%%%%%%%%%%

%%%%%%%%%%%%%%%%%%%%%%%%%%%%%%%%%%       31 decembre 2018       %%%%%%%%%%%%%%%%%%%%%%%%%
\fancyhead[EC]{\sc{Fran\c{c}ois Dubois and Pierre Lallemand}}
\fancyhead[OC]{\sc{Single lattice Boltzmann distribution for Navier Stokes equations}}
%%%%%%%%%%%%%%%%%%%%%%%%%%%%%%%%%%       31 decembre 2018       %%%%%%%%%%%%%%%%%%%%%%%%%
%%%%%%%%%%%%%%%%%%%%%%%%%%%%%%%%%%%%%%%%%%%%%  jolie numerotation des pages
\fancyfoot[C]{\oldstylenums{\thepage}}
%%%%%%%%%%%%%%%%%%%%%%%%%%%%%%%%%%%%%%%%%%%%%  fin jolie numerotation des pages

\noindent
We recall our progressive methodology to fit the parameters of D2Q9 and D2Q13 schemes
in order to approximate isothermal Navier Stokes equations of fluid dynamics for two space dimensions.

%%%%%%%%%%%%%%%%%%%%%%%%%%%%%%%%%%%%%%%%%%%%%%%%%%%%%%%%%%%%%%%%%%%%%%%%%%%%%
\smallskip \monitem  D2Q9
%%%%%%%%%%%%%%%%%%%%%%%%%%%%%%%%%%%%%%%%%%%%%%%%%%%%%%%%%%%%%%%%%%%%%%%%%%%%%

\begin{figure}    [H]  \centering
\vskip -1.3 cm 
%% \centerline {\includegraphics[width=.350\textwidth]{lbm-graphes-d2q9-fleches-01mai2020.pdf}}
\centerline {\includegraphics[width=.350\textwidth]{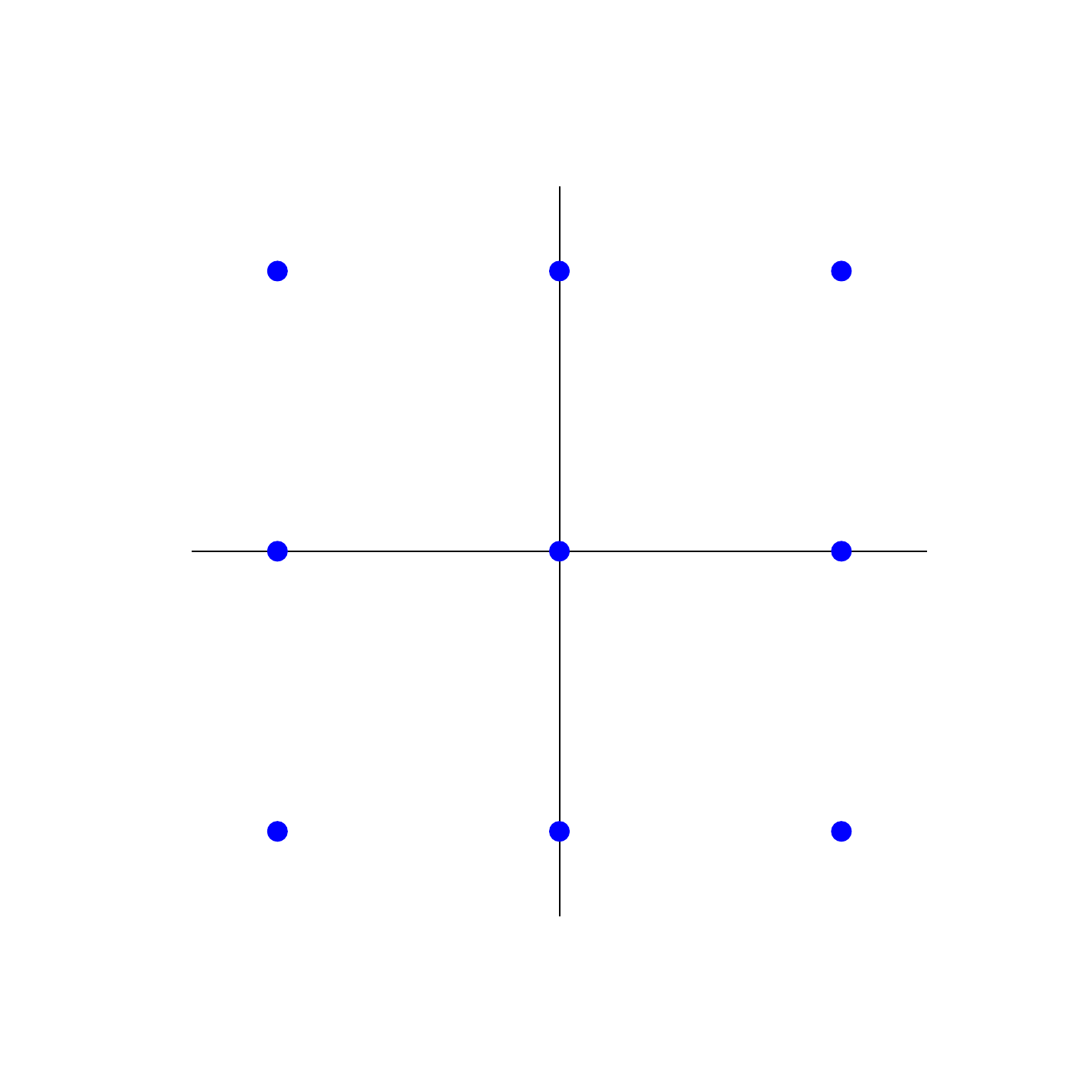}}
\vskip -.7 cm 
\caption{D2Q9 lattice Boltzmann scheme }
\label{fig-d2q9-bis} \end{figure}

\noindent
The set of velocities of the D2Q9 lattice Boltzmann scheme are recalled in Figure \ref{fig-d2q9-bis}. 

%%%%%%%%%%%%%%%%%%%%%%%%%%%%%%%%%%%%%%%%%%%%%%%%%%%%%%%%%%%%%%%%%%%%%%%%%%%%%%%%%%% figure  d2q9
%%%%%%%%%%%%%%%%%%%%%%%%%%%%%%%%%%%%%%%%%%%%%%%%%%%%%%%%%%%%%%%%%%%%%%%%%%%%%%%%%%%

%%%   \hfill $ \lambda = {{\Delta x}\over{\Delta t}} $ \qquad ~

The moments $ \, m \, $ for this classical D2Q9 scheme are given by the relation
\moneqstar
m = M \, f \, . 
\monendstar
Following {\it e.g.} \cite{LL00}, the d'Humi\`eres matrix of moments $ \, M \, $ is given by the relation  (\ref{d2q9-matrice-M}). 
\moneq  \label{d2q9-matrice-M}
M =  \left[  \!\!   \begin{array} {ccccccccc}
 1  \!&\!  1  \!&\!  1  \!&\!  1  \!&\!  1  \!&\!  1  \!&\!  1  \!&\!  1  \!&\!  1 \\
   0 \!&\! \lambda  \!&\!  0  \!&\!  -\lambda  \!&\!  0  \!&\!   \lambda  \!&\!   -\lambda   \!&\!  -\lambda  \!&\!  \lambda  \\
 0  \!&\!  0   \!&\! \lambda  \!&\!  0    \!&\!   -\lambda   \!&\!   \lambda   \!&\!   \lambda  \!&\!  -\lambda  \!&\!  -\lambda \\
-4 \lambda^2    \!&\!  -\lambda^2   \!&\!   -\lambda^2   \!&\!   -\lambda^2   \!&\!   -\lambda^2  \!&\!   2 \lambda^2   \!&\!  2 \lambda^2   \!&\!  2 \lambda^2   \!&\!  2 \lambda^2 \\
  0    \!&\!   \lambda^2    \!&\!  -\lambda^2    \!&\!   \lambda^2   \!&\!   -\lambda^2  \!&\!  0  \!&\!  0  \!&\!  0  \!&\!  0 \\
 0  \!&\!  0  \!&\!  0  \!&\!  0  \!&\!  0   \!&\!   \lambda^2   \!&\!   -\lambda^2    \!&\!   \lambda^2   \!&\!   -\lambda^2 \\
 0   \!&\! -2 \lambda^3  \!&\!  0   \!&\!  2 \lambda^3  \!&\!  0   \!&\!  \lambda^3   \!&\!   -\lambda^3   \!&\!   -\lambda^3   \!&\!    \lambda^3 \\
 0  \!&\!  0  \!&\! -2 \lambda^3  \!&\!  0  \!&\!  2 \lambda^3    \!&\!   \lambda^3    \!&\!   \lambda^3    \!&\!  -\lambda^3    \!&\!  -\lambda^3 \\
4 \lambda^4  \!&\! -2 \lambda^4  \!&\! -2 \lambda^4  \!&\! -2 \lambda^4  \!&\! -2 \lambda^4  \!&\!   \lambda^4   \!&\!  \lambda^4   \!&\!  \lambda^4  \!&\!   \lambda^4
  \end{array} \!\! \right] \, . 
\monend
The lines of this invertible matrix are chosen orthogonal, 
and it is the case for all the schemes we consider in this contribution:
\moneqstar
\sum_j M_{ij} \, M_{kj} = 0 \,\,\, {\rm if} \,\, i \not = k . 
\monendstar 
We have also 
\moneq \label{moments-polynomes}
M_{ij} = p_i(v_j) \,,\,\, 0 \leq i,\, j < q  
\monend
with $ \, q = 9 \, $ for this scheme. 
The family of polynomials $ \, p_i \, $ are given in the Appendix at the Table~\ref{polynomes-d2q9}.
The moments are named with the notation
$ \, \rho  \,,\,\,  j_x  \,,\,\,  j_y \,,\,\, \varepsilon  \,,\,\,   xx \,,\,\, xy  \,,\,\,
q_x  \,,\,\,   q_y  \,,\,   h $.
The first moment $ \, \rho \, $ is a polynomial of degree 0, the moments $ \, j_x \, $ and $ \, j_y \, $
are defined with polynomials of degree~1, the 3 moments of degree~2 are related to energy ($\varepsilon$)
and to higher moments ($xx \, $ and $ \, xy $).
Observe that the moments $ \, xx \simeq v_x^2 - v_y^2 \, $ and $ \, xy \simeq v_x \, v_y \, $
correspond to  tensors with a trace equal to zero. 
The moments $ \, q_x \, $ and $ \, q_y \, $ correspond to polynomials of degree 3; they have some relation
with heat flux \cite{LL00}. Finally, the ``square of energy'' $ \, h \, $ is associated to a polynomial of degree~4. 
For reasons that appear in the following, we divide this set of moments  into four families:
the conserved variables $ \, \rho  \,,\,\,  j_x  \,,\,\,  j_y $,
the moments   $ \,  \varepsilon  \,,\,\,   xx \,,\,\, xy \, $ associated to polynomials of degree~2
that allow the fitting  the partial differential equations at first order,
the moments  $ \, q_x  \,,\,\,   q_y  \, $ for  the control of second order
viscous terms of the equivalent  partial differential equations and 
the last family  composed here by the unique moment $ \, h \, $ that does not appear in second order
equivalent equations.   These four families are presented at Table \ref{d2q9-moments-isotherme}. 
%
%%%%%%%%%%%%%%%%%%%%%%%%%%%%%%%%%%%%%%%%%%%%%%%%%%%%%%%%%%%%%%%%%%%%%%%%%%%%%%%%%%% moments d2q9    
\begin{table}    [H]  \centering
\begin{tabular}{|c|l|c|c|} \hline
 & conserved &  $\rho \,,\,\, j_x \,,\,\, j_y $  & 3  \\  \hline
1 & fit the Euler equations  &  $ \varepsilon \,,\,\,  xx  \,,\,\,  xy    $ & 3 \\ \hline 
%2 & fit the viscous terms ? & $ q_x  \,,\,\, q_y  $  & 2  \\ \hline  
2 & fit the viscous terms  & $ q_x  \,,\,\, q_y  $  & 2  \\ \hline  
3 & without influence & $ h  $  & 1   \\ \hline 
\end{tabular} 
\caption{The four families of moments for the D2Q9  scheme  for the approximation of the isothermal Navier Stokes equations}
\label{d2q9-moments-isotherme} \end{table}
%%%%%%%%%%%%%%%%%%%%%%%%%%%%%%%%%%%%%%%%%%%%%%%%%%%%%%%%%%%%%%%%%%%%%%%%%%%%%%%%%%%
%
The momentum-velocity operator matrix $ \, \Lambda \, $ is an operator matrix defined by the relation (\ref{Lambda}). 
%
%%%%%%%% \moneqstar %%   \label{general-matrice-Lambda}
%%%%%%%%  \Lambda  \equiv  M \, {\rm diag} \big(  \sum_{1 \leq \alpha\leq 2}  v^\alpha\, \partial_\alpha \big) \, M^{-1} \, . 
%%%%%%%% \monendstar
%
It can be evaluated without difficulty for isothermal Navier Stokes with 3 conservations
and we get (see {\it e.g.} \cite{Du2022}).
\moneq  \label{d2q9-matrice-Lambda}
\Lambda_{D2Q9}^{\rm iso} =  \left( \,\, \begin{array}{|ccc|cccccc|} \hline
 0 & \!\! \partial_x & \!\! \partial_y &   0  &   0 & \!\!   0 &
\!\!   0  &    0  &    0  \\
{{2 \lambda^2}\over{3}} \, \partial_x & \!\!  0 & \!\!  0 &   {1\over6} \, \partial_x &
  {1\over2} \, \partial_x  &     \!\!  \partial_y  & \!\!   0  &   0  &   0 \\
{{2\lambda^2}\over{3}} \, \partial_y & \!\!  0 & \!\!  0 &   {1\over6} \, \partial_y &
   -{1\over2} \, \partial_y  &    \!\!  \partial_x  & \!\!   0  &
  0  &   0  \\  \vspace{-.4 cm} & & & & & & & & \\ \hline  \vspace{-.4 cm} & & & & & & & &  \\
 0  &  \!\!  \lambda^2 \, \partial_x  & \!\!   \lambda^2 \, \partial_y  &
   0  &    0  & \!\!    0  & \!\!    \partial_x  &     \partial_y  &    0  \\
  0  &  \!\!   {{\lambda^2}\over{3}} \, \partial_x  & \!\!   -{{\lambda^2}\over{3}} \, \partial_y  &   0 &   0  & \!\!
  0  & \!\!   - {1\over3} \, \partial_x  &    {1\over3} \, \partial_y   &   0  \\
 0  &  \!\!   {{2 \, \lambda^2}\over{3}} \, \partial_y  &  \!    %% ici un facteur 2 (20 fevrier 2021)
 {{2 \,\lambda^2}\over{3}} \, \partial_x  & \!\!     0  &    0  & \!\!  %% et encore ici un facteur 2 (20 fevrier 2021)
   0  & \!\!     {1\over3} \, \partial_y  &  %% faute de signe etait ici (20 fevrier 2021)
 \!\!    {1\over3} \, \partial_x  &    0   \\   0  & \!\!    0  & \!\!    0  &
     {{\lambda^2}\over{3}} \, \partial_x  &   - \lambda^2 \, \partial_x  &
\!\!    \lambda^2 \, \partial_y  & \!\!    0  &    0  &    {1\over3} \, \partial_x   \\
  0  & \!\!    0  & \!\!    0  &     {{\lambda^2}\over{3}} \, \partial_y  &     \lambda^2 \, \partial_y  &
   \lambda^2 \, \partial_x  & \!\!   0  &      0  & \!\!     {1\over3} \, \partial_y  \\
  0  & \!\!    0  & \!\!   0  &   0  &   0  &   0  &  \lambda^2 \, \partial_x  &   \lambda^2 \, \partial_y  &
  0   \\ \hline
  \end{array}  \,\,   \right) \, . 
\monend
In the relation (\ref{d2q9-matrice-Lambda}), we have emphasized the ``ABCD'' block decomposition: 
\moneqstar 
\Lambda =   \begin{pmatrix} A &   B \\   C  &
   D  \end  {pmatrix} \,,
\monendstar
with $ \, A \, $and $ \, D \, $ square matrices and $ \, B \, $ and $ \, C \, $ rectangular ones. 

%%%%%%%%%%%%%%%%%%%%%%%%%%%%%%%%%%%%%%%%%%%%%%%%%%%%%%%%%%%%%%%%%%%%%%%%%%% d2q9 isotherme  ordre 1
\smallskip \monitem 
At first order, we have  $ \, \Gamma_1 = A \, W + B \, \Phi(W) $.
After some lines of algebra, we obtain 
\moneqstar %%    \label{orsay-equili-moments-second-ordre}
\Gamma_1 =  \left(  \begin{array} {l}
\partial_x j_x + \partial_y j_y \\   \vspace{-.4 cm} \\
{2\over3} \, \lambda^2 \,  \partial_x \rho + {1\over6} \,  \partial_x  \Phi_{\varepsilon} +  {1\over2} \,  \partial_x  \Phi_{xx}
+  \partial_y  \Phi_{xy} \\  \vspace{-.4 cm} \\
{2\over3} \, \lambda^2 \,  \partial_y \rho +{1\over6} \,  \partial_y  \Phi_{\varepsilon}  -  {1\over2} \,  \partial_y  \Phi_{xx}
+  \partial_x  \Phi_{xy} 
\end{array} \right) \, .  \monendstar

\smallskip \noindent 
First, we wish to recover the first order terms of the Navier Stokes equations, {\it id est} the Euler equations of gas dynamics
obtained from (\ref{NS2d-isotherme}) with $ \, \mu = \zeta = 0 $:
\moneq \label{euler-isotherme-2d} 
\left\{  \begin{array} {l}
\partial_t  \rho +  \partial_x j_x + \partial_y j_y  = 0  \\
 \partial_t  j_x \,+\,  \partial_x \big( {{1}\over{\rho}} j_x^2 + p \big)
 \,+ \,\partial_y \big(   {{1}\over{\rho}} j_x \, j_y  \big)  = 0  \\
\partial_t  j_y  \,+ \,\partial_x \big(  {{1}\over{\rho}} j_x \, j_y \big)
 \,+\,  \partial_y \big(  {{1}\over{\rho}} j_y^2  + p \big) = 0 \, . 
\end{array} \right.  \monend
We identify the various expressions inside the space partial derivatives,
with $ \, j_x \equiv \rho \, u \, $ and $ \, j_y \equiv  \rho \, v$: 
\moneqstar \left\{  \begin{array} {l}
{2\over3} \, \lambda^2 \,  \rho +  {1\over6} \,  \Phi_{\varepsilon} +  {1\over2} \, \Phi_{xx} = \rho \, u^2 + p \\
{2\over3} \, \lambda^2 \,  \rho +  {1\over6} \,  \Phi_{\varepsilon} -  {1\over2} \, \Phi_{xx} = \rho \, v^2 + p \\ 
 \Phi_{xy}  = \rho \, u \, v \,  . 
\end{array} \right.  \monendstar
Then 
\moneqstar
\left\{  \begin{array} {l}
\Phi_{\varepsilon} = 6 \, p  -  4 \, \lambda^2 \, \rho \, + \, 3 \, \rho \, ( u^2 + v^2 )    \\  
\Phi_{xx}   =    \rho \, ( u^2 - v^2 )   \,,\,\,  \Phi_{xy}   =   \rho \, u \, v  
\end{array} \right. \monendstar
and the equilibrium value  of the second family  of moments is explicited.

%%%%%%%%%%%%%%%%%%%%%%%%%%%%%%%%%%%%%%%%%%%%%%%%%%%%%%%%%%%%%%%%%%%%%%%%%%% d2q9 isotherme second ordre 
\smallskip \monitem
Observe that the general relaxation process (\ref{moments-relaxation}) takes the
 following form for the microscopic moments of the D2Q9 scheme 
\moneq   \label{s-ordre-deux}
\left \{ \begin {array}{l} 
\varepsilon^* = \varepsilon + s_e \, (\Phi_\varepsilon - \varepsilon )  \,,\,\,      xx^* = xx + s_x \, (\Phi_{xx} - xx)
\,,\,\,   xy^* = xy + s_x \, (\Phi_{xy} - xy)    \\
 q_x^* = q_x + s_q \, (\Phi_{qx} - q_x)   \,,\,\,   q_y^* = q_y + s_q \, (\Phi_{qy} - q_y)
 \,,\,\,   h^* = h + s_h \, (\Phi_{h} - h)  \, , 
\end {array} \right. \monend
with given parameters $\, s_e $, $\, s_x $, $\, s_q \, $ and $ \, s_h $. 
The  partial set of  H\'enon parameters 
\moneq \label{michel-henon}  
\sigma_x  =  {{1}\over{s_x}} - {1\over2}  \,,\,\, \sigma_e  =  {{1}\over{s_e}} - {1\over2}  \,,
\monend
are in evidence in the partial differential equations at second order. 
They allow the construction of the H\'enon matrix $ \, \Sigma \, $ defined in all generality in (\ref{Henon})
and explicited for the D2Q9 scheme as 
\moneq \label{d2q9-henon} 
\Sigma = {\rm diag} \, ( \sigma_e \,,\, \sigma_x \,,\, \sigma_x \,,\, \sigma_q \,,\, \sigma_q \,,\, \sigma_h ) \,. 
\monend
For the determination of the second order terms, we first construct the vector $ \, \Psi_1 \, $
introduced in (\ref{Psi_1}). 
%
%%%%%%%  \moneq \label{Psi_1} 
%%%%%%%  \Psi_1 = \dd \Phi(W) .  \Gamma_1  - (C \, W + D \, \Phi(W))  \, .
%%%%%%%  \monend 
%
Then the viscous fluxes satisfy 
\moneqstar
- \Delta t \, \Gamma_2 =  - \Delta t \,  B \, \Sigma \, \Psi_1 \, . 
\monendstar
If  the  isothermal Navier Stokes equations are satisfied, these expressions must be equal to 
the physical fluxes
\moneqstar %%   \label{flux-visqueux-2d}
\left( \begin{array} {c}
  0 \\ \partial_j \tau_{xj} \\ \partial_j \tau_{yj}   \end{array} \right) \equiv
 \left( \begin{array} {c} 0 \\ 
 \partial_x ( 2 \, \mu \,   \partial_x u
+ (\zeta-\mu)  ( \partial_x u  +   \partial_y v ) )
+ \partial_y (\mu ( \partial_x v  +   \partial_y u ) )  \\
 \partial_x (\mu (   \partial_x v +    \partial_y u ) )
+ \partial_y (  (\zeta-\mu) (   \partial_x u  +   \partial_y v ) + 2 \, \mu \, \partial_y v ))
 \end{array} \right)
 \monendstar 
introduced in (\ref{flux-visqueux-2d}). 
We must now solve a linear system
with unknowns  equal to the partial derivatives relative to $ \, \rho $, $ \, u \, $ and $ \, v \, $
of the equilibrium functions $ \, \Phi_{qx} \, $ and $ \, \Phi_{qy} \, $ 
for the second family $\, q_x ,\, q_y \, $ of nonconserved moments.
This system is composed by one equation for each of the  2 moments,
relative to each  dimension,
 one equation for each of  the associated partial derivatives $ \, \partial_x \, $ and $ \, \partial_y $,
%%one equation for each of  the 3 nonconserved variales $ \, \rho $, $ \, u $, $ \, v \, $
one equation for each of  the 3 conserved variables $ \, \rho $, $ \, u $, $ \, v \, $
and one equation for each of  the 2 space partial derivatives of these variables, then a total of 
$ \, 2 \times 2 \times 2 \times 3  = 24 \, $ equations.
Observe that we have only $ \,  3 \times 2 =  6  $ unknowns
since $ \, \Phi_{qx} \, $ and $ \, \Phi_{qy} \, $  must be explicited and $ \, \Phi_{h} \, $
has no influence.
If we try to  avoid  unphysical terms like $ \, \partial_x \rho \, $ and  $ \, \partial_y \rho \, $
from the second order fluxes, we do not find any solution.
Nevertheless, when we enforce 
\moneqstar %%   \label{d2q9-pression-isotherme}
p(\rho) = {{\lambda^2}\over{3}} \, \rho
\monendstar
{\it id est} if we choose  $ \, c_s = {{\lambda}\over{\sqrt{3}}} \, $
and if we impose  the relations 
\moneq  \label{d2q9-moments-equilibre-3e-ns-isotherme} %% \left\{ \begin{array} {l}
\Phi_{qx} = - \rho \, u \, \lambda^2  + 3 \, \rho \, (u^2 + v^2)  \, u\,,\,\, %%% coquille trouvee par Mahdi en juillet 2024 
\Phi_{qy}  = - \rho \,  v \, \lambda^2  + 3 \, \rho \, (u^2 + v^2) \,  v %%% coquille trouvee par Mahdi en juillet 2024 
\monend
suggested in \cite{LL00}, we have a third order  relative to velocity.
More precisely, with  the  H\'enon matrix given by the relation (\ref{d2q9-henon}), 
we can introduce the shear viscosity $ \,  \mu \, $ and the bulk viscosity  $ \,  \zeta $. 
We have 
\moneqstar
\mu    =  {{\lambda}\over{3}} \,  \rho \, \sigma_x \, \Delta x  \,,\,\, 
\zeta    = {{\lambda}\over{3}} \,  \rho  \, \sigma_e \, \Delta x
\monendstar
and 
\moneqstar \left\{ \begin{array} {l}
- \Delta t \, \Gamma_2 = \begin {pmatrix} 0 \\ \partial_j \tau_{xj}  \\ \partial_j \tau_{yj}  \end{pmatrix}
- \, \sigma_x \, \Delta t \, \partial_x  \! \begin {pmatrix} 0 \\ u^3 \, \partial_x \rho - v^3 \, \partial_y \rho
+ 3 \, \rho \, (u^2 \, \partial_x u  - v^2 \, \partial_y v) \\
-v^3 \, \partial_x \rho - u^3 \, \partial_y \rho -3 \, \rho \, ( u^2 \, \partial_y u + v^2 \, \partial_x v )
\end{pmatrix} \\ \qquad \qquad  \qquad 
- \, \sigma_x \, \Delta t \, \partial_y   \!  \begin {pmatrix} 0 \\ -v^3 \, \partial_x \rho - u^3 \,  \partial_y \rho
-3 \, \rho \, ( u^2 \, \partial_y u + v^2 \, \partial_x v )
\\ -u^3 \,  \partial_x \rho + v^3 \,  \partial_y \rho + 3 \, \rho \, (-u^2 \, \partial_x u + v^2 \,  \partial_y v ) 
\end{pmatrix}  \, . \end{array} \right. \monendstar  
The classical ``third order'' discrepancy of the D2Q9 scheme for isothermal flows \cite{LL00}  is completely explicited.

%%% \newpage 
%%%%%%%%%%%%%%%%%%%%%%%%%%%%%%%%%%%%%%%%%%%%%%%%%%%%%%%%%%%%%%%%%%%%%%%%%%%%%%%
\smallskip \monitem  D2Q13
%%%%%%%%%%%%%%%%%%%%%%%%%%%%%%%%%%%%%%%%%%%%%%%%%%%%%%%%%%%%%%%%%%%%%%%%%%%%%%%

\begin{figure}    [H]  \centering
\vskip -1.5 cm 
%%%   \centerline { \includegraphics[width=.44\textwidth]{d2q13-18mars2021-fig.pdf}}
\centerline { \includegraphics[width=.44\textwidth]{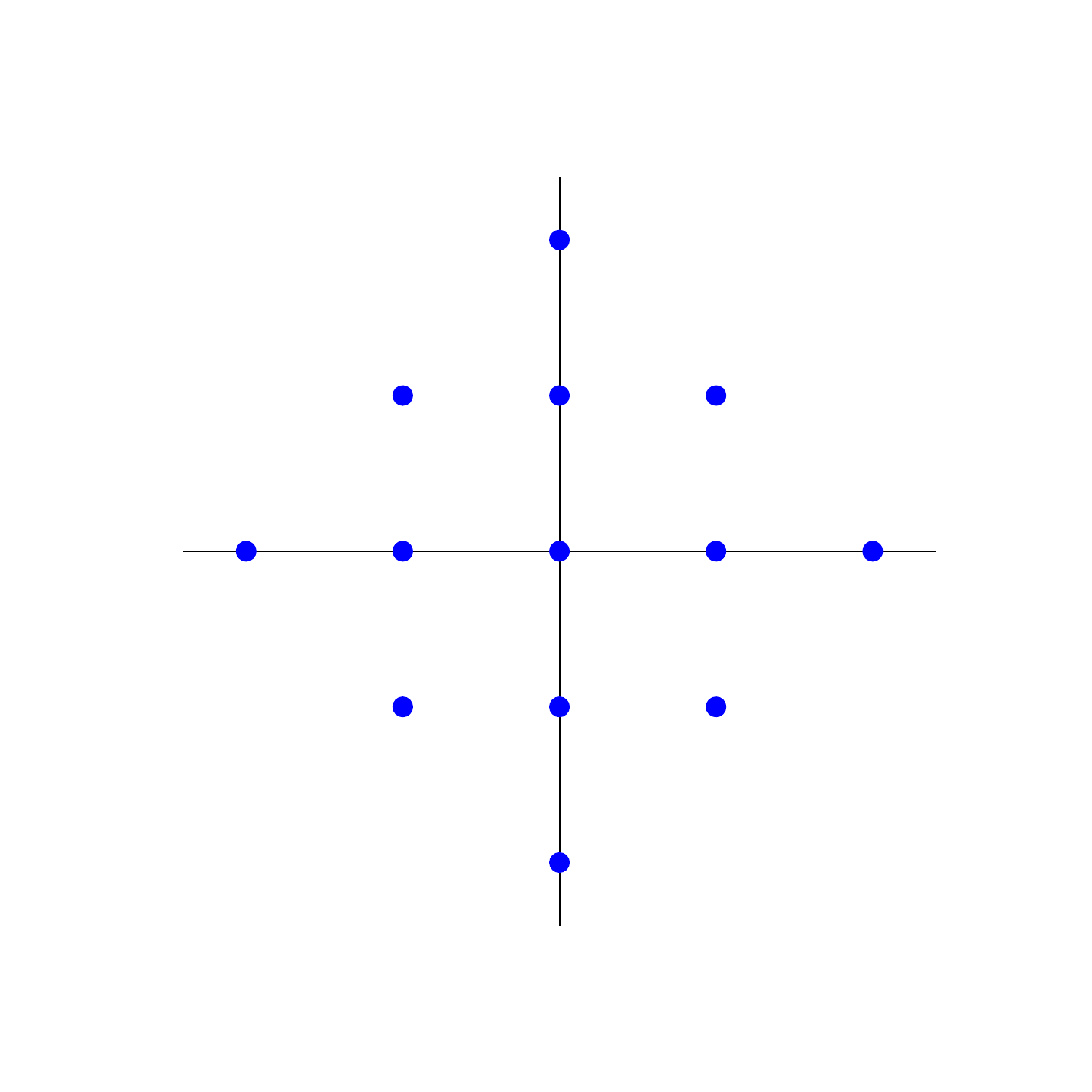}}
\vskip -1 cm 
\caption{D2Q13 lattice Boltzmann scheme }
\label{fig-d2q13} \end{figure}

\noindent
The D2Q13 scheme adds four velocities of the type $ \, (2,0) \, $ to the D2Q9 stencil,
as presented in Figure \ref{fig-d2q13}. %% \refas proposed in  \cite{Qi93}.  

%%%%%%%%%%%%%%%%%%%%%%%%%%%%%%%%%%%%%%%%%%%%%%%%%%%%%%%%%%%%%%%%%%%%%%%%%%%%%%%%%%% figure  d2q13
%%%%%%%%%%%%%%%%%%%%%%%%%%%%%%%%%%%%%%%%%%%%%%%%%%%%%%%%%%%%%%%%%%%%%%%%%%%%%%%%%%%
%
The moments are obtained through polynomials as in relation (\ref{moments-polynomes}) and the D2Q13 polynomials
for isothermal flows are defined in Table \ref{polynomes-d2q13}.
We have 3 conserved moments and the other 10 nonconserved ones can be ordered into three families as described
in Table \ref{d2q13-moments-isotherme}. 
%
%%%%%%%%%%%%%%%%%%%%%%%%%%%%%%%%%%%%%%%%%%%%%%%%%%%%%%%%%%%%%%%%%%%%%%%%%%%%%%%%%%% moments d2q13   
\begin{table}    [H]  \centering
\begin{tabular}{|c|l|c|c|} \hline
 & conserved &  $\rho \,,\,\, j_x \,,\,\, j_y $ & 3  \\  \hline
1 & fit the Euler equations  &  $ \varepsilon \,,\,\,  xx  \,,\,\,  xy  $ & 3 \\ \hline 
2 & fit the viscous terms & $ q_x  \,,\,\, q_y  \,,\,\, r_x \,,\,\,  r_y  $  & 4  \\ \hline  
3 & without influence & $ h  \,,\,\, xx_e  \,,\,\,   h_3   $  & 3   \\ \hline 
\end{tabular} 
\caption{The four families of moments for the D2Q13  scheme  for the approximation of the isothermal Navier Stokes equations}
\label{d2q13-moments-isotherme} \end{table}
%%%%%%%%%%%%%%%%%%%%%%%%%%%%%%%%%%%%%%%%%%%%%%%%%%%%%%%%%%%%%%%%%%%%%%%%%%%%%%%%%%%
% 
The momentum-velocity operator matrix $ \, \Lambda \, $ for the D2Q13 scheme is presented at relation (\ref{d213-Lambda-isotherme}). 
The numerical coefficients in front of the partial derivatives are not explicited in this document for practical
reasons  of space, 
but they are explicited in \cite{Du22git}.  
\moneq   \label{d213-Lambda-isotherme}
\Lambda_{D2Q13}^{\rm iso}  =  \left[ \begin{array} {ccccccccccccc}
0 \!\!&\!\!  * \partial_x \!\!&\!\!  * \partial_y  & %% debut B
0  \!\! & \!\! 0  \!\!&\!\! 0 \!\!&\!\!0\!\!&\!\! 0  \!\!&\!\! 0  \!\!&\!\! 0 \!\!&\!\! 0 \!\!&\!\! 0 \!\!&\!\! 0 \\
* \partial_x  \!\!&\!\! 0 \!\!&\!\! 0 &   %% debut B
 * \partial_x     \!\! &  \!\! * \partial_x  \!\!&\!\! * \partial_y  \!\!&\!\! 0 \!\!&\!\!
0 \!\!&\!\! 0 \!\!&\!\! 0 \!\!&\!\! 0 \!\!&\!\! 0 \!\!&\!\! 0 \\
* \partial_y  \!\!&\!\! 0 \!\!&\!\! 0 &  %% debut B
 * \partial_y   \!\! &  \!\! * \partial_y \!\!&\!\! * \partial_x  \!\!&\!\! 0 \!\!&\!\! 0 \!\!&\!\! 0 \!\!&\!\! 0 \!\!&\!\!
0 \!\!&\!\! 0 \!\!&\!\! 0 \\ \vspace{-4 mm} \\  %%% fin des matrices A et B
0  \!\!&\!\! * \partial_x   \!\!&\!\!  * \partial_y  & %% debut D
0  \!\!  & \!\!  0 \!\!&\!\! 0 \!\!&\!\! * \partial_x  \!\!&\!\!  * \partial_y  \!\!&\!\! 0 \!\!&\!\! 0 \!\!&\!\! 0 \!\!&\!\!
0 \!\!&\!\! 0 \\
0 \!\!&\!\! * \partial_x  \!\!&\!\! * \partial_y   & %% debut D
 0  \!\!& \!\! 0 \!\!& \!\! 0 \!\!&\!\!  * \partial_x \!\!&\!\!
 * \partial_y \!\!&\!\!  * \partial_x \!\!&\!\!
 * \partial_y \!\!&\!\! 0 \!\!&\!\! 0 \!\!&\!\! 0 \\
0 \!\!&\!\! * \partial_y  \!\!&\!\!  * \partial_x  & %% debut D
0   \!\!&\!\! 0 \!\!&\!\! 0 \!\!&\!\!  * \partial_y \!\!&\!\!  * \partial_x \!\!&\!\!  * \partial_y \!\!&\!\!
 * \partial_x \!\!&\!\! 0 \!\!&\!\! 0 \!\!&\!\! 0 \\
0 \!\!&\!\! 0 \!\!&\!\! 0  & %% debut D
 * \partial_x  \!\!&\!\!  * \partial_x \!\!&\!\!   * \partial_y \!\!&\!\! 0 \!\!&\!\! 0 \!\!&\!\!0 \!\!&\!\!0 \!\!&\!\!
 * \partial_x \!\!&\!\!  * \partial_x  \!\!&\!\! 0 \\
0\!\!&\!\! 0 \!\!&\!\! 0 & %% debut D
 * \partial_y   \!\!&\!\!  * \partial_y  \!\!&\!\!  * \partial_x \!\!&\!\! 0 \!\!&\!\! 0 \!\!&\!\! 0 \!\!&\!\! 0 \!\!&
\!\!  * \partial_y  \!\!&\!\!  * \partial_y &\!\! 0 \\
0 \!\!&\!\! 0 \!\!&\!\! 0 & %% debut D
0  \!\!&\!\!  * \partial_x  \!\!&\!\!  * \partial_y \!\!&\!\! 0 \!\!&\!\!
0 \!\!&\!\! 0 \!\!&\!\! 0 \!\!& \!\!  * \partial_x  \!\!&\!\!  * \partial_x &\!\!  * \partial_x \\
0 \!\!&\!\! 0 \!\!&\!\! 0  & %% debut D
0  \!\!&\!\!  * \partial_y  \!\!&\!\!  * \partial_x \!\!&\!\! 0 \!\!&\!\! 0 \!\!&\!\! 0 \!\!&\!\! 0 \!\!&
\!\!  * \partial_y  \!\!&\!\!  * \partial_y &\!\!  * \partial_y \\
0 \!\!&\!\! 0 \!\!&\!\! 0 & %% debut D
0   \!\!&\!\! 0 \!\!&\!\! 0 \!\!&\!\!  * \partial_x \!\!&\!\!  * \partial_y \!\!&\!\! * \partial_x \!\!&\!\!
 * \partial_y  \!\!&\!\! 0 \!\!&\!\! 0  \!\!&\!\! 0  \\
0 \!\!&\!\! 0 \!\!&\!\! 0  & %% debut D
0  \!\!&\!\! 0 \!\!&\!\! 0 \!\!&\!\!  * \partial_x \!\!&\!\!  * \partial_y \!\!&\!\! * \partial_x \!\!&\!\!
 * \partial_y  \!\!&\!\! 0 \!\!&\!\! 0  \!\!&\!\! 0   \\
0 \!\!&\!\! 0 \!\!&\!\! 0  & %% debut D
0  \!\!&\!\! 0 \!\!&\!\! 0 \!\!&\!\! 0 \!\!&\!\! 0 \!\!&\!\!  * \partial_x \!\!&\!\!  * \partial_y \!\!&\!\!
0 \!\!&\!\! 0 \!\!&\!\! 0 
\end{array} \right] \,. \monend 
%

%%%%%%%%%%%%%%%%%%%%%%%%%%%%%%%%%%%%%%%%%%%%%%%%%%%%%%%%%%%%%%%%%%%%%%%%%%%%%%%%%  d2q13  isotherme ordre 1 
\smallskip \monitem 
In order to simulate gas dynamics equations, we confront the equivalent partial differential equations
of the D2Q13 scheme with 3 conservations
\moneq \label{d2q13-edp-ordre-1}  \left \{  \begin {array}{l} 
  \partial_t \rho + \partial_x j_x + \partial_y j_y = {\rm O}(\Delta t) \\
  \partial_t j_x + \partial_x ( {14\over13} \, \lambda^2 \, \rho + {1\over26} \Phi_\varepsilon +   {1\over2}  \Phi_{xx} )
  + \partial_y   \Phi_{xy} = {\rm O}(\Delta t) \\
  \partial_t j_y + \partial_x  \Phi_{xy} + \partial_y ( {14\over13} \, \lambda^2 \, \rho + {1\over26} \Phi_\varepsilon -   {1\over2}  \Phi_{xx} )  = {\rm O}(\Delta t)
\end{array} \right.  \monend
to the Euler equations (\ref{euler-isotherme-2d}). 
%
%% \moneq \label{euler-2d}  \left \{  \begin {array}{l} 
%
The two systems (\ref{d2q13-edp-ordre-1}) and (\ref{euler-isotherme-2d}) must coincide at first order for all the solutions.
Then we have the following system of equations 
\moneqstar  \left \{  \begin {array}{l}
{14\over13} \, \lambda^2 \, \rho + {1\over26} \Phi_\varepsilon + {1\over2}  \Phi_{xx}  =  \rho \, u^2 + p \\ 
{14\over13} \, \lambda^2 \, \rho + {1\over26} \Phi_\varepsilon - {1\over2}  \Phi_{xx}  =  \rho \, v^2 + p \\
\Phi_{xy}  =  \rho \, u \, v \, . 
\end{array} \right.  \monendstar 
After two lines of elementary algebra, we obtain the values of equilibria for the first family of nonconserved moments:
\moneq \label{d2q13-moments-equilibre-2e-ns-isotherme}  
\Phi_\varepsilon =\rho \, (  13 \, | {\bf u} |^2 - 28  \, \lambda^2 + 26  \, c_s^2  ) \,,\,\,  
\Phi_{xx} =  \rho \, (  u^2 - v^2 )   \,,\,\,  \Phi_{xy} =  \rho \,  u \, v     \monend
%
%%%  \moneq \label{d2q13-moments-equilibre-2e-ns-isotherme}  \left \{ \begin {array}{l}
%%%  \Phi_\varepsilon =\rho \, (  13 \, | {\bf u} |^2 - 28  \, \lambda^2 + 26  \, c_s^2  ) \\ 
%%%  \Phi_{xx} =  \rho \, (  u^2 - v^2 )   \,,\,\,  \Phi_{xy} =  \rho \,  u \, v   \, ,  
%%%  \end {array} \right. \monend
%
with $ \, p =  \rho \, c_s^2 \, $ and an arbitrary speed of sound $  \, c_s $. 
%

%%%%%%%%%%%%%%%%%%%%%%%%%%%%%%%%%%%%%%%%%%%%%%%%%%%%%%%%%%%%%%%%%%%%%%%%%%%%%%%%%  d2q13  isotherme ordre 2
\smallskip \monitem 
The  second  order equations
are constructed with the help of
\moneqstar 
\Gamma_1 = \big( \partial_x j_x + \partial_y j_y  \,,\,\, 
\partial_x ( \rho \, u^2 + p )  + \partial_y ( \rho \, u \, v ) \,,\,\, 
\partial_x ( \rho \, u \, v ) +  \partial_x ( \rho \, v^2 + p ) \big)^{\rm t}
\monendstar
and 
\moneqstar
\Psi_1 = \dd \Phi(W) .  \Gamma_1  - (C \, W + D \, \Phi(W))   \, . 
\monendstar
The question now is to identify the two expressions of the viscous fluxes:
\moneqstar
\Delta t \, \Gamma_2 =   \Delta t \,\,  B \, \Sigma \, \Psi_1
\monendstar
on one hand and the viscous dissipation $ \, {\rm div} \, \tau \, $ detailed at relation (\ref{flux-visqueux-2d}) 
%
%%%%   \moneqstar
%%%%   - {\rm div} \, \tau  =  - \left( \begin{array} {c} 0 \\ \partial_j \tau_{xj} \equiv \partial_x ( 2 \, \mu \,   \partial_x u 
%%%%   + (\zeta-\mu)  ( \partial_x u   +   \partial_y v  ) )
%%%%   + \partial_y (\mu ( \partial_x v   +   \partial_y u  ) )  \\
%%%%   \partial_j \tau_{yj}   \equiv \partial_x (\mu (   \partial_x v  +    \partial_y u  ) )
%%%%   + \partial_y (  (\zeta-\mu) (   \partial_x u   +   \partial_y v  ) + 2 \, \mu \,  \partial_y v  ))
%%%%     \end{array} \right) \monendstar
%
on the other hand.
As for the D2Q9 scheme, 
we are confronted to a system of  $ \, 24 = 2 \times 2 \times 3 \times 2   \, $ equations.
%% 2 equations for momentum, 2 partial derivatives of the type $ \, \partial_x \, $ and  $ \, \partial_y \, $
%% for each equation, 3~nonconserved variables $ \, \rho $, $ \, u \, $ and $ \, v $,
%% 2 partial derivatives  $ \, \partial_x \, $ and  $ \, \partial_y \, $ for each of these variables.
We have 4 moments
$ \, q_x $, $ \, q_y $, $ \, r_x \, $ and $ \, r_y \, $ 
that are present in the algebraic expression of $ \, \Gamma_2 $, and
3  partial derivatives  relative to $ \, \rho $, $ \, u $,  and $ \, v \, $  for each of these moments,
then a total of $ \, 4 \times 3 = 12 \, $ unknowns.
We obtain dependent linear equations and  this system has a unique solution.
Then it is easy to integrate the partial differential equations that are linear relative to the density.
We obtain finally necessary values of equilibrium moments: 
\moneq \label{d2q13-moments-equilibre-3e-ns-isotherme}  \left \{ \begin {array}{l}
\Phi_{qx} =  \rho \, u \, \big( u^2 + v^2 + 4\, \lambda^2 \, c_s^2 - 3 \, \lambda^2 \big)  \\ 
\Phi_{qy} =  \rho \, v \, \big( u^2 + v^2 + 4\, \lambda^2 \, c_s^2 - 3 \, \lambda^2 \big)  \\ 
\Phi_{rx} =  \rho \, u \, \lambda^2  \, \big( -{7\over6} \, u^2   -7 \,  v^2  - {21\over2} \, c_s^2 + {31\over6} \, \lambda^2 \big)   \\ 
\Phi_{ry} =  \rho \, v \,  \lambda^2  \, \big( -7 \, u^2 - {7\over6} \, v^2 - {21\over2}  \, c_s^2 + {31\over6} \, \lambda^2 \big) \, . 
\end{array} \!\!\! \right.  \monend
We observe that $ \, \Phi_{qx} \, $ and $ \, \Phi_{qy} \, $ define the components of a vector
whereas it is not the case for the functions $ \, \Phi_{rx} \, $ and $ \, \Phi_{ry}  $. 
We made this  algebraic calculus,  a first  time without any software \cite{Du2015} 
and in a second step with the help of SageMath \cite{sagemath}.
Our results are coherent. 
Finally, the viscosities $ \, \mu \, $ and $ \, \zeta \, $ satisfy 
\moneqstar
\mu = \, \rho \, \sigma_x  \, \lambda \, c_s^2 \, \Delta x  \,,\,\,
\zeta = \, \rho \, \sigma_e  \, \lambda \, c_s^2 \, \Delta x  \,. 
  \monendstar
  %

%%%   \newpage 
%%%%%%%%%%%%%%%%%%%%%%%%%%%%%%%%%%%%%%%%%%%%%%%%%%%%%%%%%%%%%%%%%%%%%%%%%%%%%%%  section 5
\bigskip \bigskip    \noindent {\bf \large    5) \quad  Three-dimensional isothermal Navier Stokes}
%%%%%%%%%%%%%%%%%%%%%%%%%%%%%%%%%%%%%%%%%%%%%%%%%%%%%%%%%%%%%%%%%%%%%%%%%%%%%%%%%%%%%%%%%%%%%%%%%%%%%%%%%%

\noindent
For the three-dimensional space, two lattice Boltzmann schemes are  popular~:
the D3Q19 and the D3Q27   schemes. We  show in the following that two other schemes,
the D3Q33 and the ``D3Q27-2'' scheme have a lot of  interest. 

%%% \newpage 
%%%%%%%%%%%%%%%%%%%%%%%%%%%%%%%%%%%%%%%%%%%%%%%%%%%%%%%%%%%%%%%%%%%%%%%%%%%%%%
\smallskip \monitem  D3Q19
%%%%%%%%%%%%%%%%%%%%%%%%%%%%%%%%%%%%%%%%%%%%%%%%%%%%%%%%%%%%%%%%%%%%%%%%%%%%%%

\begin{figure}    [H]  \centering 
\vspace{-.2 cm} 
%%%%  \centerline {\includegraphics[height = .55 \textwidth] {lbm-graphes-d3q19-01juin2020.pdf}}
%%%   \centerline {\includegraphics[height = .45 \textwidth] {lbm-graphes-d3q19-15mars2021.pdf}} 
\centerline {\includegraphics[height = .33 \textwidth] {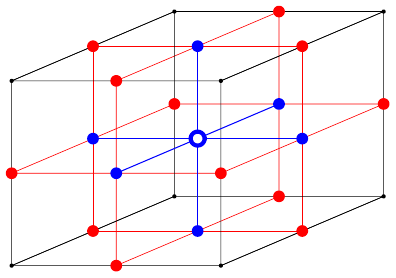}} 
\vspace{-.3 cm} 
\caption{Set of discrete velocities for the D3Q19 lattice Boltzmann scheme} 
\label{fig-d3q19} \end{figure}

\noindent
From  the cubic lattice D3Q27, we omit the corners of the cube and define in this way the
D3Q19 scheme represented in Figure \ref{fig-d3q19}. 

%%%%%%%%%%%%%%%%%%%%%%%%%%%%%%%%%%%%%%%%%%%%%%%%%%%%%%%%%%%%%%%%%%%%%%%%%%%%%%%%%%% figure d3q19
%%%%%%%%%%%%%%%%%%%%%%%%%%%%%%%%%%%%%%%%%%%%%%%%%%%%%%%%%%%%%%%%%%%%%%%%%%%%%%%%%%%
%
The definition of the moments is  detailed in Table \ref{d3q19-iso-polynomes} of the Appendix.
We have now 4 conserved moments and the nonconserved moments are 
organized into three families as in the previous schemes (see Table \ref{d3q19-moments-isotherme}).
%
%%%%%%%%%%%%%%%%%%%%%%%%%%%%%%%%%%%%%%%%%%%%%%%%%%%%%%%%%%%%%%%%%%%%%%%%%%%%%%%%%%% moments d3q19   
\begin{table}    [H]  \centering
\begin{tabular}{|c|l|c|c|} \hline
&conserved &  $\rho \,,\,\, j_x \,,\,\, j_y \,,\,\, j_z $ & 4  \\  \hline
1 & fit the Euler equations  &  $ \varepsilon \,,\,\,  xx  \,,\,\,  ww \,,\,\, xy   \,,\,\,  yz \,,\,\,  zx $ & 6 \\ \hline 
2 & fit the viscous terms? & $ q_x  \,,\,\, q_y  \,,\,\, q_z \,,\,\,  x_{yz}  \,,\,\,  y_{zx}   \,,\,\,   z_{xy} $  & 6  \\ \hline  
3 & without influence & $ h  \,,\,\, xx_e  \,,\,\,   ww_e   $  & 3   \\ \hline 
\end{tabular} 
\caption{The four families of moments for the D3Q19  scheme  for the approximation of the isothermal Navier Stokes equations}
\label{d3q19-moments-isotherme} \end{table}
%%%%%%%%%%%%%%%%%%%%%%%%%%%%%%%%%%%%%%%%%%%%%%%%%%%%%%%%%%%%%%%%%%%%%%%%%%%%%%%%%%%
%
The non null elements of the velocity-momentum operator matrix for the
D3Q19 lattice Boltzmann scheme are presented in the following relation, 
where a null symbol indicates that the corresponding element is equal to zero
and a star symbol that it is a non zero space differential operator: 
\moneqstar 
\Lambda_{D3Q19}^{\rm iso}  =   \left[   \begin{array} {cccccccccccccccccccc}
0&*&*&*&&0&0&0&0&0&0&0&0&0&0&0&0&0&0&0\\ \vspace{-.09 cm}
*&0&0&0&&*&*&0&*&0&*&0&0&0&0&0&0&0&0&0\\ \vspace{-.09 cm}
*&0&0&0&&*&*&*&*&*&0&0&0&0&0&0&0&0&0&0\\ \vspace{-.09 cm}
*&0&0&0&&*&*&*&0&*&*&0&0&0&0&0&0&0&0&0\\ \vspace{-.09 cm} \\ \vspace{-.09 cm}
0&*&*&*&&0&0&0&0&0&0&*&*&*&0&0&0&0&0&0\\ \vspace{-.09 cm}
0&*&*&*&&0&0&0&0&0&0&*&*&*&0&*&*&0&0&0\\ \vspace{-.09 cm}
0&0&*&*&&0&0&0&0&0&0&0&*&*&*&*&*&0&0&0\\ \vspace{-.09 cm}
0&*&*&0&&0&0&0&0&0&0&*&*&0&*&*&0&0&0&0\\ \vspace{-.09 cm}
0&0&*&*&&0&0&0&0&0&0&0&*&*&0&*&*&0&0&0\\ \vspace{-.09 cm}
0&*&0&*&&0&0&0&0&0&0&*&0&*&*&0&*&0&0&0\\ \vspace{-.09 cm}
0&0&0&0&&*&*&0&*&0&*&0&0&0&0&0&0&*&*&0\\ \vspace{-.09 cm}
0&0&0&0&&*&*&*&*&*&0&0&0&0&0&0&0&*&*&*\\ \vspace{-.09 cm}
0&0&0&0&&*&*&*&0&*&*&0&0&0&0&0&0&*&*&*\\ \vspace{-.09 cm}
0&0&0&0&&0&0&*&*&0&*&0&0&0&0&0&0&0&0&*\\ \vspace{-.09 cm}
0&0&0&0&&0&*&*&*&*&0&0&0&0&0&0&0&0&*&*\\ \vspace{-.09 cm}
0&0&0&0&&0&*&*&0&*&*&0&0&0&0&0&0&0&*&*\\ \vspace{-.09 cm}
0&0&0&0&&0&0&0&0&0&0&*&*&*&0&0&0&0&0&0\\ \vspace{-.09 cm}
0&0&0&0&&0&0&0&0&0&0&*&*&*&0&*&*&0&0&0\\
0&0&0&0&&0&0&0&0&0&0&0&*&*&*&*&*&0&0&0
\end{array} \right] \, .   \monendstar 
\noindent
%%%  
%%%  [                                                                   rho*uu*dx + rho*vv*dy + rho*ww*dz]
%%%  [                          1/57*Phi_ee*dx + 1/3*Phi_xx*dx + Phi_xy*dy + Phi_zx*dz + 10*rho*ll^2/19*dx]
%%%  [     Phi_xy*dx + 1/57*Phi_ee*dy + (1/(-6))*Phi_xx*dy + 1/2*Phi_ww*dy + Phi_yz*dz + 10*rho*ll^2/19*dy]
%%%  [Phi_zx*dx + Phi_yz*dy + 1/57*Phi_ee*dz + ((-1)/6)*Phi_xx*dz + ((-1)/2)*Phi_ww*dz + 10*rho*ll^2/19*dz]
%%%
All the coefficients of this matrix can be found in \cite{Du22git}.

%%%%%%%%%%%%%%%%%%%%%%%%%%%%%%%%%%%%%%%%%%%%%%%%%%%%%%%%%%%%%%%%%%%%%%%%%%%%%%%%%  d3q19  isotherme ordre 1 
\smallskip \monitem 
At first order of accuracy, the partial equivalent equations can be written as follows
\moneqstar    \left \{ \begin {array}{l}
\partial_t \rho + \partial_x j_x + \partial_y j_y  + \partial_z j_z  = 0  \\
\partial_t j_x + \partial_x \big( {1\over57} \, \Phi_\varepsilon + {1\over3} \,  \Phi_{xx} +  {10\over19} \,\rho \, \lambda^2 \big)
+  \partial_y \Phi_{xy}  + \partial_z \Phi_{zx} = 0 \\ 
\partial_t j_y + \partial_x \Phi_{xy} + \partial_y \big(  {1\over57} \, \Phi_\varepsilon  - {1\over6} \,  \Phi_{xx} + {1\over2} \,  \Phi_{ww} 
+  {10\over19} \,\rho \, \lambda^2 \big) +  \partial_z \Phi_{yz} = 0 \\
\partial_t j_z + \partial_x \Phi_{zx} + \partial_y \Phi_{yz}  + \partial_z \big(  {1\over57} \, \Phi_\varepsilon
- {1\over6} \,  \Phi_{xx} - {1\over2} \,  \Phi_{ww} +  {10\over19} \,\rho \, \lambda^2 \big) = 0 
\end {array} \right. \monendstar 
and they must be identical to the Euler equations of gas dynamics 
\moneq \label{euler-3d-sans-energie}   \left \{ \begin {array}{l}
\partial_t \rho + \partial_x j_x + \partial_y j_y  + \partial_z j_z  = 0   \\
\partial_t j_x + \partial_x ( \rho \, u^2 + p )  + \partial_y ( \rho \, u \, v )  + \partial_z ( \rho \, u \, w) = 0  \\
\partial_t j_y + \partial_x ( \rho \, u \, v ) +  \partial_y ( \rho \, v^2 + p ) + \partial_z ( \rho \, v \, w ) =  0 \\
\partial_t j_z + \partial_x ( \rho \, u \, w ) +  \partial_y ( \rho \, v \, w ) + \partial_z ( \rho \, w^2 + p ) =  0 \,. 
\end {array} \right. \monend
The equlibrium values $\, \Phi_\varepsilon $, $\,  \Phi_{xx} $, $\, \Phi_{ww} $, $\, \Phi_{xy} $, $\, \Phi_{yz} \, $
and $ \, \Phi_{zx} \, $ of the moments of the second family satisfy the set of equations 
\moneqstar    \left \{ \begin {array}{l}
{1\over57} \, \Phi_\varepsilon + {1\over3} \,  \Phi_{xx} +  {10\over19} \,\rho \, \lambda^2  = \rho \, u^2 + p \\
{1\over57} \, \Phi_\varepsilon  - {1\over6} \,  \Phi_{xx} + {1\over2} \,  \Phi_{ww} +  {10\over19} \,\rho \, \lambda^2 =  \rho \, v^2 + p \\
{1\over57} \, \Phi_\varepsilon - {1\over6} \,  \Phi_{xx} - {1\over2} \,  \Phi_{ww} +  {10\over19} \,\rho \, \lambda^2  = \rho \, w^2 + p 
\end {array} \right. \monendstar
and we have also 
\moneqstar 
\Phi_{xy} = \rho \, u \, v \,,\,\,  \Phi_{yz} =  \rho \, v \, w \,,\,\,  \Phi_{zx} = \rho \, u \, w   . 
\monendstar
With a state pressure law $ \, p = \rho \, c_s^2 \, $ and a  not  yet imposed  value of the sound velocity $\, c_s $,   
we recover the  first order isothermal equations (\ref{euler-3d-sans-energie}) by fixing the equilibrium value of the three last moments
$ \, \varepsilon $, $ \, xx \, $ and $ \, ww \, $ of the second family of moments: 
% 
%%% \moneq \label{d3q19-moments-equilibre-2e-ns-isotherme}  \left \{ \begin {array}{l}
%
\moneqstar 
\Phi_\varepsilon =\rho \, (  19 \, | {\bf u} |^2 - 30  \, \lambda^2 + 57 \, c_s^2  ) \,,\,\, 
\Phi_{xx} =  \rho \, ( 2 \, u^2 - v^2 - w^2 )   \,,\,\, \Phi_{ww} =  \rho \, ( v^2 - w^2 )  . 
%%%  \Phi_{xy} =  \rho \,  u \, v  \,,\,\,   \Phi_{yz} =  \rho \,  v \, w \,,\,\,   \Phi_{yz} =  \rho \, w \, u \, . 
%%% \end {array} \right. \monend
\monendstar
%
%

%%%%%%%%%%%%%%%%%%%%%%%%%%%%%%%%%%%%%%%%%%%%%%%%%%%%%%%%%%%%%%%%%%%%%%%%%%%%%%%%%  d3q19  isotherme ordre 2 
\smallskip \monitem 
The next step is the identification with  the second order Navier Stokes equations for
mass and momentum. We must solve a total of $ \, 108 = 3 \times 3 \times 4 \times 3 \, $
equations  to identify the second order terms of the Navier Stokes equations:
one equation for each component of the  momentum,
3 conservation terms per equation : $\, \partial_x  [**] $, $\, \partial_y  [**] \, $ and $\, \partial_z  [**] $, 
4 nonconserved variables $\, \rho $, $ \, u $, $ \, v \, $ and $ \, w \, $
and   3 partial derivatives $\, \partial_x $, $\, \partial_y \, $ and  $\, \partial_z \,$  per variable. 
We have a total of  four partial derivatives 
$\, \partial_\rho $, $ \, \partial_u $, $ \, \partial_v \, $ and $ \, \partial_w \, $ for  each equilibrium of the 6 moments
$\,  q_x $, $\, q_y   $, $\,  q_z  $, $\,  x_{yz}   $, $\, y_{zx} \, $ and $ \,  z_{xy} \, $ of the second family 
that has an influence on second order terms. Thus we have (only) 24 unknowns and have to  solve 108 equations. 
This algebraic problem has no solution.
Nevertheless,  enforcing  the usual value $ \, c_s = {{\lambda}\over{\sqrt{3}}} \, $ and with a specific choice
of the second family of moments, {\it id est} 
\moneq \label{d3q19-moments-equilibre-3e-ns-isotherme-isotropic}  \left \{ \! \begin {array}{l} 
\Phi_{q_x} =  \rho \,  u \, \xi_q \,,\,\,
\Phi_{q_y} =  \rho \,  v \, \xi_q \,,\,\,
\Phi_{q_z} =  \rho \,  w \, \xi_q \,,\,\,
%%% {\color {blue}
\xi_q \equiv 5 \,  | {\bf u} |^2  -  {2\over3} \, \lambda^2    \\ %%% } 
%%%    correction d une coquille suite a la discussion avec Mahdi le 30 mars 2023 
\Phi_{x_{yz}} = \rho \,  u \, (v^2 - w^2) \,,\,\,  
\Phi_{y_{zx}} = \rho \,  v \, (w^2 - u^2) \,,\,\,  
\Phi_{z_{xy}} = \rho \,  w \, (u^2 - v^2) \, , \\
\end {array} \right. \monend
a beginning of a resolution of the system can be stated.
But a total of  61 equations  remain unsolved. 
We deduce from the asymptotic analysis
the value of the shear and the bulk viscosities: 
\moneq \label{visosites-d3q19}
%%%  \mu =   {1\over3} \, \, \rho \,  \sigma_x \, \Delta t \, ( \lambda^2 + {\rm O}( | {\bf u} |^2 ) )  \,,\,\, 
%%%  \zeta = {2\over9} \, \rho \, \sigma_e \, \Delta t \,  ( \lambda^2 + {\rm O}( | {\bf u} |^2 ) )  \, . 
\mu =   {1\over3} \, \, \rho \,  \sigma_x \, \Delta t \, \lambda^2    \,,\,\, 
\zeta = {2\over9} \, \rho \, \sigma_e \, \Delta t \,  \lambda^2   \, . 
\monend
Then the components of the tensor of viscosities admit the classical relations (\ref{tau-3d}). 
The second order equivalent partial differential equations of the D3Q19
lattice Boltzmann scheme with 4 conserved moments can be written as 
\moneqstar \label{d3q19-ns-isotherme} \left \{ \begin {array}{l}
\partial_t \rho + {\textrm {div}} ( \rho \, {\bf u}  ) =  {\textrm O}(\Delta x^2)  \\ % \,,\,\,
\partial_t ( \rho \, {\bf u}  ) + {\textrm {div}} \,   ( \rho \, {\bf u} \otimes {\bf u}  ) + \nabla p
- \Delta t \,\,  {\textrm {div}} \, {\bf \tau} + \sigma_x \,    \Delta t  \,\,     {\textrm {div}} R =  {\textrm O}(\Delta x^2 ) \,.
\end {array} \right. \monendstar 
The  tensor of discrepancy $ \, R \, $ is symmetric. It can be explicited  through the relations
\moneqstar \left \{ \! \! \begin {array}{l}
R_{xx} =  u^3 \, \partial_x \rho - {{v^3}\over2}  \, \partial_y \rho - {{w^3}\over2}  \, \partial_z \rho 
+ 3 \, \rho \, (u^2 \, \partial_x u  - {{v^2}\over2} \, \partial_y v - {{w^2}\over2} \, \partial_z w) \\ 
R_{yy} =  - {{u^3}\over2}  \, \partial_x \rho +  v^3 \, \partial_y \rho - {{w^3}\over2}  \, \partial_z \rho  
+ 3 \, \rho \, ( - {{u^2}\over2} \, \partial_x u +  v^2 \, \partial_y v  - {{w^2}\over2} \, \partial_z w) \\ 
R_{zz} =   - {{u^3}\over2}  \, \partial_x \rho  - {{v^3}\over2}  \, \partial_y \rho +   w^3 \, \partial_z \rho
+ 3 \, \rho \, ( - {{u^2}\over2} \, \partial_x u  - {{v^2}\over2} \, \partial_y v +  w^2 \, \partial_z w) \\ 
R_{xy} = - {{v^3}\over2}  \, \partial_x \rho  - {{u^3}\over2} \, \partial_y \rho  + u \, v \, w \, \partial_z \rho \\ 
\qquad \qquad \qquad   +  \rho \, \big(  - {{3}\over2} \, (  v^2 \, \partial_x v + u^2 \,\partial_y u ) + w \, ( v \,\partial_z u + u \, \partial_z v ) \big) + u \, v \,  \partial_z w \big) \\
R_{yz} =  u \, v \, w \, \partial_x \rho  - {{w^3}\over2}  \, \partial_y \rho  - {{v^3}\over2}  \, \partial_z \rho \\ 
\qquad \qquad \qquad   +  \rho \, \big( - {{3}\over2} \, ( w^2 \, \partial_y w + v^2 \,\partial_z v )    + u \, ( w \, \partial_x v + v \,  \partial_x w) + v \, w \, \partial_x u \big) \\ 
R_{zx} =   - {{w^3}\over2}  \, \partial_x \rho +  u \, v \, w \, \partial_y \rho - {{u^3}\over2}  \, \partial_z \rho \\
\qquad \qquad \qquad   +  \rho \, \big(  - {{3}\over2} \, ( u^2 \,\partial_z u +  w^2 \,\partial_x w ) + v \, ( u \, \partial_y w  + w \, \partial_y u ) + w \, u \,  \partial_y v \big) 
\end {array} \right. \monendstar
and the vector $\,    {\textrm {div}} R \, $ admits the expression
\moneqstar 
{\textrm {div}} R = \big(\partial_x R_{xx} + \partial_y R_{xy} + \partial_z R_{zx} \,,\,\,  \partial_x R_{xy} + \partial_y R_{yy} + \partial_z R_{yz}  \,,\,\, 
\partial_x R_{zx} + \partial_y R_{yz} + \partial_z R_{zz}  \big)^{\rm t}  \, . 
 \monendstar
This result is well known and we have essentially proposed a reformulation of the results presented in  \cite{DGKLL02}.

%% \newpage 
%%%%%%%%%%%%%%%%%%%%%%%%%%%%%%%%%%%%%%%%%%%%%%%%%%%%%%%%%%%%%%%%%%%%%%%%%%%%%%%%%%%%
\smallskip \monitem D3Q27
%%%%%%%%%%%%%%%%%%%%%%%%%%%%%%%%%%%%%%%%%%%%%%%%%%%%%%%%%%%%%%%%%%%%%%%%%%%%%%%%%%%%

\begin{figure}    [H]  \centering 
\vspace{-.4 cm} 
\centerline {\includegraphics[height=.41\textwidth]{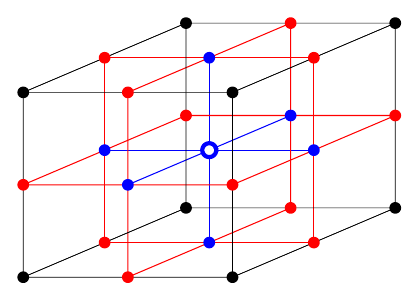}}
\vspace{-.2 cm} 
\caption{Set of discrete velocities for the D3Q27 lattice Boltzmann scheme} 
\label{fig-d3q27} \end{figure}

\noindent
Adding 8 velocities of the type $ \, (1,\, 1 ,\, 1 ) \, $ to the D3Q19 stencil, we obtain the first neighbours
of a cubic lattice, as presented in Figure \ref{fig-d3q27}. 

%%%%%%%%%%%%%%%%%%%%%%%%%%%%%%%%%%%%%%%%%%%%%%%%%%%%%%%%%%%%%%%%%%%%%%%%%%%%%%%%%%% figure d3q27
%%%%%%%%%%%%%%%%%%%%%%%%%%%%%%%%%%%%%%%%%%%%%%%%%%%%%%%%%%%%%%%%%%%%%%%%%%%%%%%%%%%

%%%%%%%%%%%%%%%%%%%%%%%%%%%%%%%%%%%%%%%%%%%%%%%%%%%%%%%%%%%%%%%%%%%%%%%%%%%%%%%%%%% moments d3q27   
\begin{table}    [H]  \centering
\begin{tabular}{|c|l|c|c|} \hline
& conserved &  $\rho \,,\,\, j_x \,,\,\, j_y \,,\,\, j_z $ & 4  \\  \hline
1 & fit the Euler equations  &  $ \varepsilon \,,\,\,  xx  \,,\,\,  ww \,,\,\, xy   \,,\,\,  yz \,,\,\,  zx $ & 6 \\ \hline 
2 & fit the viscous terms ?  & $ q_x  \,,\,\, q_y  \,,\,\, q_z \,,\,\,  x_{yz}  \,,\,\,  y_{zx}   \,,\,\,   z_{xy} \,,\,\, xyz $  & 7  \\ \hline  
3 & without influence & $ r_x  \,,\,\, r_y  \,,\,\, r_z \,,\,\, h  \,,\,\,
xx_e  \,,\,\,   ww_e  \,,\,\,   xy_e  \,,\,\,  yz_e  \,,\,\,  zx_e  \,,\,\,   h_3 $  & 10  \\ \hline 
\end{tabular} 
\caption{The four families of moments for the D3Q27  scheme  for the approximation of the isothermal Navier Stokes equations}
\label{d3q27-moments-isotherme} \end{table}
%%%%%%%%%%%%%%%%%%%%%%%%%%%%%%%%%%%%%%%%%%%%%%%%%%%%%%%%%%%%%%%%%%%%%%%%%%%%%%%%%%%
%
Among the 27 moments, we still have 4 conserved moments 
and 6 moments that are used to fit first order terms (first family of nonconserved moments). 
We have also 7 moments  to fix the viscous fluxes instead of 6 for the D3Q19 scheme (second family), and 10 moments have no direct influence on
the equivalent partial differential equations at second order (third family).
The usual names of these moments are presented in Table \ref{d3q27-moments-isotherme}. 
The algebraic formulas that define the polynomials associated to the previous moments
 with the relation (\ref{moments-polynomes}) are explicited in the Appendix at the Table \ref{d3q27-iso-polynomes}.
With this choice, the  non null elements of the momentum-velocity matrix operator  $ \, \Lambda \, $
are given at locations identified by stars in the following expression: 
\smallskip \moneqstar 
\Lambda_{D3Q27}^{\rm iso} =  \left[ \begin{array} {cccccccccccccccccccccccccccc}   \vspace{-.2 cm}
\scriptstyle 0\!\!&\!\!*\!\!&\!\!*\!\!&\!\!*\!\!&&\!\!\scriptstyle 0\!\!&\!\!\scriptstyle 0\!\!&\!\!\scriptstyle 0\!\!&\!\!\scriptstyle 0\!\!&\!\!\scriptstyle 0\!\!&\!\!\scriptstyle 0\!\!&\!\!\scriptstyle 0\!\!&\!\!\scriptstyle 0\!\!&\!\!\scriptstyle 0\!\!&\!\!\scriptstyle 0\!\!&\!\!\scriptstyle 0\!\!&\!\!\scriptstyle 0\!\!&\!\!\scriptstyle 0\!\!&\!\!\scriptstyle 0\!\!&\!\!\scriptstyle 0\!\!&\!\!\scriptstyle 0\!\!&\!\!\scriptstyle 0\!\!&\!\!\scriptstyle 0\!\!&\!\!\scriptstyle 0\!\!&\!\!\scriptstyle 0\!\!&\!\!\scriptstyle 0\!\!&\!\!\scriptstyle 0\!\!&\!\!\scriptstyle 0\\ \vspace{-.2 cm}
*\!\!&\!\!\scriptstyle 0\!\!&\!\!\scriptstyle 0\!\!&\!\!\scriptstyle 0\!\!&&\!\!*\!\!&\!\!*\!\!&\!\!\scriptstyle 0\!\!&\!\!*\!\!&\!\!\scriptstyle 0\!\!&\!\!*\!\!&\!\!\scriptstyle 0\!\!&\!\!\scriptstyle 0\!\!&\!\!\scriptstyle 0\!\!&\!\!\scriptstyle 0\!\!&\!\!\scriptstyle 0\!\!&\!\!\scriptstyle 0\!\!&\!\!\scriptstyle 0\!\!&\!\!\scriptstyle 0\!\!&\!\!\scriptstyle 0\!\!&\!\!\scriptstyle 0\!\!&\!\!\scriptstyle 0\!\!&\!\!\scriptstyle 0\!\!&\!\!\scriptstyle 0\!\!&\!\!\scriptstyle 0\!\!&\!\!\scriptstyle 0\!\!&\!\!\scriptstyle 0\!\!&\!\!\scriptstyle 0\\ \vspace{-.2 cm}
*\!\!&\!\!\scriptstyle 0\!\!&\!\!\scriptstyle 0\!\!&\!\!\scriptstyle 0\!\!&&\!\!*\!\!&\!\!*\!\!&\!\!*\!\!&\!\!*\!\!&\!\!*\!\!&\!\!\scriptstyle 0\!\!&\!\!\scriptstyle 0\!\!&\!\!\scriptstyle 0\!\!&\!\!\scriptstyle 0\!\!&\!\!\scriptstyle 0\!\!&\!\!\scriptstyle 0\!\!&\!\!\scriptstyle 0\!\!&\!\!\scriptstyle 0\!\!&\!\!\scriptstyle 0\!\!&\!\!\scriptstyle 0\!\!&\!\!\scriptstyle 0\!\!&\!\!\scriptstyle 0\!\!&\!\!\scriptstyle 0\!\!&\!\!\scriptstyle 0\!\!&\!\!\scriptstyle 0\!\!&\!\!\scriptstyle 0\!\!&\!\!\scriptstyle 0\!\!&\!\!\scriptstyle 0\\
 \vspace{-.2 cm}
*\!\!&\!\!\scriptstyle 0\!\!&\!\!\scriptstyle 0\!\!&\!\!\scriptstyle 0\!\!&&\!\!*\!\!&\!\!*\!\!&\!\!*\!\!&\!\!\scriptstyle 0\!\!&\!\!*\!\!&\!\!*\!\!&\!\!\scriptstyle 0\!\!&\!\!\scriptstyle 0\!\!&\!\!\scriptstyle 0\!\!&\!\!\scriptstyle 0\!\!&\!\!\scriptstyle 0\!\!&\!\!\scriptstyle 0\!\!&\!\!\scriptstyle 0\!\!&\!\!\scriptstyle 0\!\!&\!\!\scriptstyle 0\!\!&\!\!\scriptstyle 0\!\!&\!\!\scriptstyle 0\!\!&\!\!\scriptstyle 0\!\!&\!\!\scriptstyle 0\!\!&\!\!\scriptstyle 0\!\!&\!\!\scriptstyle 0\!\!&\!\!\scriptstyle 0\!\!&\!\!\scriptstyle 0\\
\vspace{-.3 cm} \\ \vspace{-.2 cm}
\scriptstyle 0\!\!&\!\!*\!\!&\!\!*\!\!&\!\!*\!\!&&\!\!\scriptstyle 0\!\!&\!\!\scriptstyle 0\!\!&\!\!\scriptstyle 0\!\!&\!\!\scriptstyle 0\!\!&\!\!\scriptstyle 0\!\!&\!\!\scriptstyle 0\!\!&\!\!*\!\!&\!\!*\!\!&\!\!*\!\!&\!\!\scriptstyle 0\!\!&\!\!\scriptstyle 0\!\!&\!\!\scriptstyle 0\!\!&\!\!\scriptstyle 0\!\!&\!\!\scriptstyle 0\!\!&\!\!\scriptstyle 0\!\!&\!\!\scriptstyle 0\!\!&\!\!\scriptstyle 0\!\!&\!\!\scriptstyle 0\!\!&\!\!\scriptstyle 0\!\!&\!\!\scriptstyle 0\!\!&\!\!\scriptstyle 0\!\!&\!\!\scriptstyle 0\!\!&\!\!\scriptstyle 0\\ \vspace{-.2 cm}
\scriptstyle 0\!\!&\!\!*\!\!&\!\!*\!\!&\!\!*\!\!&&\!\!\scriptstyle 0\!\!&\!\!\scriptstyle 0\!\!&\!\!\scriptstyle 0\!\!&\!\!\scriptstyle 0\!\!&\!\!\scriptstyle 0\!\!&\!\!\scriptstyle 0\!\!&\!\!*\!\!&\!\!*\!\!&\!\!*\!\!&\!\!\scriptstyle 0\!\!&\!\!*\!\!&\!\!*\!\!&\!\!\scriptstyle 0\!\!&\!\!\scriptstyle 0\!\!&\!\!\scriptstyle 0\!\!&\!\!\scriptstyle 0\!\!&\!\!\scriptstyle 0\!\!&\!\!\scriptstyle 0\!\!&\!\!\scriptstyle 0\!\!&\!\!\scriptstyle 0\!\!&\!\!\scriptstyle 0\!\!&\!\!\scriptstyle 0\!\!&\!\!\scriptstyle 0\\ \vspace{-.2 cm}
\scriptstyle 0\!\!&\!\!\scriptstyle 0\!\!&\!\!*\!\!&\!\!*\!\!&&\!\!\scriptstyle 0\!\!&\!\!\scriptstyle 0\!\!&\!\!\scriptstyle 0\!\!&\!\!\scriptstyle 0\!\!&\!\!\scriptstyle 0\!\!&\!\!\scriptstyle 0\!\!&\!\!\scriptstyle 0\!\!&\!\!*\!\!&\!\!*\!\!&\!\!*\!\!&\!\!*\!\!&\!\!*\!\!&\!\!\scriptstyle 0\!\!&\!\!\scriptstyle 0\!\!&\!\!\scriptstyle 0\!\!&\!\!\scriptstyle 0\!\!&\!\!\scriptstyle 0\!\!&\!\!\scriptstyle 0\!\!&\!\!\scriptstyle 0\!\!&\!\!\scriptstyle 0\!\!&\!\!\scriptstyle 0\!\!&\!\!\scriptstyle 0\!\!&\!\!\scriptstyle 0\\ \vspace{-.2 cm}
\scriptstyle 0\!\!&\!\!*\!\!&\!\!*\!\!&\!\!\scriptstyle 0\!\!&&\!\!\scriptstyle 0\!\!&\!\!\scriptstyle 0\!\!&\!\!\scriptstyle 0\!\!&\!\!\scriptstyle 0\!\!&\!\!\scriptstyle 0\!\!&\!\!\scriptstyle 0\!\!&\!\!*\!\!&\!\!*\!\!&\!\!\scriptstyle 0\!\!&\!\!*\!\!&\!\!*\!\!&\!\!\scriptstyle 0\!\!&\!\!*\!\!&\!\!\scriptstyle 0\!\!&\!\!\scriptstyle 0\!\!&\!\!\scriptstyle 0\!\!&\!\!\scriptstyle 0\!\!&\!\!\scriptstyle 0\!\!&\!\!\scriptstyle 0\!\!&\!\!\scriptstyle 0\!\!&\!\!\scriptstyle 0\!\!&\!\!\scriptstyle 0\!\!&\!\!\scriptstyle 0\\ \vspace{-.2 cm}
\scriptstyle 0\!\!&\!\!\scriptstyle 0\!\!&\!\!*\!\!&\!\!*\!\!&&\!\!\scriptstyle 0\!\!&\!\!\scriptstyle 0\!\!&\!\!\scriptstyle 0\!\!&\!\!\scriptstyle 0\!\!&\!\!\scriptstyle 0\!\!&\!\!\scriptstyle 0\!\!&\!\!\scriptstyle 0\!\!&\!\!*\!\!&\!\!*\!\!&\!\!\scriptstyle 0\!\!&\!\!*\!\!&\!\!*\!\!&\!\!*\!\!&\!\!\scriptstyle 0\!\!&\!\!\scriptstyle 0\!\!&\!\!\scriptstyle 0\!\!&\!\!\scriptstyle 0\!\!&\!\!\scriptstyle 0\!\!&\!\!\scriptstyle 0\!\!&\!\!\scriptstyle 0\!\!&\!\!\scriptstyle 0\!\!&\!\!\scriptstyle 0\!\!&\!\!\scriptstyle 0\\ \vspace{-.2 cm}
\scriptstyle 0\!\!&\!\!*\!\!&\!\!\scriptstyle 0\!\!&\!\!*\!\!&&\!\!\scriptstyle 0\!\!&\!\!\scriptstyle 0\!\!&\!\!\scriptstyle 0\!\!&\!\!\scriptstyle 0\!\!&\!\!\scriptstyle 0\!\!&\!\!\scriptstyle 0\!\!&\!\!*\!\!&\!\!\scriptstyle 0\!\!&\!\!*\!\!&\!\!*\!\!&\!\!\scriptstyle 0\!\!&\!\!*\!\!&\!\!*\!\!&\!\!\scriptstyle 0\!\!&\!\!\scriptstyle 0\!\!&\!\!\scriptstyle 0\!\!&\!\!\scriptstyle 0\!\!&\!\!\scriptstyle 0\!\!&\!\!\scriptstyle 0\!\!&\!\!\scriptstyle 0\!\!&\!\!\scriptstyle 0\!\!&\!\!\scriptstyle 0\!\!&\!\!\scriptstyle 0\\ \vspace{-.2 cm}
\scriptstyle 0\!\!&\!\!\scriptstyle 0\!\!&\!\!\scriptstyle 0\!\!&\!\!\scriptstyle 0\!\!&&\!\!*\!\!&\!\!*\!\!&\!\!\scriptstyle 0\!\!&\!\!*\!\!&\!\!\scriptstyle 0\!\!&\!\!*\!\!&\!\!\scriptstyle 0\!\!&\!\!\scriptstyle 0\!\!&\!\!\scriptstyle 0\!\!&\!\!\scriptstyle 0\!\!&\!\!\scriptstyle 0\!\!&\!\!\scriptstyle 0\!\!&\!\!\scriptstyle 0\!\!&\!\!\scriptstyle 0\!\!&\!\!\scriptstyle 0\!\!&\!\!\scriptstyle 0\!\!&\!\!*\!\!&\!\!*\!\!&\!\!\scriptstyle 0\!\!&\!\!*\!\!&\!\!\scriptstyle 0\!\!&\!\!*\!\!&\!\!\scriptstyle 0\\ \vspace{-.2 cm}
\scriptstyle 0\!\!&\!\!\scriptstyle 0\!\!&\!\!\scriptstyle 0\!\!&\!\!\scriptstyle 0\!\!&&\!\!*\!\!&\!\!*\!\!&\!\!*\!\!&\!\!*\!\!&\!\!*\!\!&\!\!\scriptstyle 0\!\!&\!\!\scriptstyle 0\!\!&\!\!\scriptstyle 0\!\!&\!\!\scriptstyle 0\!\!&\!\!\scriptstyle 0\!\!&\!\!\scriptstyle 0\!\!&\!\!\scriptstyle 0\!\!&\!\!\scriptstyle 0\!\!&\!\!\scriptstyle 0\!\!&\!\!\scriptstyle 0\!\!&\!\!\scriptstyle 0\!\!&\!\!*\!\!&\!\!*\!\!&\!\!*\!\!&\!\!*\!\!&\!\!*\!\!&\!\!\scriptstyle 0\!\!&\!\!\scriptstyle 0\\ \vspace{-.2 cm}
\scriptstyle 0\!\!&\!\!\scriptstyle 0\!\!&\!\!\scriptstyle 0\!\!&\!\!\scriptstyle 0\!\!&&\!\!*\!\!&\!\!*\!\!&\!\!*\!\!&\!\!\scriptstyle 0\!\!&\!\!*\!\!&\!\!*\!\!&\!\!\scriptstyle 0\!\!&\!\!\scriptstyle 0\!\!&\!\!\scriptstyle 0\!\!&\!\!\scriptstyle 0\!\!&\!\!\scriptstyle 0\!\!&\!\!\scriptstyle 0\!\!&\!\!\scriptstyle 0\!\!&\!\!\scriptstyle 0\!\!&\!\!\scriptstyle 0\!\!&\!\!\scriptstyle 0\!\!&\!\!*\!\!&\!\!*\!\!&\!\!*\!\!&\!\!\scriptstyle 0\!\!&\!\!*\!\!&\!\!*\!\!&\!\!\scriptstyle 0\\ \vspace{-.2 cm}
\scriptstyle 0\!\!&\!\!\scriptstyle 0\!\!&\!\!\scriptstyle 0\!\!&\!\!\scriptstyle 0\!\!&&\!\!\scriptstyle 0\!\!&\!\!\scriptstyle 0\!\!&\!\!*\!\!&\!\!*\!\!&\!\!\scriptstyle 0\!\!&\!\!*\!\!&\!\!\scriptstyle 0\!\!&\!\!\scriptstyle 0\!\!&\!\!\scriptstyle 0\!\!&\!\!\scriptstyle 0\!\!&\!\!\scriptstyle 0\!\!&\!\!\scriptstyle 0\!\!&\!\!\scriptstyle 0\!\!&\!\!\scriptstyle 0\!\!&\!\!\scriptstyle 0\!\!&\!\!\scriptstyle 0\!\!&\!\!\scriptstyle 0\!\!&\!\!\scriptstyle 0\!\!&\!\!*\!\!&\!\!*\!\!&\!\!\scriptstyle 0\!\!&\!\!*\!\!&\!\!\scriptstyle 0\\ \vspace{-.2 cm}
\scriptstyle 0\!\!&\!\!\scriptstyle 0\!\!&\!\!\scriptstyle 0\!\!&\!\!\scriptstyle 0\!\!&&\!\!\scriptstyle 0\!\!&\!\!*\!\!&\!\!*\!\!&\!\!*\!\!&\!\!*\!\!&\!\!\scriptstyle 0\!\!&\!\!\scriptstyle 0\!\!&\!\!\scriptstyle 0\!\!&\!\!\scriptstyle 0\!\!&\!\!\scriptstyle 0\!\!&\!\!\scriptstyle 0\!\!&\!\!\scriptstyle 0\!\!&\!\!\scriptstyle 0\!\!&\!\!\scriptstyle 0\!\!&\!\!\scriptstyle 0\!\!&\!\!\scriptstyle 0\!\!&\!\!\scriptstyle 0\!\!&\!\!*\!\!&\!\!*\!\!&\!\!*\!\!&\!\!*\!\!&\!\!\scriptstyle 0\!\!&\!\!\scriptstyle 0\\ \vspace{-.2 cm}
\scriptstyle 0\!\!&\!\!\scriptstyle 0\!\!&\!\!\scriptstyle 0\!\!&\!\!\scriptstyle 0\!\!&&\!\!\scriptstyle 0\!\!&\!\!*\!\!&\!\!*\!\!&\!\!\scriptstyle 0\!\!&\!\!*\!\!&\!\!*\!\!&\!\!\scriptstyle 0\!\!&\!\!\scriptstyle 0\!\!&\!\!\scriptstyle 0\!\!&\!\!\scriptstyle 0\!\!&\!\!\scriptstyle 0\!\!&\!\!\scriptstyle 0\!\!&\!\!\scriptstyle 0\!\!&\!\!\scriptstyle 0\!\!&\!\!\scriptstyle 0\!\!&\!\!\scriptstyle 0\!\!&\!\!\scriptstyle 0\!\!&\!\!*\!\!&\!\!*\!\!&\!\!\scriptstyle 0\!\!&\!\!*\!\!&\!\!*\!\!&\!\!\scriptstyle 0\\ \vspace{-.2 cm}
\scriptstyle 0\!\!&\!\!\scriptstyle 0\!\!&\!\!\scriptstyle 0\!\!&\!\!\scriptstyle 0\!\!&&\!\!\scriptstyle 0\!\!&\!\!\scriptstyle 0\!\!&\!\!\scriptstyle 0\!\!&\!\!*\!\!&\!\!*\!\!&\!\!*\!\!&\!\!\scriptstyle 0\!\!&\!\!\scriptstyle 0\!\!&\!\!\scriptstyle 0\!\!&\!\!\scriptstyle 0\!\!&\!\!\scriptstyle 0\!\!&\!\!\scriptstyle 0\!\!&\!\!\scriptstyle 0\!\!&\!\!\scriptstyle 0\!\!&\!\!\scriptstyle 0\!\!&\!\!\scriptstyle 0\!\!&\!\!\scriptstyle 0\!\!&\!\!\scriptstyle 0\!\!&\!\!\scriptstyle 0\!\!&\!\!*\!\!&\!\!*\!\!&\!\!*\!\!&\!\!\scriptstyle 0\\ \vspace{-.2 cm}
\scriptstyle 0\!\!&\!\!\scriptstyle 0\!\!&\!\!\scriptstyle 0\!\!&\!\!\scriptstyle 0\!\!&&\!\!\scriptstyle 0\!\!&\!\!\scriptstyle 0\!\!&\!\!\scriptstyle 0\!\!&\!\!\scriptstyle 0\!\!&\!\!\scriptstyle 0\!\!&\!\!\scriptstyle 0\!\!&\!\!\scriptstyle 0\!\!&\!\!\scriptstyle 0\!\!&\!\!\scriptstyle 0\!\!&\!\!\scriptstyle 0\!\!&\!\!\scriptstyle 0\!\!&\!\!\scriptstyle 0\!\!&\!\!\scriptstyle 0\!\!&\!\!\scriptstyle 0\!\!&\!\!\scriptstyle 0\!\!&\!\!\scriptstyle 0\!\!&\!\!*\!\!&\!\!*\!\!&\!\!\scriptstyle 0\!\!&\!\!*\!\!&\!\!\scriptstyle 0\!\!&\!\!*\!\!&\!\!*\\ \vspace{-.2 cm}
\scriptstyle 0\!\!&\!\!\scriptstyle 0\!\!&\!\!\scriptstyle 0\!\!&\!\!\scriptstyle 0\!\!&&\!\!\scriptstyle 0\!\!&\!\!\scriptstyle 0\!\!&\!\!\scriptstyle 0\!\!&\!\!\scriptstyle 0\!\!&\!\!\scriptstyle 0\!\!&\!\!\scriptstyle 0\!\!&\!\!\scriptstyle 0\!\!&\!\!\scriptstyle 0\!\!&\!\!\scriptstyle 0\!\!&\!\!\scriptstyle 0\!\!&\!\!\scriptstyle 0\!\!&\!\!\scriptstyle 0\!\!&\!\!\scriptstyle 0\!\!&\!\!\scriptstyle 0\!\!&\!\!\scriptstyle 0\!\!&\!\!\scriptstyle 0\!\!&\!\!*\!\!&\!\!*\!\!&\!\!*\!\!&\!\!*\!\!&\!\!*\!\!&\!\!\scriptstyle 0\!\!&\!\!*\\ \vspace{-.2 cm}
\scriptstyle 0\!\!&\!\!\scriptstyle 0\!\!&\!\!\scriptstyle 0\!\!&\!\!\scriptstyle 0\!\!&&\!\!\scriptstyle 0\!\!&\!\!\scriptstyle 0\!\!&\!\!\scriptstyle 0\!\!&\!\!\scriptstyle 0\!\!&\!\!\scriptstyle 0\!\!&\!\!\scriptstyle 0\!\!&\!\!\scriptstyle 0\!\!&\!\!\scriptstyle 0\!\!&\!\!\scriptstyle 0\!\!&\!\!\scriptstyle 0\!\!&\!\!\scriptstyle 0\!\!&\!\!\scriptstyle 0\!\!&\!\!\scriptstyle 0\!\!&\!\!\scriptstyle 0\!\!&\!\!\scriptstyle 0\!\!&\!\!\scriptstyle 0\!\!&\!\!*\!\!&\!\!*\!\!&\!\!*\!\!&\!\!\scriptstyle 0\!\!&\!\!*\!\!&\!\!*\!\!&\!\!*\\ \vspace{-.2 cm}
\scriptstyle 0\!\!&\!\!\scriptstyle 0\!\!&\!\!\scriptstyle 0\!\!&\!\!\scriptstyle 0\!\!&&\!\!\scriptstyle 0\!\!&\!\!\scriptstyle 0\!\!&\!\!\scriptstyle 0\!\!&\!\!\scriptstyle 0\!\!&\!\!\scriptstyle 0\!\!&\!\!\scriptstyle 0\!\!&\!\!*\!\!&\!\!*\!\!&\!\!*\!\!&\!\!\scriptstyle 0\!\!&\!\!\scriptstyle 0\!\!&\!\!\scriptstyle 0\!\!&\!\!\scriptstyle 0\!\!&\!\!*\!\!&\!\!*\!\!&\!\!*\!\!&\!\!\scriptstyle 0\!\!&\!\!\scriptstyle 0\!\!&\!\!\scriptstyle 0\!\!&\!\!\scriptstyle 0\!\!&\!\!\scriptstyle 0\!\!&\!\!\scriptstyle 0\!\!&\!\!\scriptstyle 0\\ \vspace{-.2 cm}
\scriptstyle 0\!\!&\!\!\scriptstyle 0\!\!&\!\!\scriptstyle 0\!\!&\!\!\scriptstyle 0\!\!&&\!\!\scriptstyle 0\!\!&\!\!\scriptstyle 0\!\!&\!\!\scriptstyle 0\!\!&\!\!\scriptstyle 0\!\!&\!\!\scriptstyle 0\!\!&\!\!\scriptstyle 0\!\!&\!\!*\!\!&\!\!*\!\!&\!\!*\!\!&\!\!\scriptstyle 0\!\!&\!\!*\!\!&\!\!*\!\!&\!\!\scriptstyle 0\!\!&\!\!*\!\!&\!\!*\!\!&\!\!*\!\!&\!\!\scriptstyle 0\!\!&\!\!\scriptstyle 0\!\!&\!\!\scriptstyle 0\!\!&\!\!\scriptstyle 0\!\!&\!\!\scriptstyle 0\!\!&\!\!\scriptstyle 0\!\!&\!\!\scriptstyle 0\\ \vspace{-.2 cm}
\scriptstyle 0\!\!&\!\!\scriptstyle 0\!\!&\!\!\scriptstyle 0\!\!&\!\!\scriptstyle 0\!\!&&\!\!\scriptstyle 0\!\!&\!\!\scriptstyle 0\!\!&\!\!\scriptstyle 0\!\!&\!\!\scriptstyle 0\!\!&\!\!\scriptstyle 0\!\!&\!\!\scriptstyle 0\!\!&\!\!\scriptstyle 0\!\!&\!\!*\!\!&\!\!*\!\!&\!\!*\!\!&\!\!*\!\!&\!\!*\!\!&\!\!\scriptstyle 0\!\!&\!\!\scriptstyle 0\!\!&\!\!*\!\!&\!\!*\!\!&\!\!\scriptstyle 0\!\!&\!\!\scriptstyle 0\!\!&\!\!\scriptstyle 0\!\!&\!\!\scriptstyle 0\!\!&\!\!\scriptstyle 0\!\!&\!\!\scriptstyle 0\!\!&\!\!\scriptstyle 0\\ \vspace{-.2 cm}
\scriptstyle 0\!\!&\!\!\scriptstyle 0\!\!&\!\!\scriptstyle 0\!\!&\!\!\scriptstyle 0\!\!&&\!\!\scriptstyle 0\!\!&\!\!\scriptstyle 0\!\!&\!\!\scriptstyle 0\!\!&\!\!\scriptstyle 0\!\!&\!\!\scriptstyle 0\!\!&\!\!\scriptstyle 0\!\!&\!\!*\!\!&\!\!*\!\!&\!\!\scriptstyle 0\!\!&\!\!*\!\!&\!\!*\!\!&\!\!\scriptstyle 0\!\!&\!\!*\!\!&\!\!*\!\!&\!\!*\!\!&\!\!\scriptstyle 0\!\!&\!\!\scriptstyle 0\!\!&\!\!\scriptstyle 0\!\!&\!\!\scriptstyle 0\!\!&\!\!\scriptstyle 0\!\!&\!\!\scriptstyle 0\!\!&\!\!\scriptstyle 0\!\!&\!\!\scriptstyle 0\\ \vspace{-.2 cm}
\scriptstyle 0\!\!&\!\!\scriptstyle 0\!\!&\!\!\scriptstyle 0\!\!&\!\!\scriptstyle 0\!\!&&\!\!\scriptstyle 0\!\!&\!\!\scriptstyle 0\!\!&\!\!\scriptstyle 0\!\!&\!\!\scriptstyle 0\!\!&\!\!\scriptstyle 0\!\!&\!\!\scriptstyle 0\!\!&\!\!\scriptstyle 0\!\!&\!\!*\!\!&\!\!*\!\!&\!\!\scriptstyle 0\!\!&\!\!*\!\!&\!\!*\!\!&\!\!*\!\!&\!\!\scriptstyle 0\!\!&\!\!*\!\!&\!\!*\!\!&\!\!\scriptstyle 0\!\!&\!\!\scriptstyle 0\!\!&\!\!\scriptstyle 0\!\!&\!\!\scriptstyle 0\!\!&\!\!\scriptstyle 0\!\!&\!\!\scriptstyle 0\!\!&\!\!\scriptstyle 0\\ \vspace{-.2 cm}
\scriptstyle 0\!\!&\!\!\scriptstyle 0\!\!&\!\!\scriptstyle 0\!\!&\!\!\scriptstyle 0\!\!&&\!\!\scriptstyle 0\!\!&\!\!\scriptstyle 0\!\!&\!\!\scriptstyle 0\!\!&\!\!\scriptstyle 0\!\!&\!\!\scriptstyle 0\!\!&\!\!\scriptstyle 0\!\!&\!\!*\!\!&\!\!\scriptstyle 0\!\!&\!\!*\!\!&\!\!*\!\!&\!\!\scriptstyle 0\!\!&\!\!*\!\!&\!\!*\!\!&\!\!*\!\!&\!\!\scriptstyle 0\!\!&\!\!*\!\!&\!\!\scriptstyle 0\!\!&\!\!\scriptstyle 0\!\!&\!\!\scriptstyle 0\!\!&\!\!\scriptstyle 0\!\!&\!\!\scriptstyle 0\!\!&\!\!\scriptstyle 0\!\!&\!\!\scriptstyle 0\\ \vspace{-.2 cm}
\scriptstyle 0\!\!&\!\!\scriptstyle 0\!\!&\!\!\scriptstyle 0\!\!&\!\!\scriptstyle 0\!\!&&\!\!\scriptstyle 0\!\!&\!\!\scriptstyle 0\!\!&\!\!\scriptstyle 0\!\!&\!\!\scriptstyle 0\!\!&\!\!\scriptstyle 0\!\!&\!\!\scriptstyle 0\!\!&\!\!\scriptstyle 0\!\!&\!\!\scriptstyle 0\!\!&\!\!\scriptstyle 0\!\!&\!\!\scriptstyle 0\!\!&\!\!\scriptstyle 0\!\!&\!\!\scriptstyle 0\!\!&\!\!\scriptstyle 0\!\!&\!\!*\!\!&\!\!*\!\!&\!\!*\!\!&\!\!\scriptstyle 0\!\!&\!\!\scriptstyle 0\!\!&\!\!\scriptstyle 0\!\!&\!\!\scriptstyle 0\!\!&\!\!\scriptstyle 0\!\!&\!\!\scriptstyle 0\!\!&\!\!\scriptstyle 0
\vspace{.2 cm} \end{array} \right] \, . \monendstar 
The four blocks ``ABCD''  are put in evidence.
The square matrices $ \, A \, $ and $ \, D \, $ are of order~4 and 23 respectively. The rectangular matrix $  \, B \, $
has 4 lines and 23 columns and it is the contrary for the matrix $ \, C $. 
%

%%%%%%%%%%%%%%%%%%%%%%%%%%%%%%%%%%%%%%%%%%%%%%%%%%%%%%%%%%%%%%%%%%%%%%%%%%%%%%%%%%%% d3q27 isotherme ordre 1 
\smallskip \monitem 
%%%  
%%% [                                                                rho*uu*dx + rho*vv*dy + rho*ww*dz]
%%% [                          1/3*Phi_ee*dx + 1/3*Phi_xx*dx + Phi_xy*dy + Phi_zx*dz + 2*rho*ll^2/3*dx]
%%% [     Phi_xy*dx + 1/3*Phi_ee*dy + (1/(-6))*Phi_xx*dy + 1/2*Phi_ww*dy + Phi_yz*dz + 2*rho*ll^2/3*dy]
%%% [Phi_zx*dx + Phi_yz*dy + 1/3*Phi_ee*dz + ((-1)/6)*Phi_xx*dz + ((-1)/2)*Phi_ww*dz + 2*rho*ll^2/3*dz]
%%%  
At first order of accuracy, the partial equivalent equations can be written as follows
\moneqstar    \left \{ \begin {array}{l}
\partial_t \rho + \partial_x j_x + \partial_y j_y  + \partial_z j_z  = 0  \\
\partial_t j_x + \partial_x \big( {1\over3} \, \Phi_\varepsilon + {1\over3} \,  \Phi_{xx} +  {2\over3} \,\rho \, \lambda^2 \big)
+  \partial_y \Phi_{xy}  + \partial_z \Phi_{zx} = 0 \\ 
\partial_t j_y + \partial_x \Phi_{xy} + \partial_y \big(  {1\over3} \, \Phi_\varepsilon  - {1\over6} \,  \Phi_{xx} + {1\over2} \,  \Phi_{ww} 
+  {2\over3} \,\rho \, \lambda^2 \big) +  \partial_z \Phi_{yz} = 0 \\
\partial_t j_z + \partial_x \Phi_{zx} + \partial_y \Phi_{yz}  + \partial_z \big(  {1\over3} \, \Phi_\varepsilon
- {1\over6} \,  \Phi_{xx} - {1\over2} \,  \Phi_{ww} +  {2\over3} \,\rho \, \lambda^2 \big) = 0 \, . 
\end {array} \right. \monendstar 
We compare these partial differential equations with the Euler equations (\ref{euler-3d-sans-energie}). 
The two systems are identical if 
the equlibrium functions $\, \Phi_\varepsilon $, $\,  \Phi_{xx} $, $\, \Phi_{ww} $, $\, \Phi_{xy} $, $\, \Phi_{yz} \, $
and $ \, \Phi_{zx} \, $ of the moments of the first  family satisfy the two sets of equations 
\moneqstar    \left \{ \begin {array}{l}
{1\over3} \, \Phi_\varepsilon + {1\over3} \,  \Phi_{xx} +  {2\over3} \,\rho \, \lambda^2   = \rho \, u^2 + p \\
{1\over3} \, \Phi_\varepsilon  - {1\over6} \,  \Phi_{xx} + {1\over2} \,  \Phi_{ww}   +  {2\over3} \,\rho \, \lambda^2 =  \rho \, v^2 + p \\
{1\over3} \, \Phi_\varepsilon - {1\over6} \,  \Phi_{xx} - {1\over2} \,  \Phi_{ww} +  {2\over3} \,\rho \, \lambda^2  = \rho \, w^2 + p 
\end {array} \right. \monendstar
and 
\moneqstar 
\Phi_{xy} = \rho \, u \, v \,,\,\,  \Phi_{yz} =  \rho \, v \, w \,,\,\,  \Phi_{zx} = \rho \, u \, w   . 
\monendstar
With the isothermal hypothesis $ \, p = c_s^2 \, \rho $, 
we recover the  first order isothermal equations with the following choice of equilibria for the second family
of moments:
% 
%%%%  \moneq \label{d3q27-moments-equilibre-2e-ns-isotherme}  \left \{ \begin {array}{l}
\moneqstar 
\Phi_\varepsilon =\rho \, (  | {\bf u} |^2 - 2 \, \lambda^2 + 3  \, c_s^2 ) \,,\,\, 
\Phi_{xx} =  \rho \, ( 2 \, u^2 - v^2 - w^2 )   \,,\,\, \Phi_{ww} =  \rho \, ( v^2 - w^2 )  \, . 
%%%%  \Phi_{xy} =  \rho \,  u \, v  \,,\,\,   \Phi_{yz} =  \rho \,  v \, w \,,\,\,   \Phi_{yz} =  \rho \, w \, u \, . 
%%%%  \end {array} \right. \monend 
\monendstar
%

%%%%%%%%%%%%%%%%%%%%%%%%%%%%%%%%%%%%%%%%%%%%%%%%%%%%%%%%%%%%%%%%%%%%%%%%%%%%%%%%%%%% d3q27 isotherme ordre 2 
\smallskip \monitem 
To recover the  second order equations, a total of  $ \, 3 \times  3 \times  4 \times 3 = 108 \, $ 
linear equations must be solved as in the previous D3Q19 scheme. 
%%%    The viscous terms are present in 3 equations,
%%%    we  have 3 partial derivatives for conservative second order terms in  each equation,
%%%    4 nonconserved variables $\, \rho $, $ \, u $, $ \, v \, $ and $ \, w \, $ and 
%%%    3 partial derivatives $ \, \partial_x $, $ \, \partial_y \, $ and  $ \, \partial_z $.
We have now a total of $ \, 4 \times 7 = 28 \, $ unknowns because the 7 moments of the second family
of nonequilibrium moments occur explicitly in the previous equations and we have 4 partial derivatives
$ \,  \partial_\rho $, $ \, \partial_u $, $ \, \partial_v \, $ and $ \, \partial_w \, $
per moment. Unfortunately, 
this problem has no solution. It is possible to reduce the number of unsolved equations to~56
with an isotropic choice of equilibrium moments for  the third family~:
\moneq \label{d3q27-moments-equilibre-3e-ns-isotherme-isotropic}  \left \{ \! \begin {array}{l} 
\Phi_{q_x} =  \rho \,  u \, \xi_q \,,\,\,
\Phi_{q_y} =  \rho \,  v \, \xi_q \,,\,\,
\Phi_{q_z} =  \rho \,  w \, \xi_q \,,\,\,
\xi_q \equiv 3 \, | {\bf u} |^2  -  2 \, \lambda^2   \\
\Phi_{x_{yz}} = \rho \,  u \, (v^2 - w^2) \,,\,\,  
\Phi_{y_{zx}} = \rho \,  v \, (w^2 - u^2) \,,\,\,  
\Phi_{z_{xy}} = \rho \,  w \, (u^2 - v^2) \,,\,\,   
\Phi_{xyz} = \rho \,  u \, v \, w    \\
\end {array} \right. \monend
and with a specific choice of the sound velocity: $ \, c_s = {{\lambda}\over{\sqrt{3}}} $. 
Moreover, with the following anisotropic choice
\moneq \label{d3q27-moments-equilibre-3e-ns-isotherme-anisotropic}  \left \{ \begin {array}{l}
\Phi_{q_x} =  \rho \,  u \, \xi_q \,,\,\,
\Phi_{q_y} =  \rho \,  v \, \xi_q \,,\,\,
\Phi_{q_z} =  \rho \,  w \, \xi_q \,,\,\,
\xi_q \equiv 3 \, | {\bf u} |^2  -  2 \, \lambda^2   \\
\Phi_{x_{yz}} = \rho \,  u \, \big( (v^2 - w^2) - u^2 \big) \,,\,\,  
\Phi_{y_{zx}} = \rho \,  v \,  \big( (w^2 - u^2)  + v^2 \big) \\   
\Phi_{z_{xy}} = \rho \,  w \,  \big(  (u^2 - v^2)  - w^2 \big) \,,\,\,   
\Phi_{xyz} = \rho \,  u \, v \, w  \, , \\
\end {array} \right. \monend
we have observed that only 44 equations remain unsolved.

\noindent 
With the isotropic choice (\ref{d3q27-moments-equilibre-3e-ns-isotherme-isotropic}) and the
H\'enon relations (\ref{michel-henon}), it is possible to introduce the viscosities
$\, \mu \, $ and $ \, \zeta \, $ thanks to the relation (\ref{visosites-d3q19}).
%%%  ####        munc   = rho*ll^2*sig_x / 3 
%%%  ####        zetanc = 2*rho*ll^2*sig_e / 9 
Then the viscosity tensor is defined with~(\ref{tau-3d}). 
Thus approximative isothermal Navier Stokes equations are solved at second order
\moneqstar   \left \{  \begin {array}{l}
\partial_t \rho + {\textrm {div}} ( \rho \, {\bf u}  ) = {\textrm O}(\Delta x^2  )  \\ 
\partial_t ( \rho \, {\bf u}  ) + {\textrm {div}} \,   ( \rho \, {\bf u} \otimes {\bf u}  ) + \nabla p
-  {\textrm {div}} \, {\bf \tau} +  \Delta t  \, \big(  \sigma_x \,    {\textrm {div}} R^x
+ \sigma_e \,   {\textrm {div}} R^e \big)  =    {\textrm O}(\Delta x^2 ) \,.
\end {array} \right. \monendstar 
The tensor of discrepancy has now two contributions  $ \,  R^x \, $ and $ \, R^e $.
The diagonal terms of the matrix $ \,  R^x \, $ are given by the relations 
\moneqstar \left \{  \begin {array}{l}
R^x_{xx} = u^3 \, \partial_x \rho  -{1\over2} \,  v^3 \, \partial_y \rho - {1\over6} \, \big( u^2 v  + v^3  - u^2 w  -  v^2 w +  v w^2  -2 \, w^3 \big)  \, \partial_z \rho \\
\qquad \qquad + \, \rho \, \big( 3 \, u^2 \, \partial_x u -{1\over3} \, u \, (v - w) \, \partial_z u  -{3\over2} \, v^2 \, \partial_y v
-{1\over6} \, (u^2 +  3 \, v^2  - 2 \, v \, w + w^2 ) \, \partial_z v \\ 
\qquad \qquad \qquad  + {1\over6} \, (u^2 + v^2 - 2  \, v \, w - 6 \,  w^2 )  \, \partial_z w \big) \\
R^x_{yy} = 
-{1\over2}  u^3 \, \partial_x \rho + v^3 \, \partial_y \rho
-{1\over6}\, \big(   u^2 v +   v^3  -  \, u^2 \, w  - v^2 \, w +  v \, w^2 + 2 \, w^3 \, \big) \, \partial_z \rho \\
\qquad \qquad + \, \rho \, \big(  -{3\over2} \, u^2 \, \partial_x u - {1\over3} \, u \, (v-w)  \, \partial_z u + 3 \, v^2 \, \partial_y v
-{1\over6} \, ( u^2  + 3 \, v^2  - 2 \, v \, w +  w^2 ) \, \partial_z v \\ 
\qquad \qquad \qquad  + {1\over6} \, (u^2 + v^2 - 2  \, v \, w - 6 \,  w^2 )  \, \partial_z w \big) \\
R^x_{zz} = -{1\over2} \,  u^3 \, \partial_x \rho  -{1\over2}  v^3 \, \partial_y \rho
+ {1\over3} \, ( u^2 \, v +  v^3   -  u^2 \, w   - v^2 \, w  +  v \, w^2 + 2\, w^3 ) \, \partial_z \rho \\
\qquad \qquad \, + \rho \, \big(    -{3\over2} \, u^2 \, \partial_x u +  {2\over3} \, u \, (v-w) \, \partial_z u  -{3\over2} \,  v^2 \, \partial_y v
+  ( {1\over3} \, u^2 + v^2 - {2\over3} \,  v \,  w + {1\over3} \,  w^2 ) \, \partial_z v \\ 
\qquad \qquad \qquad  - {1\over3} \,  ( u^2 +  v^2 -2 \,  v \, w - 6  \, w^2 ) \, \partial_z w \big) 
\end {array} \right. \monendstar
and the extradiagonal terms of the symmetric tensor $ \,  R^x \, $ follow:
\moneqstar \left \{ \begin {array}{l}
R^x_{xy} =   -{1\over2}  v^3 \, \partial_x \rho  -{1\over2}  u^3 \, \partial_y \rho
- {3\over2} \,\rho \, ( u^2 \, \partial_y u +  v^2 \, \partial_x v ) \\
R^x_{yz} = 
-{1\over2}  u^2 v \, \partial_y \rho -{1\over2}  v^3 \, \partial_y \rho + {1\over2} \, u^2 w \, \partial_y \rho + {1\over2} \, v^2 w \, \partial_y \rho -{1\over2}  v w^2 \, \partial_y \rho -{1\over2}  v^3 \, \partial_z \rho \\
\qquad \qquad + \, \rho \, \big(  u \, (w-v) \, \partial_y u
- {1\over2} \, (  u^2 + 3 \, v^2 - 2 \, v \, w \, \partial_y v +  w^2 ) \, \partial_y v - {3\over2} \,  v^2  \, \partial_z v\\ 
\qquad \qquad \qquad +  {1\over2} \, ( u^2  + v^2 \, \partial_y w - 2 \,  v \, w ) \, \partial_y w \big) \\
R^x_{zx} =  {1\over2} \, ( - u^2 \, v \, -  v^3  + u^2 \, w + v^2 \, w ) \, \partial_x \rho
-{1\over2}  v w^2 \, \partial_x \rho  -{1\over2} \,  u^3 \, \partial_z \rho  
+ \, \rho \, \big(  u \, (w-v) \, \partial_x u  \\
\qquad \qquad -{3\over2} \, u^2 \, \partial_z u -  {1\over2} \, ( u^2 + 3 \,  v^2 - 2 \,  v \, w  + w^2) \, \partial_x v
+ {1\over2} \, ( u^2 + v^2 - 2 \, v \, w ) \, \partial_x w \big) \, . 
\end {array} \right. \monendstar 

\smallskip \noindent 
Finally, the tensor  $ \, R^e \, $ is proportional to the unity tensor:
\moneqstar \left \{ \!\! \begin {array}{l}
R^e_{xx} = R^e_{yy}  =  R^e_{zz}  =  -{1\over3} \, ( u^2 \, v +  v^3 - u^2 \, w -  v^2 \, w +  v \, w^2 -  w^3 ) \, \partial_z \rho
  + \, \rho \, \big( -{2\over3} \,  u \, (v - w) \, \partial_z u  \\ 
\qquad \qquad   -{1\over3} \, ( u^2  + 3 \, v^2 - 2 \, v \, w  +  w^2 ) \, \partial_z \, v 
+ {1\over3} \, ( u^2 +   v^2  -  2 \, v \, w  + 3  \, w^2 ) \, \partial_z w \big) \, . 
\end {array} \right. \monendstar 
With the choices we have done for defining the moments and despite its  nice geometrical  structure,
the D3Q27 scheme is not able to produce equivalent partial equations without artefacts
that involve third order terms relative the  velocity field.

%%%%  \newpage 
%%%%%%%%%%%%%%%%%%%%%%%%%%%%%%%%%%%%%%%%%%%%%%%%%%%%%%%%%%%%%%%%%%%%%%%%%%%%%%%%%%%%
\smallskip \monitem D3Q33
%%%%%%%%%%%%%%%%%%%%%%%%%%%%%%%%%%%%%%%%%%%%%%%%%%%%%%%%%%%%%%%%%%%%%%%%%%%%%%%%%%%%

\begin{figure}    [H]  \centering 
\vspace{-.2 cm} 
\centerline {\includegraphics[height=.55\textwidth]{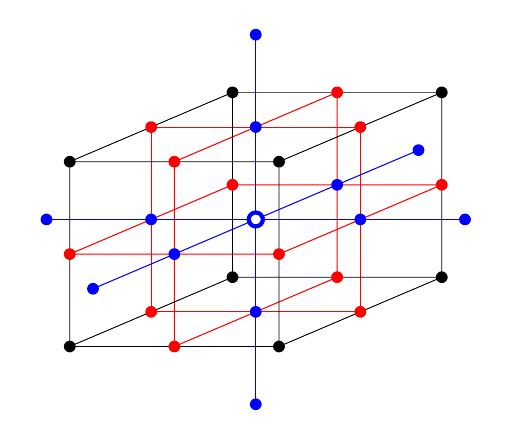}}
\vspace{-.2 cm} 
\caption{Set of discrete velocities for the D3Q33 lattice Boltzmann scheme} 
\label{fig-d3q33} \end{figure}

\noindent
We add six velocities of the type $ \, (2,\,0,\,0) \, $ to the D3Q27 lattice Boltzmann scheme
as proposed in Figure \ref{fig-d3q33}.

%
%%%%%%%%%%%%%%%%%%%%%%%%%%%%%%%%%%%%%%%%%%%%%%%%%%%%%%%%%%%%%%%%%%%%%%%%%%%%%%%%%%% figure d3q33
%%%%%%%%%%%%%%%%%%%%%%%%%%%%%%%%%%%%%%%%%%%%%%%%%%%%%%%%%%%%%%%%%%%%%%%%%%%%%%%%%%%
%
The moments are defined through polynomials with the relation (\ref{moments-polynomes}).
These polynomials are explicited 
in Appendix in the Table~\ref{d3q33-iso-polynomes}.
The four families of moments are described in the Table \ref{d3q33-moments-isotherme}.
The conserved moments and the nonconserved moments of the first family 
 are  similar  to the ones used in the
D3Q19 and D3Q27 schemes. We have now 13 moments of the second family for viscous terms 
and 10 of the third family have no influence on
second order equations. 
%
%%%%%%%%%%%%%%%%%%%%%%%%%%%%%%%%%%%%%%%%%%%%%%%%%%%%%%%%%%%%%%%%%%%%%%%%%%%%%%%%%%% moments d3q33 isotherme     
\begin{table}    [H]  \centering
\begin{tabular}{|c|l|c|c|} \hline
& conserved &  $\rho \,,\,\, j_x \,,\,\, j_y \,,\,\, j_z $ & 4  \\  \hline
1 & fit the Euler equations  &  $ \varepsilon \,,\,\,  xx  \,,\,\,  ww \,,\,\, xy   \,,\,\,  yz \,,\,\,  zx $ & 6 \\ \hline 
2 & fit the viscous terms & $ \! q_x  \,,\,\, q_y  \,,\,\, q_z \,,\,\,  x_{yz}  \,,\,\,  y_{zx}   \,,\,\,   z_{xy} \,,\,\, xyz \,,\,\,
r_x  \,,\,\, r_y  \,,\,\, r_z \,,\,\, t_x  \,,\,\, t_y  \,,\,\, t_z \! $ & 13 \\ \hline 
3 & without influence & $ xx_e  \,,\,\,   ww_e  \,,\,\,   xy_e  \,,\,\,  yz_e  \,,\,\,  zx_e  \,,\,\, h \,,\,\,
 xx_h  \,,\,\,   ww_h  \,,\,\,   h_3 \,,\,\,    h_4 $  & 10  \\ \hline 
\end{tabular} 
\caption{D3Q33 moments for isothermal  Navier Stokes equations}
\label{d3q33-moments-isotherme} \end{table}
%%%%%%%%%%%%%%%%%%%%%%%%%%%%%%%%%%%%%%%%%%%%%%%%%%%%%%%%%%%%%%%%%%%%%%%%%%%%%%%%%%%
%
%
The non  null elements of the $ \, \Lambda \, $ matrix for the D3Q33 scheme are given through the following relation.
The star symbol indicates a non-zero value: 
\moneqstar
%% $ % \begin{equation*}
\Lambda_{D3Q33}^{\rm iso} =   \left[ \begin{array} {cccccccccccccccccccccccccccccccccc}   \vspace{-.25 cm}
\scriptscriptstyle 0\!\!&\!\!\scriptstyle *\!\!&\!\!\scriptstyle *\!\!&\!\!\scriptstyle *\!\!&\!\!&\!\!\scriptscriptstyle 0\!\!&\!\!\scriptscriptstyle 0\!\!&\!\!\scriptscriptstyle 0\!\!&\!\!\scriptscriptstyle 0\!\!&\!\!\scriptscriptstyle 0\!\!&\!\!\scriptscriptstyle 0\!\!&\!\!\scriptscriptstyle 0\!\!&\!\!\scriptscriptstyle 0\!\!&\!\!\scriptscriptstyle 0\!\!&\!\!\scriptscriptstyle 0\!\!&\!\!\scriptscriptstyle 0\!\!&\!\!\scriptscriptstyle 0\!\!&\!\!\scriptscriptstyle 0\!\!&\!\!\scriptscriptstyle 0\!\!&\!\!\scriptscriptstyle 0\!\!&\!\!\scriptscriptstyle 0\!\!&\!\!\scriptscriptstyle 0\!\!&\!\!\scriptscriptstyle 0\!\!&\!\!\scriptscriptstyle 0\!\!&\!\!\scriptscriptstyle 0\!\!&\!\!\scriptscriptstyle 0\!\!&\!\!\scriptscriptstyle 0\!\!&\!\!\scriptscriptstyle 0\!\!&\!\!\scriptscriptstyle 0\!\!&\!\!\scriptscriptstyle 0\!\!&\!\!\scriptscriptstyle 0\!\!&\!\!\scriptscriptstyle 0\!\!&\!\!\scriptscriptstyle 0\!\!&\!\!\scriptscriptstyle 0\\ \vspace{-.25 cm}
\scriptstyle *\!\!&\!\!\scriptscriptstyle 0\!\!&\!\!\scriptscriptstyle 0\!\!&\!\!\scriptscriptstyle 0\!\!&\!\!\!\!&\!\!\scriptstyle *\!\!&\!\!\scriptstyle *\!\!&\!\!\scriptscriptstyle 0\!\!&\!\!\scriptstyle *\!\!&\!\!\scriptscriptstyle 0\!\!&\!\!\scriptstyle *\!\!&\!\!\scriptscriptstyle 0\!\!&\!\!\scriptscriptstyle 0\!\!&\!\!\scriptscriptstyle 0\!\!&\!\!\scriptscriptstyle 0\!\!&\!\!\scriptscriptstyle 0\!\!&\!\!\scriptscriptstyle 0\!\!&\!\!\scriptscriptstyle 0\!\!&\!\!\scriptscriptstyle 0\!\!&\!\!\scriptscriptstyle 0\!\!&\!\!\scriptscriptstyle 0\!\!&\!\!\scriptscriptstyle 0\!\!&\!\!\scriptscriptstyle 0\!\!&\!\!\scriptscriptstyle 0\!\!&\!\!\scriptscriptstyle 0\!\!&\!\!\scriptscriptstyle 0\!\!&\!\!\scriptscriptstyle 0\!\!&\!\!\scriptscriptstyle 0\!\!&\!\!\scriptscriptstyle 0\!\!&\!\!\scriptscriptstyle 0\!\!&\!\!\scriptscriptstyle 0\!\!&\!\!\scriptscriptstyle 0\!\!&\!\!\scriptscriptstyle 0\!\!&\!\!\scriptscriptstyle 0\\
\end{array} \right] \, .   \monendstar
The blocks $ \, A $, $ \, B $, $ \, C \, $ and $ \, D \, $
are of order $ \, 4 \times 4 $,  $ \, 4 \times 29 $,  $ \, 29 \times 4 \, $ and $ \, 29 \times 29 \, $
respectively. 

%%%%%%%%%%%%%%%%%%%%%%%%%%%%%%%%%%%%%%%%%%%%%%%%%%%%%%%%%%%%%%%%%%%%%%%%%%%%%%%%%%%% d3q33 isotherme ordre 1 
\smallskip \monitem 
At first order of accuracy, the partial equivalent equations are easy to produce from
the~$ \, A \, $ and $ \, B \, $ blocks of the operator matrix $\, \Lambda_{D3Q33}^{\rm iso} \, $
and the equilibria: $ \, \Gamma_1 = A \, W + B \, \Phi $. We obtain the first order equivalent equations
%
%%%%%%
%%%%%%     [                                                                   rho*uu*dx + rho*vv*dy + rho*ww*dz]
%%%%%%     [                          1/33*Phi_ee*dx + 1/3*Phi_xx*dx + Phi_xy*dy + Phi_zx*dz + 26*rho*ll^2/33*dx]
%%%%%%     [     Phi_xy*dx + 1/33*Phi_ee*dy + ((-1)/6)*Phi_xx*dy + 1/2*Phi_ww*dy + Phi_yz*dz + 26*rho*ll^2/33*dy]
%%%%%%     [Phi_zx*dx + Phi_yz*dy + 1/33*Phi_ee*dz + ((-1)/6)*Phi_xx*dz + ((-1)/2)*Phi_ww*dz + 26*rho*ll^2/33*dz]
%
\moneqstar    \left \{ \begin {array}{l}
\partial_t \rho + \partial_x j_x + \partial_y j_y  + \partial_z j_z  = 0  \\
\partial_t j_x + \partial_x \big( {1\over33} \, \Phi_\varepsilon + {1\over3} \,  \Phi_{xx} +   {26\over33} \,\rho \, \lambda^2 \big)
+  \partial_y \Phi_{xy}  + \partial_z \Phi_{zx} = 0 \\ 
\partial_t j_y + \partial_x \Phi_{xy} + \partial_y \big(  {1\over33} \, \Phi_\varepsilon  - {1\over6} \,  \Phi_{xx} + {1\over2} \,  \Phi_{ww} 
+  {26\over33} \,\rho \, \lambda^2 \big) +  \partial_z \Phi_{yz} = 0 \\
\partial_t j_z + \partial_x \Phi_{zx} + \partial_y \Phi_{yz}  + \partial_z \big(  {1\over33} \, \Phi_\varepsilon
- {1\over6} \,  \Phi_{xx} - {1\over2} \,  \Phi_{ww} +   {26\over33} \,\rho \, \lambda^2 \big) = 0 \, . 
\end {array} \right. \monendstar 
We compare these  differential equations with the gas dynamics equations (\ref{euler-3d-sans-energie}). 
The two systems are identical when 
the equilibrium functions $\, \Phi_\varepsilon $, $\,  \Phi_{xx} $, $\, \Phi_{ww} $, $\, \Phi_{xy} $, $\, \Phi_{yz} \, $
and $ \, \Phi_{zx} \, $ satisfy the  set of conditions 
\moneqstar    \left \{ \begin {array}{l}
{1\over33} \, \Phi_\varepsilon + {1\over3} \,  \Phi_{xx} +   {26\over33} \,\rho \, \lambda^2   = \rho \, u^2 + p \\
{1\over33} \, \Phi_\varepsilon  - {1\over6} \,  \Phi_{xx} + {1\over2} \,  \Phi_{ww} +  {26\over33} \,\rho \, \lambda^2   =  \rho \, v^2 + p \\
{1\over33} \, \Phi_\varepsilon - {1\over6} \,  \Phi_{xx} - {1\over2} \,  \Phi_{ww} +   {26\over33} \,\rho \, \lambda^2  = \rho \, w^2 + p \\
\Phi_{xy} = \rho \, u \, v \,,\,\,  \Phi_{yz} =  \rho \, v \, w \,,\,\,  \Phi_{zx} = \rho \, u \, w   . 
\end {array} \right. \monendstar
This system is easy to solve and we complete the list of equilibrium moments with the relations
\moneqstar %%% \label{d3q33-moments-equilibre-2e-ns-isotherme}  \left \{ \begin {array}{l}
\Phi_\varepsilon =\rho \, ( 11 \, | {\bf u} |^2 - 26 \, \lambda^2 + 33 \, c_s^2 ) \,,\,\, 
\Phi_{xx} =  \rho \, ( 2 \, u^2 - v^2 - w^2 )   \,,\,\, \Phi_{ww} =  \rho \, ( v^2 - w^2 )  \, . %%  \\
%%% \Phi_{xy} =  \rho \,  u \, v  \,,\,\,   \Phi_{yz} =  \rho \,  v \, w \,,\,\,   \Phi_{yz} =  \rho \, w \, u \, . 
%%% \end {array} \right.
\monendstar 
The  equation of state of the isothermal flow  can be written
$ \, p \equiv c_s^2  \, \rho \, $ 
and the  sound velocity~$ \,  c_s \,  $ is {\it a priori}  not  imposed.

%%%%%%%%%%%%%%%%%%%%%%%%%%%%%%%%%%%%%%%%%%%%%%%%%%%%%%%%%%%%%%%%%%%%%%%%%%%%%%%%%%%%   d3q33 isotherme ordre 2
\smallskip \monitem 
In order to fit the second order partial differential equations as the isothermal Navier Stokes equations, 
the 108 linear equations previously considered 
have to be solved for the partial derivatives,
exactly as the previous D3Q19 and D3Q27 schemes.
We have now 13 moments in the second family and  we have  a total of $ \, 4 \times 13 = 52 \, $ unknowns.
Then a  set of solutions is emerging and a simple quadrature gives necessary
algebraic nonlinear expressions at equilibrium  for the  second  family of moments
$ \,  q_x  $, $ \, q_y  $, $ \, q_z $ $ \, x_{yz} $, $\,  y_{zx} $, $ \, z_{xy} $, $ \, xyz $,
$ \, r_x  $, $ \,  r_y  $, $ \, r_z $, $ \, t_x  $, $ \, t_y  \, $ and $ \, t_z $: 
\moneq \label{d3q33-moments-equilibre-3e-ns-isotherme}  \left \{ \begin {array}{l}
\Phi_{q_x} =  \rho \,  u \, \xi_q \,,\,\,
\Phi_{q_y} =  \rho \,  v \, \xi_q \,,\,\,
\Phi_{q_z} =  \rho \,  w \, \xi_q \,,\,\,
\xi_q \equiv 13 \, | {\bf u} |^2  -  37 \, \lambda^2 + 65 \, c_s^2  \\
\Phi_{x_{yz}} = \rho \,  u \, (v^2 - w^2) \,,\,\,  
\Phi_{y_{zx}} = \rho \,  v \, (w^2 - u^2) \,,\,\,  
\Phi_{z_{xy}} = \rho \,  w \, (u^2 - v^2) \\   
\Phi_{xyz} = \rho \,  u \, v \, w  \\
\Phi_{r_x}  + {38\over13} \,  {{1}\over{\lambda^2}} \, \Phi_{t_x} =  {23\over39} \,  \rho \, u \, \lambda^2 \,
\big( 55 \, \lambda^2 - 111 \, c_s^2 - ( 7 \, u^2 + 195 \, v^2 + 195 \, w^2) \big) \\
\Phi_{r_y}  + {38\over13} \,  {{1}\over{\lambda^2}} \, \Phi_{t_y} =  {23\over39} \,  \rho \, v \, \lambda^2 \,
\big( 55 \, \lambda^2 - 111 \, c_s^2 - ( 7 \, v^2 + 195 \, w^2 + 195 \, u^2) \big) \\
\Phi_{r_z}  + {38\over13} \,  {{1}\over{\lambda^2}} \, \Phi_{t_z} =  {23\over39} \,  \rho \, w \, \lambda^2 \,
\big( 55 \, \lambda^2 - 111 \, c_s^2 - ( 7 \, w^2 + 195 \, u^2 + 195 \, v^2) \big) \, . \\ 
\end {array} \right. \monend
% 
%%%   $ {\color {red} \Phi_{rx} + {38\over13} \,  {{1}\over{\lambda^2}} \, \Phi_{tx} = - {161\over39} \,  \rho \, u^3 \,  \lambda^2
%%%    - {345\over13}\,  \rho \, u \, (v^2 + w^2)  + {1265\over39} \,  \rho \, u \,  \lambda^4 -  {{2553}\over{39}}  \,  \rho \, u \,  \lambda^2 \, c_s^2 }  $
%%%
Obeserve that the relations (\ref{d3q33-moments-equilibre-3e-ns-isotherme})
describe a whole family of possible equilibria. 
The  108 equations  are completely  solved and 
the isothermal Navier Stokes equations are formally  satisfied at second order accuracy as the mesh size tends to zero. 
The relaxation of the moments of the second family as in (\ref{s-ordre-deux}),  introduces the coefficients
$ \, \sigma_x \, $ and $ \,  \sigma_e \, $ with the H\'enon relations (\ref{michel-henon}). 
Then the shear viscosity $ \, \mu \, $ and the bulk viscosity $ \, \zeta \, $
are obtained through the relations
\moneq \label{visosites-d3q33}
\mu =  \rho \, c_s^2 \, \lambda \, \sigma_x \, \Delta x \,,\,\, \zeta = {2\over3} \, \rho \, c_s^2 \,  \lambda \, \sigma_e \, \Delta x \, . 
\monend
These shear and bulk viscosities control the viscous tensor $ \, \tau \, $ defined at the relation (\ref{tau-3d}).
Thus the isothermal Navier Stokes equations take finally the classical form  (\ref{NS3d-isotherme}), {\it id est} 
\moneqstar  \left \{ \begin {array}{l}
\partial_t \rho + {\textrm {div}} ( \rho \, {\bf u}  ) = {\textrm O}(\Delta x^2)   \\  
\partial_t ( \rho \, {\bf u}  ) + {\textrm {div}} \,   ( \rho \, {\bf u} \otimes {\bf u}  ) + \nabla p
-   {\textrm {div}} \, {\bf \tau} =   {\textrm O}(\Delta x^2) \, .  
\end {array} \right. \monendstar  
Observe finally that the equilibrium of 
the  third  family of 10 moments  $ \, xx_e $, $ \,  ww_e   $, $ \, xy_e   $, $ \,  yz_e   $, $ \,  zx_e   $,
$ \, h   $, $ \,    xx_h   $, $ \,   ww_h   $, $ \,   h_3  \, $ and $ \,   h_4 \, $
can be chosen {\it ad libitum} with the point of view of asymptotic analysis. 
 They have no influence on second order equivalent partial differential equations
but their equilibrium values can affect stability. 
  
%%%%%%%%%%%%%%%%%%%%%%%%%%%%%%%%%%%%%%%%%%%%%%%%%%%%%%%%%%%%%%%%%%%%%%%%%%%%%%%%%%%%
\smallskip \monitem D3Q27-2
%%%%%%%%%%%%%%%%%%%%%%%%%%%%%%%%%%%%%%%%%%%%%%%%%%%%%%%%%%%%%%%%%%%%%%%%%%%%%%%%%%%%

\begin{figure}    [H]  \centering 
\vspace{-.7 cm} 
\centerline {\includegraphics[height=.50\textwidth]{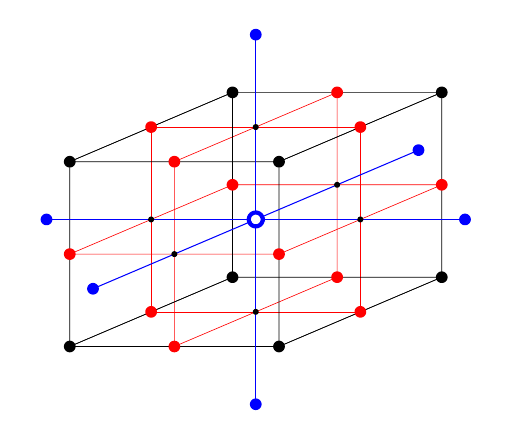}}
\vspace{-.7 cm} 
\caption{Set of discrete velocities for the D3Q27-2 lattice Boltzmann scheme} 
\label{fig-d3q27-2} \end{figure}

\noindent
From the D3Q33 scheme, we drop out the six velocities of the type $ \, (1,\,0,\, 0) \, $ and 
obtain the ``D3Q27-2'' lattice Boltzmann scheme, presented in Figure \ref{fig-d3q27-2}.
%
%%%%%%%%%%%%%%%%%%%%%%%%%%%%%%%%%%%%%%%%%%%%%%%%%%%%%%%%%%%%%%%%%%%%%%%%%%%%%%%%%%% figure d3q27-2
%%%%%%%%%%%%%%%%%%%%%%%%%%%%%%%%%%%%%%%%%%%%%%%%%%%%%%%%%%%%%%%%%%%%%%%%%%%%%%%%%%%
% 
The definition of our moments for  the D3Q27-2 scheme is 
 detailed in the Table \ref{d3q7-2-iso-polynomes}  of the Appendix. Four are conserved,
6 can be tuned in order to fit the first order Euler equations,  10 to fit the viscous terms
and  7 have no  influence on the Navier Stokes equations. They are detailed in the following Table. 
%

%%%%%%%%%%%%%%%%%%%%%%%%%%%%%%%%%%%%%%%%%%%%%%%%%%%%%%%%%%%%%%%%%%%%%%%%%%%%%%%%%%% moments d3q27-2   
\begin{table}    [H]  \centering
\begin{tabular}{|c|l|c|c|} \hline
& conserved &  $\rho \,,\,\, j_x \,,\,\, j_y \,,\,\, j_z $ & 4  \\  \hline
1 & fit the Euler equations  &  $ \varepsilon \,,\,\,  xx  \,,\,\,  ww \,,\,\, xy   \,,\,\,  yz \,,\,\,  zx $ & 6 \\ \hline 
2 & fit the viscous terms & $ q_x  \,,\,\, q_y  \,,\,\, q_z \,,\,\,  x_{yz}  \,,\,\,  y_{zx}   \,,\,\,   z_{xy} \,,\,\, xyz  \,,\,\, r_x  \,,\,\, r_y  \,,\,\, r_z  $  & 10 \\ \hline 
3 & without influence & $ h \,,\,\, xx_e  \,,\,\,   ww_e  \,,\,\,  xy_e  \,,\,\,  yz_e  \,,\,\,  zx_e  \,,\,\,    h_3 $  & 7  \\ \hline 
\end{tabular} 
\caption{D3Q27-2 moments for isothermal  Navier Stokes equations}
\label{d3q27-2-moments-isotherme} \end{table}
%%%%%%%%%%%%%%%%%%%%%%%%%%%%%%%%%%%%%%%%%%%%%%%%%%%%%%%%%%%%%%%%%%%%%%%%%%%%%%%%%%%

%%%%%%%%%%%%%%%%%%%%%%%%%%%%%%%%%%%%%%%%%%%%%%%%%%%%%%%%%%%%%%%%%%%%%%%%%%%%%%%%%%%%   d3q27-2 isotherme ordre 1
\smallskip \monitem 
At first order of accuracy, the partial equivalent equations are easy to produce from
the ABCD approach. We get 
%%%%%%    [                                                                                rho*dx*uu + rho*dy*vv + rho*dz*ww]
%%%%%%    [                                          1/9*dx*Phi_ee + 1/3*dx*Phi_xx + dy*Phi_xy + dz*Phi_zx + 8*rho*ll^2/9*dx]
%%%%%%    [     ((-1)/(-9))*dy*Phi_ee + ((-1)/6)*dy*Phi_xx + 1/2*dy*Phi_ww + dx*Phi_xy + dz*Phi_yz + ((-8*rho*ll^2)/(-9))*dy]
%%%%%%    [((-1)/(-9))*dz*Phi_ee + ((-1)/6)*dz*Phi_xx + ((-1)/2)*dz*Phi_ww + dy*Phi_yz + dx*Phi_zx + ((-8*rho*ll^2)/(-9))*dz]
%
\moneqstar    \left \{ \begin {array}{l}
\partial_t \rho + \partial_x j_x + \partial_y j_y  + \partial_z j_z  = 0  \\
\partial_t j_x + \partial_x \big( {1\over9} \, \Phi_\varepsilon + {1\over3} \,  \Phi_{xx} +   {8\over9} \,\rho \, \lambda^2 \big)
+  \partial_y \Phi_{xy}  + \partial_z \Phi_{zx} = 0 \\ 
\partial_t j_y + \partial_x \Phi_{xy} + \partial_y \big(  {1\over99} \, \Phi_\varepsilon  - {1\over6} \,  \Phi_{xx} + {1\over2} \,  \Phi_{ww} 
+  {8\over9} \,\rho \, \lambda^2 \big) +  \partial_z \Phi_{yz} = 0 \\
\partial_t j_z + \partial_x \Phi_{zx} + \partial_y \Phi_{yz}  + \partial_z \big(  {1\over99} \, \Phi_\varepsilon
- {1\over6} \,  \Phi_{xx} - {1\over2} \,  \Phi_{ww} +   {8\over9} \,\rho \, \lambda^2 \big) = 0 \, . 
\end {array} \right. \monendstar 
We identify these  differential equations with the Euler  equations (\ref{euler-3d-sans-energie}). 
The two systems are identical when  
the equlibrium values $\, \Phi_\varepsilon $, $\,  \Phi_{xx} \, $ and $\, \Phi_{ww} \, $ 
%%%% $\, \Phi_{xy} $, $\, \Phi_{yz} \, $ and $ \, \Phi_{zx} \, $
satisfy the  set of equations 
\moneqstar    \left \{ \begin {array}{l}
{1\over9} \, \Phi_\varepsilon + {1\over3} \,  \Phi_{xx} +   {8\over9} \,\rho \, \lambda^2   = \rho \, u^2 + p \\
{1\over9} \, \Phi_\varepsilon  - {1\over6} \,  \Phi_{xx} + {1\over2} \,  \Phi_{ww} +  {8\over9} \,\rho \, \lambda^2   =  \rho \, v^2 + p \\
{1\over9} \, \Phi_\varepsilon - {1\over6} \,  \Phi_{xx} - {1\over2} \,  \Phi_{ww} +   {8\over9} \,\rho \, \lambda^2  = \rho \, w^2 + p \, . 
%%% \Phi_{xy} = \rho \, u \, v \,,\,\,  \Phi_{yz} =  \rho \, v \, w \,,\,\,  \Phi_{zx} = \rho \, u \, w   . 
\end {array} \right. \monendstar
Then we obtain easily equilibrium values for moments of the first family:
\moneqstar    \left \{ \begin {array}{l}
\Phi_\varepsilon =\rho \, ( 3 \, | {\bf u} |^2 - 8 \, \lambda^2 + 9 \, c_s^2 ) \\ 
\Phi_{xx} =  \rho \, ( 2 \, u^2 - v^2 - w^2 )   \,,\,\, \Phi_{ww} =  \rho \, ( v^2 - w^2 )   \\
\Phi_{xy} = \rho \, u \, v \,,\,\,  \Phi_{yz} =  \rho \, v \, w \,,\,\,  \Phi_{zx} = \rho \, u \, w   . 
\end {array} \right. \monendstar
Observe that the  flow remains isothermal :
\moneqstar
p \equiv c_s^2  \, \rho
\monendstar 
and the sound velocity $ \, c_s \, $ is {\it a priori}  not  imposed.

%%%%%%%%%%%%%%%%%%%%%%%%%%%%%%%%%%%%%%%%%%%%%%%%%%%%%%%%%%%%%%%%%%%%%%%%%%%%%%%%%%%%   d3q27-2 isotherme ordre 2
\smallskip \monitem 
At second order, we must solve 108 equations as previously with the help of partial derivatives
among $ \, \rho $, $ \, u $, $ \, v \, $ and $ \, w \, $ of the 10  moments
$ \, q_x $, $\,  q_y $, $\,  q_z $, $\, x_{yz} $, $\, y_{zx} $, $\, z_{xy} $, $ \, xyz $,   $ \, r_x $, $\,  r_y \, $ and $ \, r_z \, $ 
of the second family, {\it id est} a total of 40 unknowns.
This system admits  a unique solution
and the equilibrium values of the second family  D3Q27-2 moments for isothermal  Navier Stokes equations
satisfy the necessary and sufficient relations:
\moneqstar \label{d3q27-2-moments-equilibre-ns-isotherme-family-3}  \left \{ \begin {array}{l}
\Phi_{q_x} =  \rho \,  u \, \xi_q \,,\,\,
\Phi_{q_y} =  \rho \,  v \, \xi_q \,,\,\,
\Phi_{q_z} =  \rho \,  w \, \xi_q \,,\,\,
\xi_q \equiv | {\bf u} |^2  -  3 \, \lambda^2 + 5 \, c_s^2  \\
\Phi_{x_{yz}} = \rho \,  u \, (v^2 - w^2) \,,\,\,  
\Phi_{y_{zx}} = \rho \,  v \, (w^2 - u^2) \,,\,\,  
\Phi_{z_{xy}} = \rho \,  w \, (u^2 - v^2) \\   
\Phi_{xyz} = \rho \,  u \, v \, w  \\
 \Phi_{r_x}   =  \rho \, u \,   \lambda^2  \, \big( 5 \, \lambda^2 - 9 \, c_s^2 - ( u^2 + 3 \, v^2 + 3 \, w^2)\big)   \\ 
%%% \Phi_{r_x} =  \rho \,  u \, \xi_r \,,\,\,
%%% \xi_r =  \lambda^2  \, (5 \, \lambda^2 - 9 \, c_s^2) -  \lambda^2 \, ( u^3 + 3 \, v^2 + 3 \, w^2)  
 \Phi_{r_y}   =  \rho \, v \, \lambda^2 \,  \big( 5 \, \lambda^2 - 9 \, c_s^2 - ( v^2 + 3 \, w^2 + 3 \, u^2) \big) \\ 
 \Phi_{r_z}   =  \rho \, w \, \lambda^2 \,  \big( 5 \, \lambda^2 - 9 \, c_s^2 - ( w^2 + 3 \, u^2 + 3 \, v^2) \big) \, . 
\end {array} \right. \monendstar
% 
%%% with $ \, | {\bf u} |^2 \equiv u^2 + v^2 + w^2 $.
We remark  that with the previous framework concerning the H\'enon coefficients (\ref{michel-henon}), 
the shear viscosity $ \, \mu \, $ and the bulk viscosity $ \, \zeta \, $
are still obtained through the relations (\ref{visosites-d3q33}). 
The D3Q27-2 lattice Boltzmann scheme allows to recover isothermal Navier Stokes
whether with our choice of moments, it is not the case for the initial D3Q27 scheme!

%%%%%%%%%%%%%%%%%%%%%%%%%%%%%%%%%%%%%%%%%%%%%%%%%%%%%%%%%%%%%%%%%%%%%%%%%%%%%%%  section 6
\bigskip \bigskip    \noindent {\bf \large    6) \quad  Two-dimensional Navier Stokes with conservation of energy}
%%%%%%%%%%%%%%%%%%%%%%%%%%%%%%%%%%%%%%%%%%%%%%%%%%%%%%%%%%%%%%%%%%%%%%%%%%%%%%%%%%%%%%%%%%%%%%%%%%%%%%

%%%%%%%%%%%%%%%%%%%%%%%%%%%%%%%%%%       31 decembre 2018       %%%%%%%%%%%%%%%%%%%%%%%%%
\fancyhead[EC]{\sc{Fran\c{c}ois Dubois and Pierre Lallemand}}
\fancyhead[OC]{\sc{Single lattice Boltzmann distribution for Navier Stokes equations}}
%%%%%%%%%%%%%%%%%%%%%%%%%%%%%%%%%%       31 decembre 2018       %%%%%%%%%%%%%%%%%%%%%%%%%
%%%%%%%%%%%%%%%%%%%%%%%%%%%%%%%%%%%%%%%%%%%%%  jolie numerotation des pages
\fancyfoot[C]{\oldstylenums{\thepage}}
%%%%%%%%%%%%%%%%%%%%%%%%%%%%%%%%%%%%%%%%%%%%%  fin jolie numerotation des pages

\noindent
In this section, we study four schemes in two space dimensions: the D2Q13 presented previously in the
isothermal case and 3 schemes with 17 velocities: D2Q17, D2V17 and D2W17.

%%%%%%%%%%%%%%%%%%%%%%%%%%%%%%%%%%%%%%%%%%%%%%%%%%%%%%%%%%%%%%%%%%%%%%%%%%%%%%%%
\smallskip \monitem D2Q13
%%%%%%%%%%%%%%%%%%%%%%%%%%%%%%%%%%%%%%%%%%%%%%%%%%%%%%%%%%%%%%%%%%%%%%%%%%%%%%%%

\noindent
The velocity set is still given according to Figure \ref{fig-d2q13}. 
The moments are defined with the polynomials
defined in (\ref{polynomes-d2q13-thermique}) in the Appendix.
They are identical to the ones presented   in  (\ref{polynomes-d2q13}).
But the  families of moments differ.
We have now 4 conserved moments (instead of 3),
4 moments are necessary to fit the Euler equations (first family of nonconserved moments),
4 moments are available for the reconstruction of the viscous terms  (second family) and
1 moment has absolutely no influence for the equivalent partial differential equations at second order (third family).
The repartition is explicited in Table~\ref{d2q13-moments-thermiquee}, that can be compared
to the Table \ref{d2q13-moments-isotherme} for isothermal Navier Stokes equations.

%%%%%%%%%%%%%%%%%%%%%%%%%%%%%%%%%%%%%%%%%%%%%%%%%%%%%%%%%%%%%%%%%%%%%%%%%%%%%%%%%%% moments d2q13 ns thermique 
\begin{table}    [H]  \centering
\begin{tabular}{|c|l|c|c|} \hline
  & conserved &  $\rho \,,\,\, j_x \,,\,\, j_y  \,,\,\, \varepsilon $ & 4  \\  \hline
1 & fit the Euler equations  &  $ xx  \,,\,\,  xy   \,,\,\,  q_x  \,,\,\, q_y  $ & 4 \\ \hline 
2 & fit the viscous terms ? & $ r_x \,,\,\,  r_y   \,,\,\, h  \,,\,\, xx_e   $  & 4  \\ \hline  
3 & without influence & $ h_3   $  & 1   \\ \hline 
\end{tabular} 
\caption{The four families of moments for the D2Q13  scheme  for the approximation of the thermal Navier Stokes equations}
\label{d2q13-moments-thermiquee} \end{table}
%%%%%%%%%%%%%%%%%%%%%%%%%%%%%%%%%%%%%%%%%%%%%%%%%%%%%%%%%%%%%%%%%%%%%%%%%%%%%%%%%%%
%
%%%  \qquad \includegraphics[width=.44\textwidth]{d2q13-2018.pdf}
The  operator matrix $ \, \Lambda \, $ is no longer given by the representation
(\ref{d213-Lambda-isotherme}). Despite the coefficients of the matrix are identical,
the ABCD decomposition is different because the number of conserved variables has changed. 
For example, the top-left $ \, A \, $ block is now a $ \, 4 \times 4 \, $ matrix whereas it is a  $ \, 3 \times 3 \, $
matrix for isothermal Navier Stokes equations. 
The operator matrix~$ \, \Lambda \, $ takes now the form 
\moneqstar %%%  \label{d213-Lambda-thermique}
\Lambda_{D2Q13}^{\rm thermal} =  \left[ \begin{array} {ccccccccccccc}
0 \!\!&\!\!  * \partial_x \!\!&\!\!  * \partial_y \!\!&\!\! 0  & %% debut B
 0 \!\!&\!\!  0 \!\!&\!\!  0 \!\!&\!\!  0 \!\!&\!\!
 0 \!\!&\!\!  0 \!\!&\!\!  0 \!\!&\!\!  0 \!\!&\!\!  0 \\
* \partial_x  \!\!&\!\! 0 \!\!&\!\! 0 \!\!&\!\! * \partial_x  & %% debut B
 * \partial_x \!\!&\!\!  * \partial_y \!\!&\!\!  0 \!\!&\!\!  0 \!\!&\!\!  0 \!\!&\!\!  0 \!\!&\!\!  0 \!\!&\!\!
 0 \!\!&\!\!  0 \\
* \partial_y  \!\!&\!\! 0 \!\!&\!\! 0 \!\!&\!\!  * \partial_y & %% debut B
 * \partial_y \!\!&\!\!  * \partial_x \!\!&\!\!  0 \!\!&\!\!  0 \!\!&\!\!  0 \!\!&\!\!  0 \!\!&\!\!
 0 \!\!&\!\!  0 \!\!&\!\!  0 \\
0 \!\!&\!\! * \partial_x \!\!&\!\! * \partial_y \!\!&\!\! 0    & %% debut B
 0 \!\!&\!\!  0 \!\!&\!\!  * \partial_x \!\!&\!\!  * \partial_y \!\!&\!\!  0 \!\!&\!\!  0 \!\!&\!\!  0 \!\!&\!\!
 0 \!\!&\!\!  0 \\ \vspace{-2 mm} \\  %%% fin des matrices A et B
0 \!\!&\!\!  * \partial_x  \!\!&\!\!  * \partial_y \!\!&\!\! 0  & %% debut D
 0 \!\!& \!\!  0 \!\!&\!\!  * \partial_x \!\!&\!\!  * \partial_y \!\!&\!\!  * \partial_x \!\!&\!\!
 * \partial_y \!\!&\!\!  0 \!\!&\!\!  0 \!\!&\!\!  0 \\
0 \!\!&\!\!  * \partial_y  \!\!&\!\!   * \partial_x  \!\!&\!\! 0   & %% debut D
 0 \!\!&\!\!  0 \!\!&\!\!  * \partial_y \!\!&\!\!  * \partial_x \!\!&\!\!  * \partial_y \!\!&\!\!
 * \partial_x \!\!&\!\!  0 \!\!&\!\!  0 \!\!&\!\!  0 \\
0 \!\!&\!\! 0 \!\!&\!\! 0 \!\!&\!\!  * \partial_x  & %% debut D
 * \partial_x \!\!&\!\!  * \partial_y \!\!&\!\!  0 \!\!&\!\!  0 \!\!&\!\! 0 \!\!&\!\! 0 \!\!&\!\!
 * \partial_x \!\!&\!\!  * \partial_x  \!\!&\!\!  0 \\
0\!\!&\!\! 0 \!\!&\!\! 0 \!\!&\!\!  * \partial_y   & %% debut D
 * \partial_y  \!\!&\!\!  * \partial_x \!\!&\!\!  0 \!\!&\!\!  0 \!\!&\!\!  0 \!\!&\!\!  0 \!\!&
\!\!  * \partial_y  \!\!&\!\!  * \partial_y &\!\!  0 \\
0 \!\!&\!\! 0 \!\!&\!\! 0 \!\!&\!\! 0  & %% debut D
 * \partial_x  \!\!&\!\!  * \partial_y \!\!&\!\!  0 \!\!&\!\!  0 \!\!&\!\!  0 \!\!&\!\!  0 \!\!&
\!\!  * \partial_x  \!\!&\!\!  * \partial_x &\!\!  * \partial_x \\
 0 \!\!&\!\!  0 \!\!&\!\!  0 \!\!&\!\!  0  & %% debut D
 * \partial_y  \!\!&\!\!  * \partial_x \!\!&\!\!  0 \!\!&\!\!  0 \!\!&\!\!  0 \!\!&\!\!  0 \!\!&
\!\!  * \partial_y  \!\!&\!\!  * \partial_y &\!\!  * \partial_y \\
 0 \!\!&\!\!  0 \!\!&\!\!  0 \!\!&\!\!  0  & %% debut D
 0 \!\!&\!\!  0 \!\!&\!\!  * \partial_x \!\!&\!\!  * \partial_y \!\!&\!\! * \partial_x \!\!&\!\!
 * \partial_y  \!\!&\!\!  0 \!\!&\!\!  0  \!\!&\!\!  0  \\
 0 \!\!&\!\!  0 \!\!&\!\!  0 \!\!&\!\!  0  & %% debut D
 0 \!\!&\!\!  0 \!\!&\!\!  * \partial_x \!\!&\!\!  * \partial_y \!\!&\!\! * \partial_x \!\!&\!\!
 * \partial_y  \!\!&\!\!  0 \!\!&\!\!  0  \!\!&\!\!  0   \\
 0 \!\!&\!\!  0 \!\!&\!\!  0 \!\!&\!\! 0  & %% debut D
 0 \!\!&\!\!  0 \!\!&\!\!  0 \!\!&\!\!  0 \!\!&\!\!  * \partial_x \!\!&\!\!   * \partial_y \!\!&\!\!
 0 \!\!&\!\!  0 \!\!&\!\! 0
\end{array} \right]  \, . \monendstar
%

%%%%%%%%%%%%%%%%%%%%%%%%%%%%%%%%%%%%%%%%%%%%%%%%%%%%%%%%%%%%%%%%%%%%%%%%%%%%%%%%%%%     d2q13 ns ordre 1  
\smallskip \monitem
The equivalent first order equations of the D2Q13 lattice Boltzmann scheme
\moneqstar   \left \{ \begin {array}{l}
\partial_t \rho + \partial_x j_x + \partial_y j_y = 0 \\
\partial_t j_x + \partial_x ( {14\over13} \, \lambda^2 \, \rho + {1\over26} \varepsilon +   {1\over2}  \Phi_{xx} )  + \partial_y   \Phi_{xy} = 0 \\ 
\partial_t j_y  + \partial_x   \Phi_{xy} + \partial_y ( {14\over13} \, \lambda^2 \, \rho + {1\over26} \varepsilon -   {1\over2}  \Phi_{xx} )  = 0 \\
\partial_t \varepsilon  +  11 \, \lambda^2 \, (\partial_x j_x + \partial_y j_y ) + 13 \, ( \partial_x \Phi_{qx}  + \partial_y \Phi_{qy} ) = 0
\end {array} \right. \monendstar 
must be identical to the Euler equations with the conservation of total energy for two space dimensions 
\moneq \label{euler-2d-avec-energie}   \left \{ \begin {array}{l}
\partial_t \rho + \partial_x j_x + \partial_y j_y  = 0 \\
\partial_t j_x + \partial_x ( \rho \, u^2 + p )  + \partial_y ( \rho \, u \, v ) = 0 \\
\partial_t j_y + \partial_x ( \rho \, u \, v ) +  \partial_y ( \rho \, v^2 + p )  = 0 \\
\partial_t E + \partial_x ( E \, u + p \, u )  + \partial_y ( E \, v + p \, v ) = 0 \, . 
\end {array} \right. \monend
The total energy $ \, E \, $ is given by the relation
\moneqstar
E =  {1\over2} \, \rho \,  |{\bf u} |^2 + \rho \, e \, ,  
\monendstar
with $ \, |{\bf u} |^2 = u^2 + v^2 $. 
Then we identify the expressions under the partial derivatives for the three first relations
and have to solve the system of equations 
\moneqstar  \left \{ \begin {array}{l}
%% \moneq \label{d2q13-euler-energie}    \left \{ \begin {array}{l}
 {14\over13} \, \lambda^2 \, \rho + {1\over26} \varepsilon + {1\over2}  \Phi_{xx}  =  \rho \, u^2 + p \\ 
 {14\over13} \, \lambda^2 \, \rho + {1\over26} \varepsilon - {1\over2}  \Phi_{xx}  =  \rho \, v^2 + p \\ 
 \Phi_{xy}  =  \rho \, u \, v  \,.  %%  \\ 
%%  11 \, \lambda^2 \, \rho  \, u + 13 \, \Phi_{qx} =  E \, u + p \, u  \\ 
%%  11 \, \lambda^2 \, \rho  \, v + 13 \, \Phi_{qy} =  E \, v + p \, v  \, .
\end {array} \right. \monendstar 
We add the  first two equations:
\moneqstar  
 {28\over13} \, \lambda^2 \, \rho + {1\over13} \varepsilon =  \rho \,  |{\bf u} |^2  + 2 \, p = 2 \, E + 2 \, (p - \, \rho \, e)  \,. 
\monendstar
The left hand side of the previous expression is a conserved quantity.
Then it is also the case for the right hand side. The total energy $ \, E \, $ is a conserved quantity.
But it is {\it a priori} not the case for pressure
or volumic internal energy. In consequence, the two last terms must compensate and we  have to take into account 
the constraint 
\moneq  \label{gamma-egale-deux}  
p =  \rho \, e  \, . 
\monend
In other terms, the ratio of specific heats must be equal to 2 : $ \, \gamma = 2 $.
We observe that the energy $ \, \varepsilon \, $ for the D2Q13 lattice Boltzmann scheme
is not exactly equal to the total energy of Physics. We have precisely 
\moneqstar  
\varepsilon = 26 \, E - 28  \, \lambda^2 \, \rho =  13 \, \rho \,  |{\bf u} |^2 + 26 \, \rho \, e - 28 \, \rho \, \lambda^2 \, . 
\monendstar
Then we can write the last relations relative to the energy equation:
\moneqstar \left \{ \begin {array}{l}
11 \,  \lambda^2 \, \rho \, u + 13 \, \Phi_{qx} = 26 \, \big[   \big(  {1\over2} \rho \, (u^2+v^2)  + \rho \, e  \big) \, u + p \, u  \big]  - 28 \, \lambda^2 \, u \\
11 \,  \lambda^2 \, \rho \, v + 13 \, \Phi_{qy} = 26 \,  \big[ \big(  {1\over2} \rho \, (u^2+v^2)  + \rho \, e  \big) \, v + p \, v   \big] - 28 \, \lambda^2 \, v 
\end {array} \right. \monendstar
and we deduce from the previous equations %%% (\ref{d2q13-euler-energie})
\moneqstar \left \{ \begin {array}{l}
\Phi_{xx} = \rho \, (u^2 - v^2) \,,\,\, \Phi_{xy} = \rho \, u \, v \\ 
\Phi_{qx} =  \rho \, u \, \big( u^2 +v^2 + 4 \, e - 3 \, \lambda^2 \big)  \,,\,\, 
\Phi_{qy} =  \rho \, v \, \big( u^2 +v^2 + 4 \, e - 3 \, \lambda^2 \big)  \, .
\end {array} \right. \monendstar 
%%%%  Phiqxnc  = rho*uu*(uu**2 + vv**2 + 4*ee - 3*ll**2) 
%%%%  Phiqync  = rho*vv*(uu**2 + vv**2 + 4*ee - 3*ll**2)
The equilibrium values for the first  family of moments is now completely explicited. 

%%%%%%%%%%%%%%%%%%%%%%%%%%%%%%%%%%%%%%%%%%%%%%%%%%%%%%%%%%%%%%%%%%%%%%%%%%%%%%%%%%%%%%%%%%%%%%   d2q13 ns  ordre 2 
\smallskip \monitem 
For second  order equations, we explicit the vector 
$ \, \Psi_1 = \dd \Phi(W) .  \Gamma_1  - (C \, W + D \, \Phi(W))   \, $ from the previous results and the vector
$ \, \Phi \, $ of moments at equilibrium. The equilibrium functions of the second family 
$ \, \Phi_{rx} $,  $\, \Phi_{ry} $, $ \,  \Phi_{h} \, $ and $ \,  \Phi_{xxe} \, $
%%%
%%%   $ \, \Phi = \big( {\color {blue} \Phi_{xx} ,\, \Phi_{xy} ,\,  \Phi_{qx} \,, \Phi_{qy} } \,,\,
%%%   {\color {red} \Phi_{rx} ,\, \Phi_{ry} ,\,  \Phi_{h}  ,\,  \Phi_{xxe}  } \,,\,    {\color {vertfonce}  \Phi_{h2} } \big)^{\textrm t} $
%%%
have to be determined and the function $ \, \Phi_{h2} \, $ has no influence on the result.
Then we compare the viscous fluxes
\moneqstar
- \Delta t \, \Gamma_2 =  - \Delta t \,  B \, \Sigma \, \Psi_1 
\monendstar
computed with the lattice Boltzmann scheme and the physical viscous  fluxes
given by the relations (\ref{flux-visqueus-ns-thermique-2d}). 
%
%%%%%%  \moneq \label{flux-visqueus-ns-thermique-2d}
%%%%%%  \Phi_{\rm NS}^{2D} = \left( \begin {array}{l}
%%%%%%    \partial_j \sigma_{xj} \equiv \partial_x ( 2 \, \mu \,   \partial_x u
%%%%%%  + (\zeta-\mu)  ( \partial_x u  +   \partial_y v ) )
%%%%%%  + \partial_y (\mu ( \partial_x v  +   \partial_y u ) )  \\
%%%%%%  \partial_j \sigma_{yj}  \equiv \partial_x (\mu (   \partial_x v +    \partial_y u ) )
%%%%%%  + \partial_y (  (\zeta-\mu) (   \partial_x u  +   \partial_y v ) + 2 \, \mu \,   \partial_y v )) \\
%%%%%% \partial_j ( u_i \, \sigma_{ij} ) + {{\gamma}\over{Pr}} \, \big( \partial_x ( \mu \, \partial_x e  ) + \partial_y ( \mu \,  \partial_y e  ) \big) %% ]
%%%%%% \end{array} \right) \, .  \monend
%
We have to solve a total of
3 equations (two for the components of the momentum and one for the energy), 
2 conservation terms $ \, \partial_x [**] \, $ and $ \, \partial_y [**] $ per equation, 
4 nonconserved variables $ \, \rho $, $ \, u $, $ \, v \, $ and internal energy $ \, e $, 
2 partial derivatives  $ \, \partial_x \, $ and $ \, \partial_y  \, $  per variable,
thus $ \, 3 \times 2  \times 4  \times 2 = 48 \, $ equations.
On the other hand, we have 4 moments $ \, \Phi_{rx} $,  $\, \Phi_{ry} $, $ \,  \Phi_{h} \, $ and $ \,  \Phi_{xxe} \, $
(see Table \ref{d2q13-moments-thermiquee})  that generate a total of 16 unknowns
through their partial derivatives relative to $ \, \rho $, $\, u $, $ \, v \, $ and $ \, e $. 
%
%%%%  4 × 3 × 5 × 3 = 180 equations to solve to identify
%%%%  the second order terms of the Navier Stokes equations
%%%%  4 equations for momentum and energy
%%%%  3 conservation terms per equation: ∂x [∗∗], ∂y [∗∗] and ∂z [∗∗]
%%%%  5 nonconserved variables ρ, u, v , w and e
%%%%  3 partial derivatives ∂x , ∂y and ∂z per variable
The corresponding linear system of 48 equations and 16 unknowns 
has no solution. Nevertheless, %%%  with the notation $ \, |{\bf u} |^2 = u^2 + v^2 $,
we have put in evidence  that the following relations 
\moneqstar \left \{ \begin {array}{l}
\Phi_{rx} = \rho \, u \, \lambda^2 \, \big( {31\over6}  \, \lambda^2 - {7\over6} \, (u^2 + 6\, v^2) -  {21\over2} \, e  \big) \\   
\Phi_{ry}  = \rho \, v \, \lambda^2 \, \big( {31\over6}  \, \lambda^2 - {7\over6} \, (6 \, u^2 + v^2) -  {21\over2} \, e  \big) \\ 
 \Phi_{h} \,\, = \rho \, \big(  {77\over2} \,  |{\bf u} |^4 + 308 \, ( |{\bf u} |^2  +e) \, e - 361 \, \lambda^2 \, (  e + |{\bf u} |^2 ) 
 +   140 \, \lambda^4 \big) \\
 \Phi_{xxe} =  \rho \, (u^2-v^2) \, \big( {17\over12} \, |{\bf u} |^2 +  {17\over2} \, e - {65\over12} \, \lambda^2 \big) 
\end {array} \right. \monendstar
reduce to 22 the number of unsolved equations. 
The viscosities can also be explicited:
\moneqstar
\mu = 2 \, \rho \, e \, \sigma_x \, \Delta t \,,\,\, \zeta = \mu \, .
\monendstar
We note  these  important remaining discrepancies
and conclude that the D2Q13 scheme is not formally appropriate for thermal Navier Stokes equations
at second order accuracy  with only one particle distribution.

%%%%%%%    rxrouu =  35*ll^2/6*ee
%%%%%%%    rxrovv =  0
%%%%%%%    rxroee =  0
%%%%%%%    ryrouu =  0
%%%%%%%    ryrovv =  35*ll^2/6*ee
%%%%%%%    ryroee =  ((-21*ll^2)/2)*vv
%%%%%%%    hhrouu =  77*uu^3 + 231*uu*vv^2 + 616*uu*ee + (-154*ll^2)*uu
%%%%%%%    hhrovv =  462*uu^2*vv + 154*vv^3 + 924*vv*ee + (-308*ll^2)*vv
%%%%%%%    hhroee =  308*ee
%%%%%%%    xerouu =  17/3*uu*ee
%%%%%%%    xerovv =  ((-17)/2)*uu^2*vv + ((-17)/6)*vv^3 + (-17)*vv*ee + 17*ll^2/3*vv
%%%%%%%    xeroee =  0
%%%%%%%  eqq05d =  2*rho*ee*sig_x*Delta_t - mu  version 8 juillet
%%%%%%%  eqq09d =  2*rho*ee*sig_x*Delta_t - zeta

%%%%%%%%%%%%%%%%%%%%%%%%%%%%%%%%%%%%%%%%%%%%%%%%%%%%%%%%%%%%%%%%%%%%%%%%%%%%%%%%
\smallskip \monitem D2Q17  %%%  of Qian and Zhou \cite{QZ98}
%%%%%%%%%%%%%%%%%%%%%%%%%%%%%%%%%%%%%%%%%%%%%%%%%%%%%%%%%%%%%%%%%%%%%%%%%%%%%%%%

\begin{figure}    [H]  \centering 
\vspace{-1.2 cm} 
%%%%  \centerline {\includegraphics[height=.55\textwidth]{d2q17-2018.pdf}}
\centerline {\includegraphics[height=.55\textwidth]{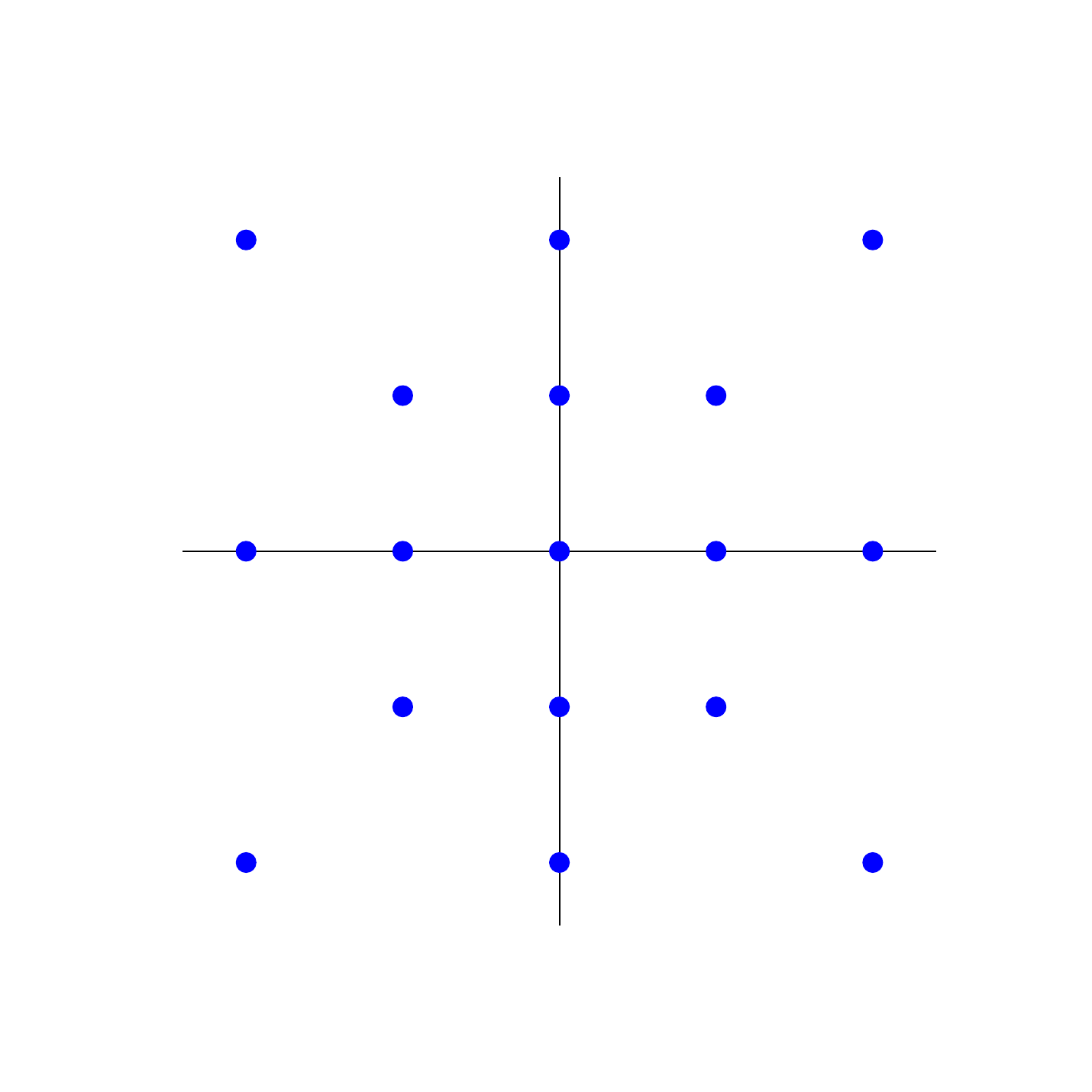}}
\vspace{-.9 cm} 
\caption{Set of discrete velocities for the D2Q17 lattice Boltzmann scheme} 
\label{fig-d2q17} \end{figure}

\noindent
The scheme D2Q17 adds four velocities of the type $ \, (2,_, 2) \, $ to the D2Q13 stencil %% \cite{QZ98},
as presented at Figure \ref{fig-d2q17}.

%%%%%%%%%%%%%%%%%%%%%%%%%%%%%%%%%%%%%%%%%%%%%%%%%%%%%%%%%%%%%%%%%%%%%%%%%%%%%%%%%%% figure d2q17
%%%%%%%%%%%%%%%%%%%%%%%%%%%%%%%%%%%%%%%%%%%%%%%%%%%%%%%%%%%%%%%%%%%%%%%%%%%%%%%%%% 
%
%
The moments are identical to the ones presented  polynomials in  (\ref{d2q17-ns-polynomes}) in the Appendix.
As in the previous D2Q13 scheme, we have four conserved moments 
and four moments are necessary to fit the Euler equations (first family).
But we have now 7 moments for  the reconstruction of the viscous terms  (second family) and
2 moments have  no influence for the equivalent partial differential equations at second order (third  family).
This  repartition is explicited in Table \ref{d2q17-moments-thermique}.
%
%%%%%%%%%%%%%%%%%%%%%%%%%%%%%%%%%%%%%%%%%%%%%%%%%%%%%%%%%%%%%%%%%%%%%%%%%%%%%%%%%%% moments d2q17 ns thermique 
\begin{table}    [H]  \centering
\begin{tabular}{|c|l|c|c|} \hline
0 & conserved &  $\rho \,,\,\, j_x \,,\,\, j_y  \,,\,\, \varepsilon $ & 4  \\  \hline
1 & fit the Euler equations  &  $ xx  \,,\,\,  xy   \,,\,\,  q_x  \,,\,\, q_y  $ & 4 \\ \hline 
2 & fit the viscous terms  & $ r_x \,,\,\,  r_y   \,,\,\, t_x \,,\,\,  t_y   \,,\,\, h  \,,\,\, xx_e   \,,\,\, xy_e  $  & 7  \\ \hline  
3 & without influence & $ h_3 \,,\,\,  h_4  $  & 2   \\ \hline 
\end{tabular} 
\caption{The four families of moments for the D2Q17  scheme  for the approximation of the thermal Navier Stokes equations}
\label{d2q17-moments-thermique} \end{table}
%%%%%%%%%%%%%%%%%%%%%%%%%%%%%%%%%%%%%%%%%%%%%%%%%%%%%%%%%%%%%%%%%%%%%%%%%%%%%%%%%%%
%
The non zero elements of the momentum-velocity operator matrix  $ \, \Lambda \, $
are located in the relation presented below.
The four blocks $\, A \, $(top left), $ \, B \,$(top right), $\, C \, $(bottom left) and
$ \, D \, $(bottom right) are put in evidence:  
\moneq \label{lambda-d2q17-thermique} \Lambda_{D2Q17}^{\rm thermal} =  \left[   \begin{array} {cccccccccccccccccc}
0&*&*&0&&0&0&0&0&0&0&0&0&0&0&0&0&0\\
*&0&0&*&&*&*&0&0&0&0&0&0&0&0&0&0&0\\
*&0&0&*&&*&*&0&0&0&0&0&0&0&0&0&0&0\\
0&*&*&0&&0&0&*&*&0&0&0&0&0&0&0&0&0\\ \vspace{-.3 cm} \\
0&*&*&0&&0&0&*&*&*&*&*&*&0&0&0&0&0\\
0&*&*&0&&0&0&*&*&*&*&*&*&0&0&0&0&0\\
0&0&0&*&&*&*&0&0&0&0&0&0&*&*&*&0&0\\
0&0&0&*&&*&*&0&0&0&0&0&0&*&*&*&0&0\\
0&0&0&0&&*&*&0&0&0&0&0&0&*&*&*&*&0\\
0&0&0&0&&*&*&0&0&0&0&0&0&*&*&*&*&0\\
0&0&0&0&&*&*&0&0&0&0&0&0&0&*&*&*&*\\
0&0&0&0&&*&*&0&0&0&0&0&0&0&*&*&*&*\\
0&0&0&0&&0&0&*&*&*&*&0&0&0&0&0&0&0\\
0&0&0&0&&0&0&*&*&*&*&*&*&0&0&0&0&0\\
0&0&0&0&&0&0&*&*&*&*&*&*&0&0&0&0&0\\
0&0&0&0&&0&0&0&0&*&*&*&*&0&0&0&0&0\\
0&0&0&0&&0&0&0&0&0&0&*&*&0&0&0&0&0 
\end{array} \right]\, .  \monend  
%

%%%%%%%%%%%%%%%%%%%%%%%%%%%%%%%%%%%%%%%%%%%%%%%%%%%%%%%%%%%%%%%%%%%%%%%%%%%%%%%%%%%%%%%%%%%%%%   d2q17 thermique ordre 1 
\smallskip \monitem 
With $ \, \Gamma_1(W) = A \, W + B \, \Phi(W) $,
the equivalent first order equations  
\moneqstar
\partial_t W + \Gamma_1(W) = {\rm O}(\Delta t)
\monendstar 
take the form
\moneqstar   \left \{ \begin {array}{l}
\partial_t \rho + \partial_x j_x + \partial_y j_y = 0 \\
\partial_t j_x + \partial_x ( {30\over17} \, \lambda^2 \, \rho + {1\over34} \varepsilon +   {1\over2}  \Phi_{xx} ) + \partial_y   \Phi_{xy} = 0 \\
\partial_t j_y + \partial_x  \Phi_{xy} + \partial_y ( {30\over17} \, \lambda^2 \, \rho + {1\over34} \varepsilon -   {1\over2}  \Phi_{xx} )   = 0 \\
\partial_t \varepsilon  +  {109\over3} \, \lambda^2 \, (\partial_x j_x + \partial_y j_y ) + {17\over3} \, ( \partial_x \Phi_{qx} + \partial_y \Phi_{qy} ) = 0  \, . 
\end {array} \right. \monendstar 
They have to be confronted to the Euler equations (\ref{euler-2d-avec-energie}). Then the momentum equations
imply
\moneqstar   \left \{ \begin {array}{l}
{30\over17} \, \lambda^2 \, \rho + {1\over34} \varepsilon + {1\over2}  \Phi_{xx}  =  \rho \, u^2 + p \\ 
\Phi_{xy}  =  \rho \, u \, v  \\
{30\over17} \, \lambda^2 \, \rho + {1\over34} \varepsilon - {1\over2}  \Phi_{xx}  =  \rho \, v^2 + p 
\end {array} \right. \monendstar 
and we deduce that 
\moneqstar
\Phi_{xx}  =  \rho \, (u^2 - v^2 ) \, . 
\monendstar 
We observe also that the quantity 
\moneqstar
{60\over17} \, \lambda^2 \, \rho + {1\over17} \, \varepsilon = \rho \, (u^2 + v^2 ) +  2 \, p = 2 \, E + 2 \, (p-\rho \, e) 
\monendstar 
is conserved.
Then the relation (\ref{gamma-egale-deux}) %%  $ \, p = \rho \, e \, $
is again enforced and we must have $ \, \gamma = 2 $. The fourth moment of the D2Q17 lattice Boltzmann scheme $\, \varepsilon \, $
is a linear combination of density $ \, \rho \, $ and of total energy $ \, E $: 
\moneqstar
\varepsilon = 34 \, E - 60  \, \lambda^2 \, \rho = \rho  \, ( 17 \, |{\bf u}|^2 + 34 \, e - 60 \, \lambda^2 ) \, .
\monendstar 
The conservation of the fourth conserved quantity $ \, \varepsilon \, $ is now confronted to
the energy conservation issued from Physics after multiplication by  34, 
with deduction of $ \, 60 \, \lambda^2 \,$  times the mass conservation equation:
\moneqstar
\partial_t (34 \, E - 60  \, \lambda^2 \, \rho) + \partial_x ( \varepsilon \, u + 34 \, p \, u )
+ \partial_y (  \varepsilon \, v + 34 \, p \,  v ) = 0  \,. 
\monendstar 
We deduce the necessary relations 
\moneqstar   %%% \left \{ \begin {array}{l}
{109\over3} \, \lambda^2 \, j_x +  {17\over3} \,  \Phi_{qx} =  \varepsilon \, u + 34 \, p \, u  \,,\,\, 
{109\over3} \, \lambda^2 \, j_y +  {17\over3} \,  \Phi_{qy} =  \varepsilon \, v + 34 \, p \, v  
\monendstar 
and finally
\moneqstar %%%% avril 2021 
\Phi_{qx} = \rho \, u \, \big( 3 \, |{\bf u}|^2 + 12 \, e - 17 \, \lambda^2 \big)  \,,\,\, 
\Phi_{qy} = \rho \, v \, \big( 3 \, |{\bf u}|^2 + 12 \, e - 17 \, \lambda^2 \big) \, . 
\monendstar 
%

%%%%%%%%%%%%%%%%%%%%%%%%%%%%%%%%%%%%%%%%%%%%%%%%%%%%%%%%%%%%%%%%%%%%%%%%%%%%%%%%%%%%%%%%%%%%%%   d2q17 thermique ordre 2 
\smallskip \monitem 
With second order terms, analogously to the D2Q13 scheme, 
we have to solve a set of  48 equations and we have
7 equilibrium unknown functions for the moments
$ \, \Phi_{rx} $, $ \, \Phi_{ry} $, $\, \Phi_{tx} $, $ \, \Phi_{ty} $, $\, \Phi_{h} $, $ \, \Phi_{xxe} \, $ and $\, \Phi_{xye} \, $ 
of the second family, {\it id est} a total of 28 associated partial derivatives that define the unknowns of the algebraic problem. 
If we satisfy the constraint $ \,  \sigma_x  =  \sigma_q $, 
this set of equations admits a family of solutions, characterized by the relations
%
%%%   \moneqstar   \left \{ \begin {array}{l}
%%%   \Phi_{sx} + {31\over2} \, \lambda^2 \, \Phi_{rx} =  {221\over4} \, \rho \, \lambda^6 \, u
%%%   -  {101\over4} \, \rho \, \lambda^4 \, u^3 +   {27\over2} \, \rho \, \lambda^4 \, u \, v^2 -  {249\over4} \, \rho \, \lambda^4 \, u \, e \\
%%%   \Phi_{sy} + {31\over2} \, \lambda^2 \, \Phi_{ry} =  {221\over4} \, \rho \, \lambda^6 \, v
%%%   -  {101\over4} \, \rho \, \lambda^4 \, v^3 +   {27\over2} \, \rho \, \lambda^4 \, u^2 \, v -  {249\over4} \, \rho \, \lambda^4 \, v \, e \\
%%%   \Phi_{h} =  620 \, \rho \, \lambda^4 + {109\over2} \, \rho \, (u^2 + v^2)^2 + 436 \,\rho \, e \, ( u^2 + v^2  +  e )  
%%%   - {969\over2}  \, \rho \, \lambda^2 \,(u^2 + v^2 + 2 e) \\
%%%   \Phi_{xxe} =   -{65\over12} \, \rho \, \lambda^2 \,  (u^2 - v^2) + {17\over12} \, \rho \, (u^4 - v^4) +  {17\over2} \, \rho \, (u^2 - v^2) \, e  \\ 
%%%    \Phi_{xye} =   -{65\over12} \, \rho \, \lambda^2 \, u \, v +  {17\over24} \, \rho \, (u^2 + v^2)\, u \, v \, . 
%%%   \end {array} \right. \monendstar
%
\moneqstar   \left \{ \begin {array}{l} 
 \Phi_{rx} + {2\over31}  {{\Phi_{tx}}\over{\lambda^2}}  =   {1\over62} \, \rho \,  u \, \lambda^2 \big(  221 \, \lambda^2  
-  101 \, u^2 + 54 \, v^2 -  249 \, u \, e \big) \\ 
\Phi_{ry} + {2\over31}  {{\Phi_{ty}}\over{\lambda^2}}  =   {1\over62} \, \rho \,  v \, \lambda^2 \big(  221 \, \lambda^2  
  -  101 \, v^2 + 54 \, u^2 -  249 \, u \, e \big) \\ 
\Phi_{h} =  \rho \, \big( 620 \, \lambda^4 + {109\over2} \, |{\bf u}|^4 + 436 \, e \, ( |{\bf u}|^2 + e )
- {969\over2} \, \lambda^2 \,( |{\bf u}|^2 + 2 e) \big) \\
 \Phi_{xxe} =   \rho \, (u^2 - v^2) \, \big(  -{65\over12} \,  \lambda^2 + {17\over12} \,  |{\bf u}|^2  +  {17\over2} \, e \big) \\ 
 \Phi_{xye} =  \rho \, u \, v \, \big(    -{65\over12} \, \lambda^2  +  {17\over24} \,  |{\bf u}|^2 +  {17\over4} \, e \big)  \, . 
\end {array} \right. \monendstar
%
%%%%  with $ \,  |{\bf u}|^2 = u^2 + v^2 $.
The two viscosities and the Prandtl number satisfy
\moneq \label{mu-zeta-pr-d2q17} 
\mu = \rho \, e \, \sigma_x \, \Delta t  \,,\,\, \zeta = 0 \,,\,\,  Pr = 1  \,.
\monend 
The thermal Navier Stokes equations are formally approachable with the D2Q17 scheme but with a set of physical
constraints: the ratio $ \, \gamma \, $ of the specific heats is equal to 2 and the Prandtl number
is equal to 1. 
%%%  and we must satisfy the constraints $ \,  \sigma_x  =  \sigma_q $. 

%%%%%%%%%%%%%%%%%%%%%%%%%%%%%%%%%%%%%%%%%%%%%%%%%%%%%%%%%%%%%%%%%%%%%%%%%%%%%%%%
\smallskip \monitem  D2V17   %%% of Philippi et al  \cite {PH06}
%%%%%%%%%%%%%%%%%%%%%%%%%%%%%%%%%%%%%%%%%%%%%%%%%%%%%%%%%%%%%%%%%%%%%%%%%%%%%%%%
%

\begin{figure}    [H]  \centering 
\vspace{-1.2 cm} 
%%%%  \centerline {\includegraphics[height=.55\textwidth]{d2q17-2018.pdf}}
\centerline {\includegraphics[height=.55\textwidth]{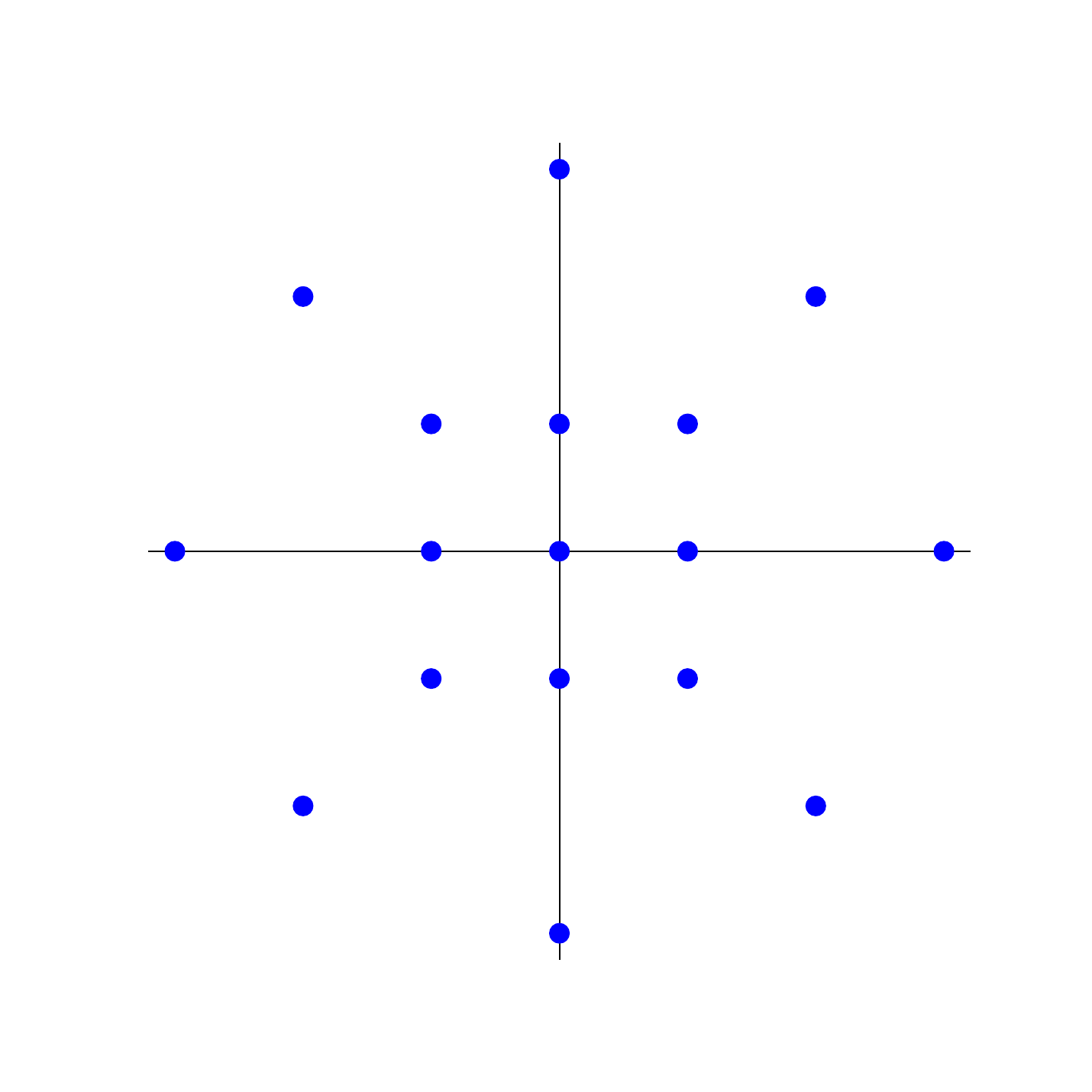}}
\vspace{-.9 cm} 
\caption{Set of discrete velocities for the D2V17 lattice Boltzmann scheme} 
\label{fig-d2v17} \end{figure}

\noindent
The scheme D2V17 eliminates  four velocities of the type $ \, (2, \,0) \, $ from the D2Q13 scheme 
and adds four velocities of the type $ \, (3, \,0) \, $ and four of the type $ \, (2, \,2) \, $
to the stencil, as  proposed at Figure \ref{fig-d2v17}.
%
%%%%%%%%%%%%%%%%%%%%%%%%%%%%%%%%%%%%%%%%%%%%%%%%%%%%%%%%%%%%%%%%%%%%%%%%%%%%%%%%%%% figure d2v17
%%%%%%%%%%%%%%%%%%%%%%%%%%%%%%%%%%%%%%%%%%%%%%%%%%%%%%%%%%%%%%%%%%%%%%%%%%%%%%%%%%%
%
%
The moments are  described with the  polynomials explicited in Table  (\ref{d2v17-ns-polynomes}) of the Appendix.
As in the previous D2Q17 scheme, we have four conserved moments 
and four moments are necessary to fit the Euler equations (first  family), 
7 moments for  the reconstruction of the viscous terms  (second family) and
2 moments have  no influence for the equivalent partial differential equations at second order (third family).
The repartition is explicited in Table~\ref{d2v17-moments-thermique}, similar to the Table~\ref{d2q17-moments-thermique}.
%
%%%%%%%%%%%%%%%%%%%%%%%%%%%%%%%%%%%%%%%%%%%%%%%%%%%%%%%%%%%%%%%%%%%%%%%%%%%%%%%%%%% moments d2v17 ns thermique 
\begin{table}    [H]  \centering
\begin{tabular}{|c|l|c|c|} \hline
  & conserved &  $\rho \,,\,\, j_x \,,\,\, j_y  \,,\,\, \varepsilon $ & 4  \\  \hline
1 & fit the Euler equations  &  $ xx  \,,\,\,  xy   \,,\,\,  q_x  \,,\,\, q_y  $ & 4 \\ \hline 
2 & fit the viscous terms  & $ r_x \,,\,\,  r_y   \,,\,\, t_x \,,\,\,  t_y   \,,\,\, h  \,,\,\, xx_e   \,,\,\, xy_e  $  & 7  \\ \hline  
3 & without influence & $ h_3 \,,\,\,  h_4  $  & 2   \\ \hline 
\end{tabular} 
\caption{The four families of moments for the D2V17  scheme  for the approximation of the thermal Navier Stokes equations}
\label{d2v17-moments-thermique} \end{table}
%%%%%%%%%%%%%%%%%%%%%%%%%%%%%%%%%%%%%%%%%%%%%%%%%%%%%%%%%%%%%%%%%%%%%%%%%%%%%%%%%%%
%
The non zero elements of the $ \, \Lambda \, $ matrix
are located in the relation presented below.
The four blocks $\, A \, $(top left), $ \, B \,$(top right), $\, C \, $(bottom left) and
$ \, D \, $(bottom right) are put in evidence in the following relation: 
\moneqstar \Lambda_{D2V17}^{\rm thermal} = \left[   \begin{array} {cccccccccccccccccc}
0&*&*&0&&0&0&0&0&0&0&0&0&0&0&0&0&0\\
*&0&0&*&&*&*&0&0&0&0&0&0&0&0&0&0&0\\
*&0&0&*&&*&*&0&0&0&0&0&0&0&0&0&0&0\\
0&*&*&0&&0&0&*&*&0&0&0&0&0&0&0&0&0\\
&&&&&&&&&&&&&&&&&\\ \vspace{-.8 cm} \\
0&*&*&0&&0&0&*&*&*&*&*&*&0&0&0&0&0\\
0&*&*&0&&0&0&*&*&*&*&*&*&0&0&0&0&0\\
0&0&0&*&&*&*&0&0&0&0&0&0&*&*&*&0&0\\
0&0&0&*&&*&*&0&0&0&0&0&0&*&*&*&0&0\\
0&0&0&0&&*&*&0&0&0&0&0&0&*&*&*&*&0\\
0&0&0&0&&*&*&0&0&0&0&0&0&*&*&*&*&0\\
0&0&0&0&&*&*&0&0&0&0&0&0&0&*&*&*&*\\
0&0&0&0&&*&*&0&0&0&0&0&0&0&*&*&*&*\\
0&0&0&0&&0&0&*&*&*&*&0&0&0&0&0&0&0\\
0&0&0&0&&0&0&*&*&*&*&*&*&0&0&0&0&0\\
0&0&0&0&&0&0&*&*&*&*&*&*&0&0&0&0&0\\
0&0&0&0&&0&0&0&0&*&*&*&*&0&0&0&0&0\\
0&0&0&0&&0&0&0&0&0&0&*&*&0&0&0&0&0
\end{array} \right] \, . \monendstar 
Its structure is analogous to the one of the operator matrix $ \, \Lambda_{D2Q17}^{\rm thermal} \, $ presented at the relation~(\ref{lambda-d2q17-thermique}).
As previously, all the non zero coefficients are explicited in \cite{Du22git}. 

%%%%%%%%%%%%%%%%%%%%%%%%%%%%%%%%%%%%%%%%%%%%%%%%%%%%%%%%%%%%%%%%%%%%%%%%%%%%%%%%%%%%%%%%%%%%%%   d2v17 thermique ordre 1 
\smallskip \monitem 
The  first order equivalent equations  take the form 
\moneqstar   \left \{ \begin {array}{l}
\partial_t \rho + \partial_x j_x + \partial_y j_y = 0 \\
\partial_t j_x + \partial_x (  {1\over2} \,\Phi_{xx} + {1\over34} \,  \varepsilon + {40\over17} \, \rho \, \lambda^2 )
+ \partial_y   \Phi_{xy} = 0  \\ 
\partial_t j_y + \partial_x  \Phi_{xy}  + \partial_y ( - {1\over2} \,\Phi_{xx} + {1\over34} \,  \varepsilon
+ {40\over17} \, \rho \, \lambda^2 )  -   {1\over2}  \Phi_{xx} )   = 0 \\ 
\partial_t \varepsilon  +  {95\over2} \, \lambda^2 \, (\partial_x j_x + \partial_y j_y )
+ {17\over2} \, ( \partial_x \Phi_{qx}  + \partial_y \Phi_{qy}  ) = 0 \, .
\end {array} \right. \monendstar
They are compared to the Euler equations (\ref{euler-2d-avec-energie}).
The momentum equations imply
\moneqstar   \left \{ \begin {array}{l}
{40\over17} \, \lambda^2 \, \rho + {1\over34} \varepsilon + {1\over2}  \Phi_{xx}  =  \rho \, u^2 + p \\ 
\Phi_{xy}  =  \rho \, u \, v  \\
{40\over17} \, \lambda^2 \, \rho + {1\over34} \varepsilon - {1\over2}  \Phi_{xx}  =  \rho \, v^2 + p \, . 
\end {array} \right. \monendstar 
Then $\,  \Phi_{xx}  =  \rho \, (u^2 - v^2 ) \, $ and
$ \,  {1\over34} \,  \varepsilon + {40\over17} \, \rho \, \lambda^2  =  {1\over2} \, \rho \,  |{\bf u}|^2  + p $.
As developed previously, the quantity 
 $ \,\, {1\over2} \, \rho \,  |{\bf u}|^2 + p \, \, $ must be  conserved  and this is possible only when $ \,  \gamma = 2 $. 
In consequence, $  {1\over34} \,  \varepsilon + {40\over17} \, \rho \, \lambda^2  = E \, $ %% \equiv  {1\over2} \, \rho \, (u^2+v^2) + \rho \, e $
and 
\moneqstar
\varepsilon = 34 \, E - 80  \, \lambda^2 \, \rho  =   \rho \, ( 17 \, | {\bf u} |^2 + 34 \, e - 80 \, \lambda^2 ) \, .
\monendstar 
%
%%% with $ \,  |{\bf u}|^2 = u^2 + v^2 $.  
The comparison between the two forms of conservation of energy 
leads to the equations
\moneqstar
 {95\over2} \, \lambda^2 \, j_x + {17\over2} \, \Phi_{qx}  = \varepsilon \, u + 34 \, p \, u \,,\,\, 
 {95\over2} \, \lambda^2 \, j_y + {17\over2} \, \Phi_{qy}  = \varepsilon \, v + 34 \, p \, v \,  
\monendstar
and we have
\moneqstar
\Phi_{qx}  =  \rho \, u   \, ( 2 \, | {\bf u} |^2  + 8 \, e - 15 \, \lambda^2 ) \,,\,\, 
\Phi_{qy}  =  \rho \, v   \,( 2 \, | {\bf u} |^2  + 8 \, e - 15 \, \lambda^2 )    \, . 
\monendstar 
%

%%%%%%%%%%%%%%%%%%%%%%%%%%%%%%%%%%%%%%%%%%%%%%%%%%%%%%%%%%%%%%%%%%%%%%%%%%%%%%%%%%%%%%%%%%%%%%   d2v17 thermique ordre 2
\smallskip \monitem 
The identification of the second order terms with the dissipative terms of Navier Stokes equations 
leads to a system of 28 equations as presented for the D2Q13 and D2Q17 schemes.
We have 7 moments in the second family, then 28 independent partial derivatives. 
Again, if $ \, \sigma_x  =  \sigma_q $, 
it is possible to operate a reconstruction of the nonlinear equilibrium functions for their corresponding moments: 
\moneqstar   \left \{ \begin {array}{l}
\Phi_{rx} + {5\over3} \, {{\Phi_{sx}}\over{\lambda^2}} =  \rho \, u \, \lambda^2 \,   \big( {1\over9} \,  u^2 
-  {8\over3} \, v^2  -   {7\over3} \, e +  {35\over9} \, \lambda^2 \big)  \\ 
 \Phi_{ry} + {5\over3} \, {{\Phi_{sy}}\over{\lambda^2}} =  \rho \, v \, \lambda^2 \,   \big( {1\over9} \,  v^2 
-  {8\over3} \, u^2  -   {7\over3} \, e +  {35\over9} \, \lambda^2 \big)  \\ 
\Phi_{h} = \rho \, \big( {19\over2} \, |{\bf u}|^4 + 76 \,e \, (  |{\bf u}|^2 +  e)  - 185 \, \lambda^2 \, ( {1\over2} |{\bf u}|^2 + e ) 
+ 100 \, \lambda^4 \big)  \\ 
 \Phi_{xxe} =  \rho \, (u^2 - v^2) \, \big(   {41\over36} \,  |{\bf u}|^2  + {41\over6} \, e -  {365\over36} \,\lambda^2  \big) \\ 
\Phi_{xye} =   \rho \,  u \, v \, \big( {17\over24} \,  |{\bf u}|^2 +  {17\over4} \, e - {65\over12} \, \lambda^2 \big) \, . 
\end {array} \right. \monendstar 
The viscosities and the Prandtl number (\ref{mu-zeta-pr-d2q17}) are analogous to the ones obtained for the  D2Q17 lattice Boltzmann scheme:
$ \, \mu = \rho \, e \, \sigma_x \, \Delta t $, $ \, \zeta = 0 \, $ and $ \, Pr = 1 $.
With our point of view of formal asymptotic analysis, the D2V17 scheme has the same qualities (and defects!)
than the D2Q17 scheme studied in the previous subsection.

%%%%%%%%%%%%%%%%%%%%%%%%%%%%%%%%%%%%%%%%%%%%%%%%%%%%%%%%%%%%%%%%%%%%%%%%%%%%%%%%%%%%%%%%% 
\smallskip \monitem D2W17 %%   of Pierre Lallemand \cite{La10} 
%%%%%%%%%%%%%%%%%%%%%%%%%%%%%%%%%%%%%%%%%%%%%%%%%%%%%%%%%%%%%%%%%%%%%%%%%%%%%%%%%%%%%%%%% 

\begin{figure}    [H]  \centering 
\vspace{-1.4 cm} 
\centerline {\includegraphics[height=.55\textwidth]{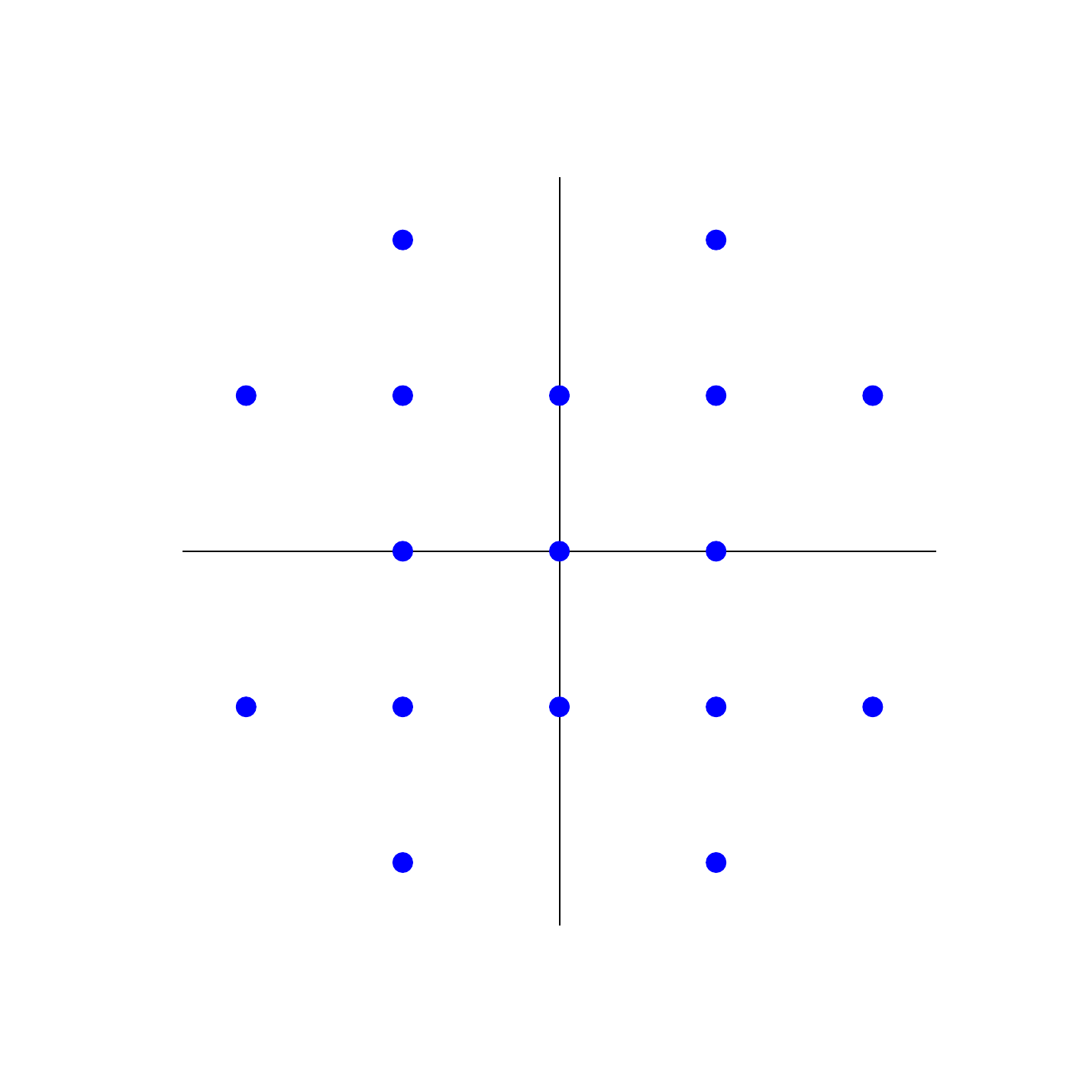}}
%%% \centerline {\includegraphics[height=.55\textwidth]{d2v17-18mars2021-fig.pdf}}
\vspace{-1.1 cm} 
\caption{Set of discrete velocities for the D2W17 lattice Boltzmann scheme} 
\label{fig-d2w17} \end{figure}

\noindent
The scheme D2W17 adds  8 velocities  of the type $ \, (2, \,1) \, $ 
to the D2Q9 scheme,  as a horse moving  at chess play. It is presented at Figure \ref{fig-d2w17}.
%
%%%%%%%%%%%%%%%%%%%%%%%%%%%%%%%%%%%%%%%%%%%%%%%%%%%%%%%%%%%%%%%%%%%%%%%%%%%%%%%%%%% figure d2w17
%%%%%%%%%%%%%%%%%%%%%%%%%%%%%%%%%%%%%%%%%%%%%%%%%%%%%%%%%%%%%%%%%%%%%%%%%%%%%%%%%%%
%
The construction of moments with polynomials is  presented at Table \ref{d2w17-ns-polynomes} of the Appendix. 
It  is a bit different from the choice done previously. 
The two first families of moments are identical to the ones of D2Q17 and D2V17 schemes.
We adopt the nomenclature $ \, r_x $, $ \, r_y $,  $ \, xy_2  $,  $ \, yx_2  $,   $ \, h $, $ \, xx_e  \, $ and  $ \, xy_e  \, $ 
for the  7 moments of the second  family.
The two moments of the  third family are $ \, xx_{xy} \, $ and $ \, h_3 $. 
%
%%%%%%%%%%%%%%%%%%%%%%%%%%%%%%%%%%%%%%%%%%%%%%%%%%%%%%%%%%%%%%%%%%%%%%%%%%%%%%%%%%%   moments d2w17 ns thermique 
\begin{table}    [H]  \centering
\begin{tabular}{|c|l|c|c|} \hline
  & conserved &  $\rho \,,\,\, j_x \,,\,\, j_y  \,,\,\, \varepsilon $ & 4  \\  \hline
1 & fit the Euler equations  &  $ xx  \,,\,\,  xy   \,,\,\,  q_x  \,,\,\, q_y  $ & 4 \\ \hline 
2 & fit the viscous terms  & $ r_x \,,\,\,  r_y   \,,\,\, xy_2  \,,\,\, yx_2  \,,\,\, h  \,,\,\, xx_e   \,,\,\, xy_e  $  & 7  \\ \hline  
3 & without influence & $ xx_{xy}   \,,\,\,  h_3  $  & 2   \\ \hline 
\end{tabular} 
\caption{The four families of moments for the D2W17  scheme  for the approximation of the thermal Navier Stokes equations}
\label{d2w17-moments-thermique} \end{table}
%%%%%%%%%%%%%%%%%%%%%%%%%%%%%%%%%%%%%%%%%%%%%%%%%%%%%%%%%%%%%%%%%%%%%%%%%%%%%%%%%%%
%

%%%%%%%%%%%%%%%%%%%%%%%%%%%%%%%%%%%%%%%%%%%%%%%%%%%%%%%%%%%%%%%%%%%%%%%%%%%%%%%%%%%%%%%%%%%%%%   d2w17 thermique ordre 1
\smallskip \monitem 
The  first order equivalent equations follow the relations 
\moneqstar   \left \{ \begin {array}{l}
\partial_t \rho + \partial_x j_x + \partial_y j_y = 0 \\
\partial_t j_x + \partial_x (  {1\over2} \,\Phi_{xx} + {1\over34} \,  \varepsilon + {26\over17} \, \rho \, \lambda^2 )
+ \partial_y   \Phi_{xy} = 0  \\
\partial_t j_y + \partial_x  \Phi_{xy}  + \partial_y ( - {1\over2} \,\Phi_{xx} + {1\over34} \,  \varepsilon + {26\over17} \, \rho \, \lambda^2 )
-   {1\over2}  \Phi_{xx} )   = 0 \\
\partial_t \varepsilon  +  {259\over13} \, \lambda^2 \, (\partial_x j_x + \partial_y j_y )
+ {17\over13} \, ( \partial_x \Phi_{qx}  + \partial_y \Phi_{qy}  ) = 0 \, .
\end {array} \right. \monendstar
They are compared to the Euler equations (\ref{euler-2d-avec-energie}).
The momentum equations induce 
\moneqstar   \left \{ \begin {array}{l}
{26\over17} \, \lambda^2 \, \rho + {1\over34} \varepsilon + {1\over2}  \Phi_{xx}  =  \rho \, u^2 + p \\ 
\Phi_{xy}  =  \rho \, u \, v  \\
{26\over17} \, \lambda^2 \, \rho + {1\over34} \varepsilon - {1\over2}  \Phi_{xx}  =  \rho \, v^2 + p \, . 
\end {array} \right. \monendstar 
Then $\,  \Phi_{xx}  =  \rho \, (u^2 - v^2 ) \, $ and
$ \,  {1\over34} \,  \varepsilon + {26\over17} \, \rho \, \lambda^2  =  {1\over2} \, \rho \, (u^2+v^2)  + p $.
In consequence, the quantity 
$ \,\, {1\over2} \, \rho \,  |{\bf u}|^2 + p \, \, $ must be  conserved and this enforces again the relations
$ \, p \equiv \rho\, e \, $    and $ \,  \gamma = 2 $. 
We have now  $  {1\over34} \,  \varepsilon + {26\over17} \, \rho \, \lambda^2  = E \, $ %% \equiv  {1\over2} \, \rho \, (u^2+v^2) + \rho \, e $
and 
\moneqstar 
\varepsilon = 34 \, E - 52  \, \lambda^2 \, \rho =  \rho \, ( 17 \, | {\bf u} |^2 + 34 \, e - 52 \, \lambda^2 ) \, . 
\monendstar
The conservation of energy equation is compared to 34 times the physical energy conservation minus $ \, 52 \, \lambda^2 \,$  times the mass conservation,
{\it id est}
$ \,  \partial_t  \varepsilon  + \partial_x (  \varepsilon \, u + 34 \, p \, u )  + \partial_y ( \varepsilon \, v + 34 \,  p \, v ) = 0 $. 
In consequence, we have
$\,  {259\over13} \, \lambda^2 \, j_x + {17\over13} \, \Phi_{qx}  = \varepsilon \, u + 34 \, p \, u  \, $ and
$ \,   {259\over13} \, \lambda^2 \, j_y + {17\over13} \, \Phi_{qy}  = \varepsilon \, v + 34 \, p \, v   $. 
Then 
\moneqstar   \left \{ \begin {array}{l}
\Phi_{qx}  = {13\over17} \, \varepsilon \, u + 26 \, p \, u -  {259\over17} \, \lambda^2 \, j_x =
 \rho \, u \, ( 13 \, | {\bf u} |^2  + 52 \, e - 55 \, \lambda^2 )  \\ 
\Phi_{qy}  = {13\over17} \, \varepsilon \, v + 26 \, p \, v -  {259\over17} \, \lambda^2 \, j_y = 
\rho \, v \, ( 13 \, | {\bf u} |^2  + 52 \, e - 55 \, \lambda^2 ) \, . 
\end {array} \right. \monendstar 
%

%%%%%%%%%%%%%%%%%%%%%%%%%%%%%%%%%%%%%%%%%%%%%%%%%%%%%%%%%%%%%%%%%%%%%%%%%%%%%%%%%%%%%%%%%%%%%%   d2w17 thermique ordre 2
\smallskip \monitem 
The second order equations are treated similarly to the D2Q17 and D2V17 lattice Boltzmann schemes.
We must constrain the scheme with the relation $ \,  \sigma_x  =  \sigma_q \, $
and we have necessary relations between the 7 moments of the second family:
\moneqstar   \left \{ \begin {array}{l}
\Phi_{rx} + {171\over2} \, \Phi_{xy2} =  \rho \, u \, \lambda^2 \, \big( {85\over4} \,  u^2 
-  50  \, v^2  +   {55\over4} \, e +  {35\over4} \, \lambda^2 \big) \\
\Phi_{ry} + {171\over2} \, \Phi_{yx2} =   \rho \, u \, \lambda^2 \, \big( -50 \, u^2 +  {85\over4} \, v^2 
+  {55\over4}  \, e +  {35\over4} \, \lambda^2 \big) \\
\Phi_{h} =  \rho \,  \big( {259\over2} \, |{\bf u}|^4 + 1036 \, (  |{\bf u}|^2 +  e) \, e
- 1543 \, \lambda^2 \, (  {1\over2} \, |{\bf u}|^2 + e ) + 684 \, \lambda^4 \big) \\
\Phi_{xxe} =  \rho \, (u^2 - v^2) \, \big(   {19\over12} \,  |{\bf u}|^2  + {19\over2} \, e  -  {91\over12} \, \lambda^2 \big) \\
\Phi_{xye} =   \rho \, u \, v \, \big(  {3\over2} \,  |{\bf u}|^2 +  9 \, e - 7 \, \lambda^2 \big) \, . 
\end {array} \right. \monendstar 
Finally, we have the relation  (\ref{mu-zeta-pr-d2q17})
for the viscosities and the Prandtl number, as for the previous D2Q17 and D2V17 schemes. 
%####         munc   = rho*ee*sig_x
%####         zetanc = 0
%####         sigqnc = sig_x
%####         gspnc  = 2
%
%####      \moneqstar 
%####      \mu = \rho \, e \, \sigma_x \, \Delta t  \,,\,\,  \zeta = 0 \,,\,\,  Pr = 1
%####      \monendstar 
%
%%%  $ \, \mu = \rho \, e \, \sigma_x \, \Delta t $, $\, \zeta = 0 \, $ and $ \, Pr = 1 \, $

\smallskip \noindent
At two space dimensions, the three lattice Boltzmann schemes D2Q17, D2V17 and D2W17 give completely analogous results
for the thermal Navier Stokes equations.
The precision is of second order formal accuracy but  physical constraints must be enforced: 
the ratio $ \, \gamma \, $ of the specific heats is equal to 2 and the Prandtl number
is equal to 1.

%%%%%%%%%%%%%%%%%%%%%%%%%%%%%%%%%%%%%%%%%%%%%%%%%%%%%%%%%%%%%%%%%%%%%%%%%%%%%%%  section 7
\bigskip \bigskip    \noindent {\bf \large    7) \quad  Three-dimensional Navier Stokes with conservation of energy}
%%%%%%%%%%%%%%%%%%%%%%%%%%%%%%%%%%%%%%%%%%%%%%%%%%%%%%%%%%%%%%%%%%%%%%%%%%%%%%%%%%%%%%%%%%%%%%%%%%%%%%

\noindent
Since the two classic schemes D3Q19 and D3Q27 are not able to recover exactly  isothermal Navier Stokes equations
at second order accuracy with our framework,
in consequence, they are not considered in this section. 
%%%%%%  , we do not consider these two schemes in this
%%%%%%  section. We describe the two schemes D3Q33 and D3Q27-2.
%%%%%%  In three space dimensions, observe first that the popular D3Q19 and D3Q27 schemes, not exact at second order
%%%%%%  for isothermal Navier Stokes, are excluded of this study.
Therefore, we consider two schemes in three space dimensions: D3Q33 and D3Q27-2, already considered
in the isothermal case.

%%%%%%%%%%%%%%%%%%%%%%%%%%%%%%%%%%%%%%%%%%%%%%%%%%%%%%%%%%%%%%%%%%%%%%%%%%%%%%%%%%%%%%%%%%%%%%%%%%%%%%
\smallskip \monitem  D3Q33
%%%%%%%%%%%%%%%%%%%%%%%%%%%%%%%%%%%%%%%%%%%%%%%%%%%%%%%%%%%%%%%%%%%%%%%%%%%%%%%%%%%%%%%%%%%%%%%%%%%%%%

\noindent
The stencil of velocities has been presented previously at Figure \ref{fig-d3q33}. 
The moments are  explicited through the polynomials presented at the Table \ref{d3q33-ns-polynomes}
in the Appendix.
They are identical to the ones of the Table \ref{d3q33-iso-polynomes}
but their repartition into four families is changed.
We have  now 5 conserved moments and
a total of 8 moments constitute the first family of nonconserved moments that allows to treat first order terms.
Moreover, 16 moments are available for the second family that influences  the equivalent
partial differential equations at second order accuracy.
The third  family is composed by 4 moments which have no impact on second order terms.
This repartition is summarized in Table \ref{d3q33-moments-thermique}. 
%
%%%%%%%%%%%%%%%%%%%%%%%%%%%%%%%%%%%%%%%%%%%%%%%%%%%%%%%%%%%%%%%%%%%%%%%%%%%%%%%%%%% moments d3q33  thermique       
\begin{table}    [H]  \centering
\begin{tabular}{|c|l|c|c|} \hline
 & conserved &  $\rho \,,\,\, j_x \,,\,\, j_y \,,\,\, j_z \,,\,\, \varepsilon  $ & 5  \\  \hline
1 & fit the Euler equations  &  $ xx  \,,\,\,  ww \,,\,\, xy   \,,\,\,  yz \,,\,\,  zx  \,,\,\, q_x  \,,\,\, q_y  \,,\,\, q_z $ & 8 \\ \hline
2 & fit the viscous terms & $ x_{yz}  \,,\,\,  y_{zx}   \,,\,\,   z_{xy} \,,\,\, xyz \,,\,\, r_x  \,,\,\, r_y   \,,\,\, r_z   \,,\,\,  t_x  \,,\,\, t_y  \,,\,\, t_z  $  & \\
&  & $ xx_e  \,,\,\,   ww_e  \,,\,\,   xy_e  \,,\,\,  yz_e  \,,\,\,  zx_e  \,,\,\, h  $  & 16 \\ \hline 
3 & without influence & $       xx_h  \,,\,\,   ww_h  \,,\,\,   h_3 \,,\,\,    h_4 $  & 4  \\ \hline 
\end{tabular} 
\caption{The four families of D3Q33 moments for thermal  Navier Stokes equations}
\label{d3q33-moments-thermique} \end{table}
%%%%%%%%%%%%%%%%%%%%%%%%%%%%%%%%%%%%%%%%%%%%%%%%%%%%%%%%%%%%%%%%%%%%%%%%%%%%%%%%%%%
%
The moments of the two first families constitute exactly the set of Grad's 13 moments model \cite{Gr58,ST03}.
In this model, the density, the 3 components of momentum, the 6 components of second order moments
and the 3 components of  heat flux vector are conserved quantities.
In the present study, only 5 moments are conserved and the other 8 of order 2 and 3 are function of the conserved ones
to ensure  the consistency with the Euler equations.

%%%%%%%%%%%%%%%%%%%%%%%%%%%%%%%%%%%%%%%%%%%%%%%%%%%%%%%%%%%%%%%%%%%%%%%%%%%%%%%%%%%%%%%%%%%%%%   d3q33 thermique ordre 1 
\smallskip \monitem 
At first order, the reference equations are the Euler equations of gas dynamics. %%  (see {\it e.g.} \cite{GHP01}). 
The pressure is a function of density and energy; the total volumic energy $ \, E \, $ is still given by
the relation $ \, E = {1\over2} \, \rho \, |{\bf u}|^2 + \rho \, e \, $ 
with $ \,  |{\bf u}|^2 = u^2 + v^2 + w^2 $. We set also
$ \, j_x = \rho\, u $, $ \, j_y = \rho\, v \, $ and $ \, j_z = \rho\, w  \, $ and we have 
\moneq \label{euler-3d-avec-energie}   \left \{ \begin {array}{l}
\partial_t \rho + \partial_x j_x + \partial_y j_y  + \partial_z j_z  = 0   \\
\partial_t j_x + \partial_x ( \rho \, u^2 + p )  + \partial_y ( \rho \, u \, v )  + \partial_z ( \rho \, u \, w) = 0  \\
\partial_t j_y + \partial_x ( \rho \, u \, v ) +  \partial_y ( \rho \, v^2 + p ) + \partial_z ( \rho \, v \, w ) =  0 \\
\partial_t j_z + \partial_x ( \rho \, u \, w ) +  \partial_y ( \rho \, v \, w ) + \partial_z ( \rho \, w^2 + p ) =  0 \\
\partial_t E + \partial_x ( E \, u + p \, u )  + \partial_y ( E \, v + p \, v )  + \partial_z ( E \, w + p \, w ) = 0 \,. 
\end {array} \right. \monend
The equivalent partial differential equations at first order follow the expressions
%
%%%   $ flimxx = 1/3*Phi_xx + 1/33*eps + 26*rho*ll^2/33 $
%%%   $ flimyy = ((-1)/6)*Phi_xx + 1/2*Phi_ww + 1/33*eps + 26*rho*ll^2/33 $
%%%   $ flimzz = ((-1)/6)*Phi_xx + ((-1)/2)*Phi_ww + 1/33*eps + 26*rho*ll^2/33 $ 
%%%   $ 11/13*Phi_qx*dx + 69*rho*ll^2/13*uu*dx + 11/13*Phi_qy*dy + 69*rho*ll^2/13*vv*dy + 11/13*Phi_qz*dz + 69*rho*ll^2/13*ww*dz] $ 
%
\moneqstar   \left \{ \begin {array}{l}
\partial_t \rho + \partial_x j_x + \partial_y j_y + \partial_z j_z = 0  \\ 
\partial_t j_x + \partial_x ( %%%  {1\over3} \, \rho \, |{\bf u}|^2 + {1\over3} \, \Phi_{xx} + {2\over3} \, \rho \, e )
{1\over33} \, \varepsilon +  {1\over3} \, \Phi_{xx} + {26\over33} \,\rho \, \lambda^2) 
+ \partial_y  \Phi_{xy} + \partial_z  \Phi_{zx} = 0 \\
\partial_t j_y  + \partial_x \Phi_{xy} + \partial_y ( %%% {1\over3} \, \rho \, |{\bf u}|^2 - {1\over6} \, \Phi_{xx} + {1\over2} \, \Phi_{ww}+ {2\over3} \, \rho \, e )
{1\over33} \, \varepsilon - {1\over6} \, \Phi_{xx} +  {1\over2} \, \Phi_{ww}  +   {26\over33} \,\rho \, \lambda^2) 
+ \partial_z \Phi_{yz} = 0 \\ 
\partial_t j_z  + \partial_x \Phi_{zx}  + \partial_y \Phi_{yz}
+ \partial_z ( %%% {1\over3} \, \rho \, |{\bf u}|^2 - {1\over6} \, \Phi_{xx}  - {1\over2} \, \Phi_{ww} + {2\over3} \, \rho \, e   = 0 \\
{1\over33} \, \varepsilon - {1\over6} \, \Phi_{xx} -  {1\over2} \, \Phi_{ww}  +  {26\over33} \,\rho \, \lambda^2) = 0 \\ 
\partial_t \varepsilon  +  \partial_x \big( %%% 3 \, \Phi_{qx} + \rho \, u \, \lambda^2)
{11\over13} \, \Phi_{qx} +  {69\over13} \, \rho \, u \, \lambda^2 \big) + 
\partial_y \big(  {11\over13} \, \Phi_{qy} +  {69\over13} \, \rho \, v \, \lambda^2 \big)  + 
\partial_z \big(  {11\over13} \, \Phi_{qz} +  {69\over13} \, \rho \, w \, \lambda^2 \big)
%% +  \partial_y ( 3 \, \Phi_{qy} + \rho \, v \, \lambda^2)  +  \partial_z ( 3 \, \Phi_{qz} + \rho \, w \, \lambda^2)
= 0 \, . 
\end {array} \right. \monendstar 
The comparison  of the 3 momentum relations leads to the identities
\moneqstar   \left \{ \begin {array}{l}
{1\over33} \, \varepsilon + {1\over3} \, \Phi_{xx} +  {26\over33} \,\rho \, \lambda^2 =  \rho \, u^2 + p \\ 
{1\over33} \, \varepsilon - {1\over6} \, \Phi_{xx} +  {1\over2} \, \Phi_{ww}  +   {26\over33} \,\rho \, \lambda^2  =  \rho \, v^2 + p \\ 
{1\over33} \, \varepsilon - {1\over6} \, \Phi_{xx} -  {1\over2} \, \Phi_{ww}  +  {26\over33} \,\rho \, \lambda^2 =  \rho \, w^2 + p \\
\Phi_{xy} =  \rho \, u \, v \,,\,\, \Phi_{yz} =  \rho \, v \, w \,,\,\, \Phi_{zx} =  \rho \, u \, w  \, . 
\end {array} \right. \monendstar
By summation of the 3 first relations,  we observe that the quantity
\moneqstar
{1\over11} \, \varepsilon +  {26\over11} \,\rho \, \lambda^2 =  \rho \,  |{\bf u}|^2 + 3 \, p = 2 \, E + (  3 \, p - 2 \, \rho \, e )
\monendstar    
is conserved. Then we must have 
%%%  $ \, \rho \, |{\bf u}|^2   + 2  \, \rho \, e =  \rho \,  |{\bf u}|^2 + 3 \, p \, $ and we deduce that  
%
\moneqstar 
p = {2\over3} \, \rho \, e \,,\,\, \gamma \equiv {{c_p}\over{c_v}} = {5\over3} \, . 
\monendstar
We remark that for three space dimensions, this value of the ratio of specific heats $ \, \gamma \, $
corresponds to a monoatomic gas. The 3 degrees of freedom associated to translation
are the only taken into account for the constitution of the thermodynamical equilibrium.
Our result is analogous to the classical one for continuous Boltzmann equation.
To incorporate degrees of freedom associated to rotation and vibration, more complex models
have to be considered for this mesoscopic description,
intensively developed since the book of Hirschfelder, Curtiss and Bird \cite{HCF54}.
The Boltzmann approach with a finite set of velocities admits the same constraints than
the  continuous Boltzmann equation.

\smallskip \noindent  
The fifth  momentum $ \, \varepsilon \, $  of the D3Q33 lattice scheme is related to the total physical energy~$ \, E \, $ through the relation  
\moneqstar
\varepsilon = 22 \, E - 26 \, \lambda^2 \, \rho =  \rho \,  (11 \, |{\bf u}|^2  + 22 \, e  - 26 \, \lambda^2 )  \, . 
\monendstar
The 3 first equations for momentum take now the following form
\moneqstar   \left \{ \begin {array}{l}
{1\over3} \, \rho \, |{\bf u}|^2 + {1\over3} \, \Phi_{xx} + {2\over3} \, \rho \, e =  \rho \, u^2 + p \\ 
{1\over3} \, \rho \, |{\bf u}|^2 - {1\over6} \, \Phi_{xx}  + {1\over2} \, \Phi_{ww}+ {2\over3} \, \rho \, e  =  \rho \, v^2 + p \\ 
{1\over3} \, \rho \, |{\bf u}|^2 - {1\over6} \, \Phi_{xx}  - {1\over2} \, \Phi_{ww}+ {2\over3} \, \rho \, e  =  \rho \, w^2 + p 
\end {array} \right. \monendstar
and we deduce
$ \,\, \Phi_{xx}  =   \rho \, ( 2 \, u^2 - v^2 - w^2) \, $ and $ \, \,  \Phi_{ww}  =   \rho \, ( v^2 - w^2)  $. 
For the energy equation, we compare the equation
\moneqstar 
\partial_t \varepsilon  +  \partial_x \Big( %%% 3 \, \Phi_{qx} + \rho \, u \, \lambda^2)
{11\over13} \, \Phi_{qx} +  {69\over13} \, \rho \, u \, \lambda^2 \Big) + 
\partial_y \Big(  {11\over13} \, \Phi_{qy} +  {69\over13} \, \rho \, v \, \lambda^2 \Big)  + 
\partial_z \Big(  {11\over13} \, \Phi_{qz} +  {69\over13} \, \rho \, w \, \lambda^2 \Big) = 0
\monendstar
and 22 times the energy conservation minus $ \, 26  \, \lambda^2 \, $ times
the mass conservation, {\it id est}
\moneqstar 
\partial_t \varepsilon  + 
\partial_x ( \varepsilon  \, u + 22 \, p \, u )  + \partial_y ( \varepsilon  \, v +  22 \, p \, v )  + \partial_z ( \varepsilon  \, w +  22 \, p \, w ) = 0 \, .
\monendstar
Then we have the  relation
\moneqstar 
 {11\over13} \, \Phi_{qx} +  {69\over13} \, \rho \, u \, \lambda^2  =  \varepsilon  \, u + 22 \, p \, u
= \big( 11 \, \rho \, |{\bf u}|^2  + 22 \, \rho \, e  - 26 \, \lambda^2 \, \rho \big) \, u + 22 \, p \, u  
\monendstar
and analogous relations for $ \,  \Phi_{qy} \, $ and $ \,  \Phi_{qz} $.
The equilibrium values of the moments of the first family are finally given by 
\moneqstar   \left \{ \begin {array}{l}
\Phi_{xx}  =   \rho \, ( 2 \, u^2 - v^2 - w^2) \,, \,\,  \Phi_{ww}  =   \rho \, ( v^2 - w^2)  \\
\Phi_{xy}  =   \rho \, u \, v \,, \,\,  \Phi_{yz}  =   \rho \, v \,  w \,, \,\,  \Phi_{zx}  =   \rho \, w \, u  \\ 
\Phi_{qx}  = \rho \, u \, \xi_q \,,\,\, \Phi_{qy}  = \rho \, v \, \xi_q \,,\,\, \Phi_{qz}  = \rho \, v \, \xi_q \,,\,\, 
\xi_q =  13 \, | {\bf u} |^2  + {130\over3} \, e - 37 \, \lambda^2  \, . 
\end {array} \right. \monendstar

%%%%%%%%%%%%%%%%%%%%%%%%%%%%%%%%%%%%%%%%%%%%%%%%%%%%%%%%%%%%%%%%%%%%%%%%%%%%%%%%%%%%%%%%%%%%%%   d3q33 thermique ordre 2
\smallskip \monitem 
The second order accuracy is studied by comparing the second order term $\, \Gamma_2 = B \, \Sigma \, \Psi_1 \, $
and  the second order terms of the thermal Navier Stokes equations (\ref{flux-visqueus-ns-thermique-3d}), {\it id est} 
\moneqstar %%  \label{flux-visqueus-ns-thermique-3d}
\Phi_{\rm NS}^{3D} \! = \!\! \left( \!\! \begin {array}{l}
  \partial_j \sigma_{xj} \equiv \partial_x \big( 2 \, \mu \,   \partial_x u \!+\! (\zeta-{3\over3}\, \mu) \,  {\rm div}\,{\bf u} \big) 
\!+\! \partial_y (\mu ( \partial_x v  \!+\!   \partial_y u ) ) \!+\! \partial_z (\mu ( \partial_x w  \!+\!   \partial_z u ) )  \\
\partial_j \sigma_{yj} \equiv \partial_x (\mu ( \partial_x v  \!+\!   \partial_y u ) ) 
\!+\! \partial_y  \big( 2 \, \mu \,   \partial_y v \!+\! (\zeta-{3\over3}\, \mu) \,  {\rm div}\,{\bf u} \big) 
 \!+\! \partial_z (\mu ( \partial_y w  \!+\!   \partial_z v ) )  \\
 \partial_j \sigma_{zj} \equiv \partial_x (\mu ( \partial_x w  \!+\!   \partial_z u ) )
 \!+\! \partial_y (\mu ( \partial_y w  \!+\!   \partial_z v ) )  
\!+\! \partial_z \big( 2 \, \mu \,   \partial_z w \!+\! (\zeta-{3\over3}\, \mu) \,  {\rm div}\,{\bf u} \big)  \\
%% 26 \, [
\partial_j ( u_i \, \sigma_{ij} ) + {{\gamma}\over{Pr}} \, \big( \partial_x ( \mu \, \partial_x e  ) + \partial_y ( \mu \,  \partial_y e  )
+ \partial_z ( \mu \,  \partial_z e  )\big) %% ]
\end{array} \!\! \right)  \monendstar 
with $ \, {\rm div}\,{\bf u} \equiv \partial_x u +  \partial_y v +  \partial_z w $.
We have to solve a system of $ \, 180 = 4 \times 3 \times 5 \times 3 \, $ equations:
4 equations contain second order terms, each equation contains 3 space partial derivatives $\, \partial_x $ $\, \partial_y \, $
and $ \, \partial_z $, we have 5 nonconservative variables $ \, \rho  $,  $ \, u $, $ \, v $, $ \, w \, $ and $ \, e \, $ that define the viscous fluxes
and each of these fields occurs with  3 space partial derivatives.
The unknowns concern the 5  derivatives relative the nonconserved variables $ \, \rho  $,  $ \, u $, $ \, v $, $ \, w \, $ and $ \, e \, $
of the 16   moments 
$ \,  x_{yz} $, $ \, y_{zx}  $, $ \, z_{xy}  $, $ \,xyz  $, $ \, r_x   $, $ \,r_y    $, $ \, r_z  $,
$ \,   t_x  $, $ \, t_y   $, $ \,  t_z $, $ \, xx_e  $, $ \,  ww_e   $, $ \,   xy_e   $, $  yz_e  $, $\,  zx_e  \, $ and $ \, h \, $ 
of the  second family, then a total of 80 unknowns.
Due to the degeneracy of this set of equations, it is possible to find a family of solutions for the equilibria: 
\moneqstar   \left \{ \begin {array}{l}
\Phi_{x \,\, yz}  = \rho \,  u \, (v^2 - w^2) \,,\,\, \Phi_{y\,\, zx}  = \rho \,  v \, (w^2 - u^2) \,,\,\,
\Phi_{z\,\, xy}  = \rho \,  w \, (u^2 - v^2)  \\
\Phi_{xyz}  = \rho \,  u \, v \, w \\   
\Phi_{rx} + {38\over13} \,  {1\over{\lambda^2}} \,   \Phi_{tx}  =
\rho \, u  \, \lambda^2 \, \big( -  {161\over39} \, u^2 -  {345\over13} \, (v^2+w^2) -   {1702\over39} \, e  +  {1265\over39} \, \lambda^2 \big) \\
\Phi_{ry} +  {38\over13} \,  {1\over{\lambda^2}} \,   \Phi_{ty}  =
\rho \, v  \, \lambda^2 \, \big( -  {161\over39} \, v^2 -  {345\over13} \, (w^2+u^2) -   {1702\over39} \, e  +  {1265\over39} \, \lambda^2 \big) \\
\Phi_{rz} +  {38\over13} \,  {1\over{\lambda^2}} \,   \Phi_{tz} =
\rho \, w  \, \lambda^2 \, \big( -  {161\over39} \, w^2 -  {345\over13} \, (u^2+v^2) -   {1702\over39} \, e  +  {1265\over39} \, \lambda^2 \big) \\  
\Phi_{xxe}  = \rho \,  (2 \, u^2 - v^2 - w^2 ) \, \big( 38 \, | {\bf u} |^2  + {266\over3} \, e - 38 \, \lambda^2 \big)  \\
\Phi_{wwe}  =  \rho \,  (v^2 - w^2 ) \, \big( 38 \, | {\bf u} |^2  + {266\over3} \, e - 38 \, \lambda^2 \big) \\ 
\Phi_{xye} = \rho \,  u \,v \, \xi_e \,,\,\, \Phi_{yze} = \rho \,  v \,w \, \xi_e \,,\,\,\Phi_{zxe} = \rho \,  w \, u \, \xi_e \,,\,\,
\xi_e = 3  \, |{\bf u}|^2 + 14 \, e - 8 \, \lambda^2 \\   
\Phi_{h} = \rho \, \big( {69\over2} \, | {\bf u} |^4 + 230 \, ( | {\bf u} |^2 + e ) \, e
- 325 \, \lambda^2 \, ( {1\over2} \,  | {\bf u} |^2 +  e)  + 152 \, \lambda^4  \big) \, . 
\end {array} \right. \monendstar
Moreover, under the condition $ \, \sigma_x  =  \sigma_q  $, it is possible to match the viscous fluxes
(\ref{flux-visqueus-ns-thermique-3d}) of the thermal Navier Stokes equations, with the following viscous parameters
\moneq \label{viscosite-d3q33-thermique} 
%%%%%%%%%%  eqq01f =  zeta
%%%%%%%%%%  eqq04f =  (-rho)*ee*sig_x*Delta_t + mu
\mu = {2\over3} \, \rho \, e \, \sigma_x \, \Delta t \,,\,\,    \zeta = 0  \,,\,\,   Pr = 1 \, .
\monend
The D3Q33 lattice Boltzmann scheme is compatible with three-dimensional thermal Navier Stokes equations.
Observe that the number of velocities remains moderate compared to the number of velocities proposed {\it e.g.}
in \cite{SPP09} and used in \cite{WKBRF21}.

%%%%%%%%%%%%%%%%%%%%%%%%%%%%%%%%%%%%%%%%%%%%%%%%%%%%%%%%%%%%%%%%%%%%%%%%%%%%%%%%%%%%
\smallskip \monitem  D3Q27-2 %%% Surprising
%%%%%%%%%%%%%%%%%%%%%%%%%%%%%%%%%%%%%%%%%%%%%%%%%%%%%%%%%%%%%%%%%%%%%%%%%%%%%%%%%%%%

\noindent
The set of velocities of the D3Q27-2 lattice Boltzmann scheme is  described in the Figure~\ref{fig-d3q27-2}. 
The construction of  moments from polynomials is detailed in the Table \ref{d3q27-2-ns-polynomes} of the Appendix.
The first 5 moments  are conserved,
8 can be tuned in order to fit the first order Euler equations,  13 can be adjusted for the viscous terms
and  1  has no  influence on the Navier Stokes equations,
as summarized in Table \ref{d3q27-2-moments-thermique}. 
%

%%%%%%%%%%%%%%%%%%%%%%%%%%%%%%%%%%%%%%%%%%%%%%%%%%%%%%%%%%%%%%%%%%%%%%%%%%%%%%%%%%% moments d3q27-2   
\begin{table}    [H]  \centering
\begin{tabular}{|c|l|c|c|} \hline
 & conserved &  $\rho \,,\,\, j_x \,,\,\, j_y \,,\,\, j_z  \,,\,\, \varepsilon $ & 5  \\  \hline
 1 & fit the Euler equations   &  $ xx  \,,\,\,  ww \,,\,\, xy   \,,\,\,  yz \,,\,\,  zx \,,\,\,  q_x  \,,\,\, q_y  \,,\,\, q_z$ & 8 \\ \hline 
2 & fit the viscous terms & $ x_{yz}  \,,\,\,  y_{zx}   \,,\,\,   z_{xy} \,,\,\, xyz  \,,\,\, r_x  \,,\,\, r_y  \,,\,\, r_z   $  &  \\
&  & $ h  \,,\,\,  xx_e  \,,\,\,   ww_e  \,,\,\,  xy_e  \,,\,\,  yz_e  \,,\,\,  zx_e $  & 13 \\ \hline 
3 & without influence & $    h_3 $  & 1  \\ \hline 
\end{tabular} 
\caption{D3Q27-2 moments for thermal  Navier Stokes equations}
\label{d3q27-2-moments-thermique} \end{table}
\smallskip \monitem
The equivalent first order partial differential equations 
\moneqstar   \left \{ \begin {array}{l}
 \partial_t \rho + \partial_x j_x + \partial_y j_y + \partial_z j_z = 0 \\ 
\partial_t j_x + \partial_x ({1\over9} \, \varepsilon + {1\over3} \, \Phi_{xx} +  {8\over9} \, \rho \, \lambda^2 )
+ \partial_y  \Phi_{xy} + \partial_z  \Phi_{zx} = 0 \\
\partial_t j_y  + \partial_x \Phi_{xy}   + \partial_y ({1\over9} \, \varepsilon  - {1\over6} \, \Phi_{xx}  + {1\over2} \, \Phi_{ww}
+  {8\over9} \, \rho \, \lambda^2 )  + \partial_z \Phi_{yz} = 0 \\ 
\partial_t j_z  + \partial_x \Phi_{zx}  + \partial_y \Phi_{yz}   + \partial_z ( {1\over9} \, \varepsilon  - {1\over6} \, \Phi_{xx}
- {1\over2} \, \Phi_{ww}  + {8\over9} \, \rho \, \lambda^2  ) = 0 \\
\partial_t \varepsilon  +  \partial_x ( 3 \, \Phi_{qx} + \rho \, u \, \lambda^2)
 +  \partial_y ( 3 \, \Phi_{qy} + \rho \, v \, \lambda^2)  +  \partial_z ( 3 \, \Phi_{qz} + \rho \, w \, \lambda^2)  = 0
\end {array} \right. \monendstar
are compared to the Euler equations of gas dynamics (\ref{euler-3d-avec-energie}). 
The equality of the 3 equations for momentum implies   
\moneqstar   \left \{ \begin {array}{l}
{1\over9} \, \varepsilon + {1\over3} \, \Phi_{xx} +  {8\over9} \, \rho \, \lambda^2  =   \rho \, u^2 + p \\ \vspace{-.45 cm} \\
 {1\over9} \, \varepsilon   - {1\over6} \, \Phi_{xx}  + {1\over2} \, \Phi_{ww}  + {8\over9} \, \rho \, \lambda^2 =   \rho \, v^2 + p \\ \vspace{-.45 cm} \\
 {1\over9} \, \varepsilon   - {1\over6} \, \Phi_{xx}  - {1\over2} \, \Phi_{ww}  + {8\over9} \, \rho \, \lambda^2  = \rho \, w^2 + p \, . 
\end {array} \right. \monendstar
We deduce $ \,  {1\over3} \, \varepsilon + {8\over3} \, \rho \, \lambda^2   =  \rho \, |{\bf u}|^2 + 3 \, p = 2 \, E   - 2 \, \rho \, e  + 3 \, p $.
Analogously to the D3Q33 scheme, the left hand side is a conserved quantity and it is also the case for the physical total energy~$ \, E $.
Then we have the constraint  $ \, p = {2\over3} \, \rho \, e \, $ and $ \,  \gamma \equiv {{c_p}\over{c_v}} = {5\over3} $.
We have also the relation 
\moneqstar
\varepsilon = 6 \, E - 8 \, \rho \, \lambda^2 = \rho \, ( 3 \,  |{\bf u}|^2 + 6 \, e -  8  \, \lambda^2 ) \, . 
\monendstar
It is then elementary to fit the 8 moments of the first family:
\moneqstar   \left \{ \begin {array}{l}
\Phi_{xx}  =   \rho \, ( 2 \, u^2 - v^2 - w^2) \,, \,\,  \Phi_{ww}  =   \rho \, ( v^2 - w^2)  \\
\Phi_{xy}  =   \rho \, u \, v \,, \,\,  \Phi_{yz}  =   \rho \, v \,  w \,, \,\,  \Phi_{zx}  =   \rho \, w \, u   \\ 
\Phi_{qx}  =  \rho \, u \, \xi_q \,,\,\, \Phi_{qy}  =  \rho \, v \, \xi_q \,,\,\, \Phi_{qz}  =  \rho \, w \, \xi_q \,,\,\,
\xi_q = ( | {\bf u} |^2  + {10\over3} \, e - 3 \, \lambda^2 ) \, . 
\end {array} \right. \monendstar

%%%%%%%%%%%%%%%%%%%%%%%%%%%%%%%%%%%%%%%%%%%%%%%%%%%%%%%%%%%%%%%%%%%%%%%%%%%%%%%%%%%%%%%%%%%%%%   d3q27-2 thermique ordre 2
\smallskip \monitem
The calibration  of viscous fluxes demands the resolution of 180 linear equations as observed for the D3Q33 scheme.
We have 13 moments in the second family of nonconserved moments, then $ \, 5 \times 13 = 65 \, $ unknowns.
If we suppose $ \, \sigma_x = \sigma_q $, we observe that
the equilibrium function is determined for all  microscopic moments of the second family and 
the D3Q27-2  scheme is available for thermal Navier Stokes: 
\moneqstar   \left \{ \begin {array}{l} 
\Phi_{x \,\, yz}  \, = \rho \,  u \, (v^2 - w^2) \,,\,\,  \Phi_{y\,\, zx}  \,  = \rho \,  v \, (w^2 - u^2)  \,,\,\, 
\Phi_{z\,\, xy}  \,  = \rho \,  w \, (u^2 - v^2)  \\
\Phi_{xyz}  = \rho \,  u \, v \, w  \\
\Phi_{rx}   =   \rho \, u \, \lambda^2 \, \big( - (  u^2 +  3 \, v^2 + 3 \, w^2 ) -   6  \, e +  5\, \lambda^2 \big)  \\ 
\Phi_{ry}   =   \rho \, v \, \lambda^2 \, \big( - (  3 \, u^2 +  v^2 + 3 \, w^2 ) -   6  \,  e +  5\, \lambda^2 \big) \\ 
\Phi_{rz}   =   \rho \, w \, \lambda^2 \, \big( - (  3 \, u^2 +  3 \, v^2 + w^2 ) -   6  \, e +  5\, \lambda^2 \big) \\ 
 \Phi_{h}   = \rho \, \big( {3\over2}  \, | {\bf u} |^4 + 10 \, ( | {\bf u} |^2 + e ) \, e
- 15 \, \lambda^2 \, ( {1\over2}  \, | {\bf u} |^2 + e)  + 8 \, \lambda^4  \big) \\ 
\Phi_{xxe}   = \rho \,  (2 \, u^2 - v^2 - w^2 ) \, \big( {9\over8}  \, | {\bf u} |^2  + {21\over4}  \, e - {17\over4} \, \lambda^2 \big) \\  
\Phi_{wwe}  \,\, =  \rho \,  (v^2 - w^2 ) \, \big(  {9\over8}  \, | {\bf u} |^2  + {21\over4}  \, e - {17\over4} \, \lambda^2 \big) \\  
\Phi_{xye}   = \rho \,  u \,v \, \xi_e \,,\,\, \Phi_{yze}   = \rho \,  v \,w \, \xi_e \,,\,\, \Phi_{zxe}   = \rho \,  w \, u \, \xi_e 
\,,\,\, \xi_e = (3  \, |{\bf u}|^2 + 14 \, e - 8 \, \lambda^2 )   \, .  
\end {array} \right. \monendstar
The viscosities and the Prandtl number are given by the relations (\ref{viscosite-d3q33-thermique}),
as for the D3Q33 scheme.
%%%% {\color{red} 
We observe that the constraint for the Prandtl number results from the resolution of the system
of 180 equations to enforce the thermal Navier Stokes equations at second order accuracy.
If we tune the relaxation coefficients of the
multi-resolution-times lattice Boltzmann scheme, we can change the Prandtl number
but in that case, at least one of the previous equations is not satisfied 
and we do not have a formal approximation of the Navier Stokes
equations at second order accuracy. 
The theoretical ability of the D3Q27-2 scheme to approximate formally the thermal Navier Stokes equations at second order 
accuracy is a good surprise of this  work.

%%% \newpage 
%%%%%%%%%%%%%%%%%%%%%%%%%%%%%%%%%%%%%%%%%%%%%%%%%%%%%%%%%%%%%%%%%%%%%%%%%%%%%%%  section 8
\bigskip \bigskip    \noindent {\bf \large    8) \quad  Conclusion }
%%%%%%%%%%%%%%%%%%%%%%%%%%%%%%%%%%%%%%%%%%%%%%%%%%%%%%%%%%%%%%%%%%%%%%%%%%%%%%%%%%%%%%%%%%%%%%%%%%%%%%

\noindent
In this contribution, we have used the Taylor expansion method
for multiresolution times lattice Boltzmann schemes in the framework of a single
particle distribution for the approximation of compressible Navier Stokes equations.
Equilibrium functions have been explicited  in order to fit formally
the second order equivalent partial differential equations
with the Navier Stokes equations.

\smallskip \noindent
A summary of our general conclusions is proposed at Table \ref{synthese}.

\smallskip \noindent
In the isothermal case, we recover known results concerning 
the impossibility to recover exactly the  Navier Stokes equations
with the D2Q9 scheme. Moreover, the ability of the  D2Q13 scheme
to fit the  Navier Stokes equations is demonstrated and
appropriate equilibrium function are  explicited. 
At three space dimensions, the defects of the D3Q19 scheme were  known.
But we have been surprised by very  analogous defects exhibited for the
classical D3Q27 scheme.
A natural question for future studies is to understand if it is possible to find an other  family of
moments for the D3Q27 scheme that would allow the recover the entire set of  isothermal viscous equations.
With the D3Q33 scheme and the version ``D3Q27-2''
scheme with velocities of the type $ \, (1,\, 0,\, 0) \, $ replaced
by $  \, (2,\, 0,\, 0) $, isothermal  Navier Stokes equations can be formally captured
by the asymptotic expansion and the associated equilibria are presented.
% %%  \smallskip \noindent
Observe also that for each scheme, we have chosen and fixed a moment matrix.
%% A natural question for future studies is the variation of the present results
%% when this matrix is changed.
More generally  a natural question for future studies  is to understand how the choice of the moments
affect precisely the structure  of the momentum-velocity operator matrix. %%  $ \, \Lambda $. 

%%%%%%%%%%%%%%%%%%%%%%%%%%%%%%%%%%%%%%%%%%%%%%%%%%%%%%%%%%%%%%%%%%%%%%%%%%%%%%%%%%%  synthese graphique des resultats 
\begin{table}    [H]  \centering 
\centerline{  \begin{tabular}{|c|c|c|} \hline
$\!\!$energy$\!\!$&  not appropriate &  possible \\  \hline
&  {\includegraphics[width=.18\textwidth]{d2q9-18mars2021-fig.pdf}} &  {\includegraphics[width=.18\textwidth]{d2q13-mars2021-fig.pdf}}  \\ 
& {\footnotesize  D2Q9} & {\footnotesize  D2Q13} \\    
$\!\!$without$\!\!$& & \\
& {\includegraphics[width=.15\textwidth]{lbm-graphes-d3q19-01mars2022.pdf}} {\includegraphics[width=.15\textwidth]{lbm-graphes-d3q27-30mai2020.pdf}} &
{\includegraphics[width=.18\textwidth]{lbm-graphes-d3q33-01juin2020.pdf}} {\includegraphics[width=.18\textwidth]{lbm-graphes-d3q272-15mars2021.pdf}} \\
&  {\footnotesize  D3Q19} \qquad \qquad \quad  {\footnotesize D3Q27} & {\footnotesize D3Q33}  \qquad \qquad \quad  {\footnotesize D3Q27-2}  \\     \hline
& {\includegraphics[width=.18\textwidth]{d2q13-mars2021-fig.pdf}} & 
$\!\!\!\!${\includegraphics[width=.18\textwidth]{d2q17-18mars2021-fig.pdf}}$\!\!\!\!\!\!\!$ {\includegraphics[width=.18\textwidth]{d2v17-18mars2021-fig.pdf}}$\!\!\!\!\!\!\!$  
{\includegraphics[width=.18\textwidth]{d2w17-23mars2021-fig.pdf}}$\!\!\!\!$\\ 
&  {\footnotesize D2Q13} &  {\footnotesize D2Q17} \qquad  \quad {\footnotesize D2V17} \quad \qquad {\footnotesize D2W17} \\   
$\!\!$with$\!\!$& & \\
&$\!\!${\includegraphics[width=.15\textwidth]{lbm-graphes-d3q19-01mars2022.pdf}} {\includegraphics[width=.15\textwidth]{lbm-graphes-d3q27-30mai2020.pdf}} &
{\includegraphics[width=.18\textwidth]{lbm-graphes-d3q33-01juin2020.pdf}} {\includegraphics[width=.18\textwidth]{lbm-graphes-d3q272-15mars2021.pdf}} \\
&  {\footnotesize D3Q19} \qquad \qquad \quad  {\footnotesize D3Q27} & {\footnotesize D3Q33}  \qquad \qquad \quad  {\footnotesize D3Q27-2}  \\ \hline
\end{tabular}  }
\caption{Mutiresolution times lattice Boltzmann schemes with a single particle distribution: summary of our results for compressible Navier Stokes equations}
\label{synthese} \end{table}
%%%%%%%%%%%%%%%%%%%%%%%%%%%%%%%%%%%%%%%%%%%%%%%%%%%%%%%%%%%%%%%%%%%%%%%%%%%%%%%%%%%

\smallskip \noindent
In the thermal case, the D2Q13 scheme has signifiant   defects and we recommend
to test one of the D2Q17, D2V17 or D2W17 schemes which have no formal incoherence 
for the approximation of Navier Stokes equations including the conservation of energy.
At three space dimensions, the previously studied  D3Q33 and D3Q27-2
schemes have the ability to fit exactly the thermal viscous fluid equations.
We remark  also that the relations we have obtained for the equilibrium function
are pure discrete relations and have no direct relation with the physical
gaussian equilibrium. 
Observe finally that some restrictions are imposed to the physical model.
In particular, the undelying gas is monoatomic ($ \gamma = {5\over3} $)
and the Prandtl number is restricted to~1.
In consequence,  it is probably necessary to revise the present paradigm
to overcome the observed physical limitations.
%% of $ \, \gamma = 2 \, $ for two space dimensions, 
%% $ \, \gamma = {5\over3} \, $ for three  space dimensions and Prandtl number equal to unity ? 
%
%%%% {\color{red}
A possibility is to adapt  the present paradigm to take into account 
the gas kinetic Shakhov model or the ellipsoidal statistical model (see {\it e.g.} \cite{BY23,CXC15,WWHLLZ20}).
%%% }

\smallskip \noindent
In our analysis, we have  introduced three  families of nonconserved moments.
The first family is devoted to the adjustement of first order Euler inviscid  fluxes.
%%% {\color{red}
These moments are homogeneous polynomials of degree 2 relative to the velocities.
If the chosen moments do not satisfy this property, we observe that all the structure we have described in
the Tables~\ref{d2q9-moments-isotherme} to~\ref{d3q27-2-moments-thermique} of our contribution is no more valid.
%%%  }
%
The second family allows the tuning of second order viscous fluxes.

\smallskip \noindent
Finally, the third family has no impact on the equivalent partial differential
equations.
%
%%% {\color{red}
It is a naturel idea to enforce the nullity at equilibrium for the third family of moments.
But a lattice Boltzmann scheme does not reduce to its equivalent partial differential equation
at second order accuracy! To illustrate this point, we can refer to the D2Q9 scheme studied
in detail in the references \cite{De14,LL00}. The last moment $ \, h \, $ in our present nomenclature (see Table~\ref{polynomes-d2q9})
has no influence on the partial differential equations at second order.
Nevertheless, to enforce Galilean invariance and strengthen the stability of the scheme, 
its equilibrium value is taken classically to
$ \, h^{\rm eq} = -2 \, \lambda^2 \, \rho + 3 \, |u|^2 \, \rho \, $
as proposed in  \cite{LL00}. The question of an optimal equilibrium values for microscopic moments of the third
family is one of the questions to be discussed in future works.  
%%%  }

\smallskip  \noindent  %%%   {\color{red}
The scheme has also to be coupled with an appropriate set of boundary conditions.
In particular, we have to develop an extension of the bounce-back algorithm
to take into account null velocity and given temperature on the boundary or
null velocity and given flux on the boundary.
A study of the type has been presented for the D2Q13 scheme in  \cite{LD15}.
%%%  }

\smallskip \noindent
Last but not least, a fundamental question concerns  stability. Our study is purely algebraical and
both linear and nonlinear stabilities have to be experimented
for specific schemes and test cases.
%
%%%  {\color{red}
An other approach recently into development~\cite{DST22} is purely experimental. It consists in a global optimization
of the parameters of a muti-resolution-times scheme 
in order to maintain experimentally stability in a given interval of time. 
%%%  }

\newpage 
%%%%%%%%%%%%%%%%%%%%%%%%%%%%%%%%%%%%%%%%%%%%%%%%%%%%%%%%%%%%%%%%%%%%%%%%%%%%%%%  09 novembre 2021
\bigskip \bigskip    \noindent {\bf \large     Appendix. \quad Velocity polynomials for the d'Humi\`eres matrix}
%%%%%%%%%%%%%%%%%%%%%%%%%%%%%%%%%%%%%%%%%%%%%%%%%%%%%%%%%%%%%%%%%%%%%%%%%%%%%%%  09 novembre 2021

The matrix of moments $ \,  M_{k j} \, $ is defined thanks to polynomials $ \, p_k \,$
with the relation
\moneqstar 
M_{k j} = p_k(v_{j,\,x} ,\, v_{j,\,y} ,\, v_{j,\,z}) \, .
\monendstar
Moreover, the moments have been chosen orthogonal in the following sense:
\moneqstar 
\sum_j \, M_{k j} \, M_{\ell j} = 0 \,\,\, {\rm  when}  \,\,\,\,  k \not = \ell 
\monendstar
as suggested in the classic article \cite{LL00}. 

%%% \newpage 
%%%%%%%%%%%%%%%%%%%%%%%%%%%%%%%%%%%%%%%%%%%%%%%%%%%%%%%%%%%%%%%%%%%%%%%%%%%%%%% 
\smallskip \monitem D2Q9 for isothermal Navier Stokes 

\smallskip

%%%%%%%%%%%%%%%%%%%%%%%%%%%%%%%%%%%%%%%%%%%%%%%%%%%%%%%%%%%%%%%%%%%%%%%%%%%%%%%%%%% moments d2q9   
\begin{table}    [H]  \centering
\centerline{
\begin{tabular}{|l|c|c|} \hline
$ k $ &  $m_k$ &  $p_k (v)$ \\  \hline
$ 0 $ &  $\rho$ &  $ 1 $ \\  
$ 1 $ &  $j_x $ &  $ v_x $ \\  
$ 2 $ &  $j_y $ &  $ v_y $ \\   \hline
$ 3 $ &  $ \varepsilon $ &  $  3 \, (v_x^2 + v_y^2) - 4 \, \lambda^2 $ \\  
$ 4 $ &  $ xx $ &  $  v_x^2 - v_y^2 $ \\  
$ 5 $ &  $ xy $ &  $  v_x \, v_y $ \\   \hline
$ 6 $ &  $ q_x $ &  $  [3 \, (v_x^2 + v_y^2) - 5  \, \lambda^2 ] \, v_x  $ \\  
$ 7 $ &  $ q_y $ &  $  [3 \, (v_x^2 + v_y^2) - 5  \, \lambda^2 ]  v_y  $ \\   \hline
$ 8 $ &  $ h $ &  $  [ {9\over2} \, (v_x^2 + v_y^2)^2 - {21\over2} \,\lambda^2\,(v_x^2 + v_y^2) + 4 \,\lambda^4 ] $ \\ 
%%%  \vskip -.5cm  &  &  \\ 
\hline \end{tabular}}
  
%%%   nouvelle edition des polynomes
%%%   0   1
%%%   1   vx
%%%   2   vy
%%%   3   3*vx^2 + 3*vy^2 - 4*ll^2
%%%   4   vx^2 - vy^2
%%%   5   vx*vy
%%%   6   3*vx^3 + 3*vx*vy^2 + (-5*ll^2)*vx
%%%   7   3*vx^2*vy + 3*vy^3 + (-5*ll^2)*vy
%%%   8   9/2*vx^4 + 9*vx^2*vy^2 + 9/2*vy^4 + (21*ll^2/(-2))*vx^2 + (21*ll^2/(-2))*vy^2 + 4*ll^4

\caption{Moment polynomials for  D2Q9 isothermal Navier Stokes}
\label{polynomes-d2q9} \end{table}
%%%%%%%%%%%%%%%%%%%%%%%%%%%%%%%%%%%%%%%%%%%%%%%%%%%%%%%%%%%%%%%%%%%%%%%%%%%%%%%%%%%

\bigskip %%% \newpage 
%%%%%%%%%%%%%%%%%%%%%%%%%%%%%%%%%%%%%%%%%%%%%%%%%%%%%%%%%%%%%%%%%%%%%%%%%%%%%%% 
\monitem D2Q13 for isothermal Navier Stokes 

\smallskip

%%%%%%%%%%%%%%%%%%%%%%%%%%%%%%%%%%%%%%%%%%%%%%%%%%%%%%%%%%%%%%%%%%%%%%%%%%%%%%%%%%% moments d2q13  
\begin{table}    [H]  \centering
\centerline{
\begin{tabular}{|l|c|c|} \hline
$ k $ &  $m_k$ &  $p_k (v)$ \\  \hline
$ 0 $ &  $\rho$ &  $ 1 $ \\  
$ 1 $ &  $j_x $ &  $ v_x $ \\  
$ 2 $ &  $j_y $ &  $ v_y $ \\   \hline
$ 3 $ &  $ \varepsilon $ &  $  13 \, (v_x^2 + v_y^2) - 28 \, \lambda^2 $ \\  
$ 4 $ &  $ xx $ &  $  v_x^2 - v_y^2 $ \\  
$ 5 $ &  $ xy $ &  $  v_x \, v_y $ \\   \hline 
$ 6 $ &  $ q_x $ &  $  (v_x^2 + v_y^2 - 3  \, \lambda^2 )\, v_x  $ \\  
$ 7 $ &  $ q_y $ &  $  (v_x^2 + v_y^2 - 3  \, \lambda^2 )\, v_y  $ \\  
$ 8 $ &  $ r_x $ &  $  [ {35\over12} \, (v_x^2 + v_y^2)^2 - {63\over4} \,\lambda^2\,(v_x^2 + v_y^2) + {101\over6}\,\lambda^4 ] \, v_x  $ \\  
$ 9 $ &  $ r_y $ &  $  [ {35\over12} \, (v_x^2 + v_y^2)^2 - {63\over4} \,\lambda^2\,(v_x^2 + v_y^2) + {101\over6}\,\lambda^4 ] \, v_y  $ \\   \hline 
$ 10 $ &  $ h $ &  $   {77\over2} \, (v_x^2 + v_y^2)^2 -  {361\over2} \,\lambda^2\,(v_x^2 + v_y^2) + 140 \,\lambda^4   $ \\  
$ 11 $ &  $ xx_e  $ &  $   [ {17\over12} \, (v_x^2 + v_y^2) -  {65\over12} \,\lambda^2 ] \, (v_x^2 - v_y^2)  $ \\  
$ 12 $ &  $ h_3 $ &  $   {137\over24} \, (v_x^2 + v_y^2)^3 -  {273\over8} \,\lambda^2 \,(v_x^2 + v_y^2)^2   +  {581\over12} \,\lambda^4 \,(v_x^2 + v_y^2)
                         -12 \,\lambda^6      $ \\  
                         \hline \end{tabular}}

\caption{Moment polynomials for  D2Q13 isothermal Navier Stokes}
\label{polynomes-d2q13} \end{table}

%%%   polynomes definitifs qui prennent en compte l orthogonalisation
%%%   0   1
%%%   1   vx
%%%   2   vy
%%%   3   13*vx^2 + 13*vy^2 - 28*ll^2
%%%   4   vx^2 - vy^2
%%%   5   vx*vy
%%%   6   vx^3 + vx*vy^2 + (-3*ll^2)*vx
%%%   7   vx^2*vy + vy^3 + (-3*ll^2)*vy
%%%   8   35/12*vx^5 + 35/6*vx^3*vy^2 + 35/12*vx*vy^4 + (63*ll^2/(-4))*vx^3 + (63*ll^2/(-4))*vx*vy^2 + 101*ll^4/6*vx
%%%   9   35/12*vx^4*vy + 35/6*vx^2*vy^3 + 35/12*vy^5 + (63*ll^2/(-4))*vx^2*vy + (63*ll^2/(-4))*vy^3 + 101*ll^4/6*vy
%%%   10  77/2*vx^4 + 77*vx^2*vy^2 + 77/2*vy^4 + (361*ll^2/(-2))*vx^2 + (361*ll^2/(-2))*vy^2 + 140*ll^4
%%%   11  17/12*vx^4 + ((-17)/12)*vy^4 + (65*ll^2/(-12))*vx^2 + ((-65*ll^2)/(-12))*vy^2
%%%   12  137/24*vx^6 + 137/8*vx^4*vy^2 + 137/8*vx^2*vy^4 + 137/24*vy^6 + (273*ll^2/(-8))*vx^4 + (273*ll^2/(-4))*vx^2*vy^2 + (273*ll^2/(-8))*vy^4 + ((-581*ll^4)/(-12))*vx^2 + ((-581*ll^4)/(-12))*vy^2 - 12*ll^6

%%%%%%%%%%%%%%%%%%%%%%%%%%%%%%%%%%%%%%%%%%%%%%%%%%%%%%%%%%%%%%%%%%%%%%%%%%%%%%%%%%%

\newpage 
%%%%%%%%%%%%%%%%%%%%%%%%%%%%%%%%%%%%%%%%%%%%%%%%%%%%%%%%%%%%%%%%%%%%%%%%%%%%%%% 
\monitem D3Q19 for isothermal Navier Stokes 

\smallskip

%%%%%%%%%%%%%%%%%%%%%%%%%%%%%%%%%%%%%%%%%%%%%%%%%%%%%%%%%%%%%%%%%%%%%%%%%%%%%%%%%%% moments d3q19  
\begin{table}    [H]  \centering
\centerline{
\begin{tabular}{|l|c|c|} \hline
$ k $ &  $m_k$ &  $p_k (v)$ \\  \hline
$ 0 $ &  $\rho$ &  $ 1 $ \\  
$ 1 $ &  $j_x $ &  $ v_x $ \\  
$ 2 $ &  $j_y $ &  $ v_y $ \\  
$ 3 $ &  $j_z $ &  $ v_z $ \\   \hline
$ 4 $ &  $ \varepsilon $ &  $  19 \, (v_x^2 + v_y^2 + v_z^2) - 30 \, \lambda^2 $ \\  
$ 5 $ &  $ xx $ &  $  2 \, v_x^2 - v_y^2 - v_z^2  $ \\  
$ 6 $ &  $ ww $ &  $  v_y^2 - v_z^2 $ \\  
$ 7 $ &  $ xy $ &  $  v_x \, v_y $ \\   
$ 8 $ &  $ yz $ &  $  v_y \, v_z $ \\  
$ 9 $ &  $ zx $ &  $  v_z \, v_x $ \\   \hline 
$ 10 $ &  $ q_x $ &  $  [ 5 \, (v_x^2 + v_y^2 + v_z^2)- 9  \, \lambda^2 ]\, v_x  $ \\  
$ 11 $ &  $ q_y $ &  $   [ 5 \, (v_x^2 + v_y^2 + v_z^2)- 9  \, \lambda^2 ]\, v_y  $ \\  
$ 12 $ &  $ q_z $ &  $  [ 5 \, (v_x^2 + v_y^2 + v_z^2)- 9  \, \lambda^2 ]\, v_z  $ \\  
$ 13 $ &  $ x_{yz} $ &  $  v_x \, ( v_y^2 - v_z^2 )  $ \\  
$ 14 $ &  $ y_{zx} $ &  $  v_y \, ( v_z^2 - v_x^2 )  $ \\  
$ 15 $ &  $ z_{xy} $ &  $  v_z \, ( v_x^2 - v_y^2 )  $ \\    \hline 
$ 16 $ &  $ h $ &  $   {21\over2} \, (v_x^2 + v_y^2 + v_z^2)^2 -  {53\over2} \,\lambda^2 \,(v_x^2 + v_y^2 + v_z^2) + 12 \,\lambda^4   $ \\  
$ 17 $ &  $ xx_e  $ &  $   [ 3 \, (v_x^2 + v_y^2 + v_z^2) - 5 \,\lambda^2 ] \, ( 2 \, v_x^2 - v_y^2 - v_z^2 )  $ \\  
$ 18 $ &  $ ww_e  $ &  $  [ 3 \, (v_x^2 + v_y^2 + v_z^2) - 5 \,\lambda^2 ] \, (  v_y^2 - v_z^2  )  $ \\  
\hline \end{tabular}}

\caption{Moment polynomials for  D3Q19 isothermal Navier Stokes}
\label{d3q19-iso-polynomes} \end{table}

%%%    print ("d3q19 ordre des moments") 
%%%    print ("        conserves                   rho, JJx, JJy, JJz")
%%%    print ("        regles au premier ordre     ee, xx, ww, xy, yz, zx")
%%%    print ("        vecteurs de chaleur         qx, qy, qz")
%%%    print ("        ordre trois antisymetrique  x-yz, y-zx, z-xy")
%%%    print ("        energie d ordre superieur   hh")
%%%    print ("        tenseurs d ordre superieur  xxe, wwe") 

%%%    polynomes apres orthogonalisation 
%%%    0    1
%%%    1    vx
%%%    2    vy
%%%    3    vz
%%%    4    19*vx^2 + 19*vy^2 + 19*vz^2 - 30*ll^2
%%%    5    2*vx^2 - vy^2 - vz^2
%%%    6    vy^2 - vz^2
%%%    7    vx*vy
%%%    8    vy*vz
%%%    9    vx*vz
%%%    10   5*vx^3 + 5*vx*vy^2 + 5*vx*vz^2 + (-9*ll^2)*vx
%%%    11   5*vx^2*vy + 5*vy^3 + 5*vy*vz^2 + (-9*ll^2)*vy
%%%    12   5*vx^2*vz + 5*vy^2*vz + 5*vz^3 + (-9*ll^2)*vz
%%%    13   vx*vy^2 - vx*vz^2
%%%    14   -vx^2*vy + vy*vz^2
%%%    15   vx^2*vz - vy^2*vz
%%%    16   21/2*vx^4 + 21*vx^2*vy^2 + 21/2*vy^4 + 21*vx^2*vz^2 + 21*vy^2*vz^2 + 21/2*vz^4 + (53*ll^2/(-2))*vx^2 + (53*ll^2/(-2))*vy^2 + (53*ll^2/(-2))*vz^2 + 12*ll^4
%%%    17   6*vx^4 + 3*vx^2*vy^2 + (-3)*vy^4 + 3*vx^2*vz^2 + (-6)*vy^2*vz^2 + (-3)*vz^4 + (-10*ll^2)*vx^2 + 5*ll^2*vy^2 + 5*ll^2*vz^2
%%%    18   3*vx^2*vy^2 + 3*vy^4 + (-3)*vx^2*vz^2 + (-3)*vz^4 + (-5*ll^2)*vy^2 + 5*ll^2*vz^2

%%%%%%%%%%%%%%%%%%%%%%%%%%%%%%%%%%%%%%%%%%%%%%%%%%%%%%%%%%%%%%%%%%%%%%%%%%%%%%%%%%%

\newpage 
%%%%%%%%%%%%%%%%%%%%%%%%%%%%%%%%%%%%%%%%%%%%%%%%%%%%%%%%%%%%%%%%%%%%%%%%%%%%%%% 
\monitem  D3Q27 for isothermal Navier Stokes 

\smallskip

%%%%%%%%%%%%%%%%%%%%%%%%%%%%%%%%%%%%%%%%%%%%%%%%%%%%%%%%%%%%%%%%%%%%%%%%%%%%%%%%%%% moments d3q27  
\begin{table}    [H]  \centering
\centerline{
\begin{tabular}{|l|c|c|} \hline
$ k $ &  $m_k$ &  $p_k (v)$ \\  \hline
$ 0 $ &  $\rho$ &  $ 1 $ \\  
$ 1 $ &  $j_x $ &  $ v_x $ \\  
$ 2 $ &  $j_y $ &  $ v_y $ \\  
$ 3 $ &  $j_z $ &  $ v_z $ \\   \hline
$ 4 $ &  $ \varepsilon $ &  $  v_x^2 + v_y^2 + v_z^2 - 2 \, \lambda^2 $ \\  
$ 5 $ &  $ xx $ &  $  2 \, v_x^2 - v_y^2 - v_z^2  $ \\  
$ 6 $ &  $ ww $ &  $  v_y^2 - v_z^2 $ \\
\hline \end{tabular}}

\caption{Moment polynomials for  D3Q27 isothermal Navier Stokes}
\label{d3q27-iso-polynomes} \end{table}
%%%%%%%%%%%%%%%%%%%%%%%%%%%%%%%%%%%%%%%%%%%%%%%%%%%%%%%%%%%%%%%%%%%%%%%%%%%%%%%%%%%

\newpage 
%%%%%%%%%%%%%%%%%%%%%%%%%%%%%%%%%%%%%%%%%%%%%%%%%%%%%%%%%%%%%%%%%%%%%%%%%%%%%%% 
\monitem  D3Q33 for isothermal Navier Stokes 

\smallskip

%%%%%%%%%%%%%%%%%%%%%%%%%%%%%%%%%%%%%%%%%%%%%%%%%%%%%%%%%%%%%%%%%%%%%%%%%%%%%%%%%%% moments d3q33  
\begin{table}    [H]  \centering
\centerline{
\begin{tabular}{|l|c|c|} \hline
$ k $ &  $m_k$ &  $p_k (v)$ \\  \hline
$ 0 $ &  $\rho$ &  $ 1 $ \\  
$ 1 $ &  $j_x $ &  $ v_x $ \\  
$ 2 $ &  $j_y $ &  $ v_y $ \\  
$ 3 $ &  $j_z $ &  $ v_z $ \\   \hline
$ 4 $ &  $ \varepsilon $ &  $  11 \, (v_x^2 + v_y^2 + v_z^2)  - 26 \, \lambda^2 $ \\  
$ 5 $ &  $ xx $ &  $  2 \, v_x^2 - v_y^2 - v_z^2  $ \\  
$ 6 $ &  $ ww $ &  $  v_y^2 - v_z^2 $ \\
\hline \end{tabular}}

\caption{Moment polynomials for  D3Q33 isothermal Navier Stokes}
\label{d3q33-iso-polynomes} \end{table}

\newpage 
%%%%%%%%%%%%%%%%%%%%%%%%%%%%%%%%%%%%%%%%%%%%%%%%%%%%%%%%%%%%%%%%%%%%%%%%%%%%%%% 
\monitem D3Q27-2 for isothermal Navier Stokes 

\smallskip

%%%%%%%%%%%%%%%%%%%%%%%%%%%%%%%%%%%%%%%%%%%%%%%%%%%%%%%%%%%%%%%%%%%%%%%%%%%%%%%%%%% moments d3q27-2 
\begin{table}    [H]  \centering
\centerline{
\begin{tabular}{|l|c|c|} \hline
$ k $ &  $m_k$ &  $p_k (v)$ \\  \hline
$ 0 $ &  $\rho$ &  $ 1 $ \\  
$ 1 $ &  $j_x $ &  $ v_x $ \\  
$ 2 $ &  $j_y $ &  $ v_y $ \\  
$ 3 $ &  $j_z $ &  $ v_z $ \\   \hline
$ 4 $ &  $ \varepsilon $ &  $  3 \, (v_x^2 + v_y^2 + v_z^2)  - 8 \, \lambda^2 $ \\  
$ 5 $ &  $ xx $ &  $  2 \, v_x^2 - v_y^2 - v_z^2  $ \\  
$ 6 $ &  $ ww $ &  $  v_y^2 - v_z^2 $ \\
\hline \end{tabular}}

\caption{Moment polynomials for  D3Q27-2 isothermal Navier Stokes}
\label{d3q7-2-iso-polynomes} \end{table}
%%%%%%%%%%%%%%%%%%%%%%%%%%%%%%%%%%%%%%%%%%%%%%%%%%%%%%%%%%%%%%%%%%%%%%%%%%%%%%% 

\bigskip \newpage 
%%%%%%%%%%%%%%%%%%%%%%%%%%%%%%%%%%%%%%%%%%%%%%%%%%%%%%%%%%%%%%%%%%%%%%%%%%%%%%% 
\monitem D2Q13 for thermal Navier Stokes  

\smallskip

%%%%%%%%%%%%%%%%%%%%%%%%%%%%%%%%%%%%%%%%%%%%%%%%%%%%%%%%%%%%%%%%%%%%%%%%%%%%%%%%%%% moments d2q13  
\begin{table}    [H]  \centering
\centerline{
\begin{tabular}{|l|c|c|} \hline
$ k $ &  $m_k$ &  $p_k (v)$ \\  \hline
$ 0 $ &  $\rho$ &  $ 1 $ \\  
$ 1 $ &  $j_x $ &  $ v_x $ \\  
$ 2 $ &  $j_y $ &  $ v_y $ \\  
$ 3 $ &  $ \varepsilon $ &  $  13 \, (v_x^2 + v_y^2) - 28 \, \lambda^2 $ \\   \hline
$ 4 $ &  $ xx $ &  $  v_x^2 - v_y^2 $ \\  
$ 5 $ &  $ xy $ &  $  v_x \, v_y $ \\  
$ 6 $ &  $ q_x $ &  $  (v_x^2 + v_y^2 - 3  \, \lambda^2 )\, v_x  $ \\  
$ 7 $ &  $ q_y $ &  $  (v_x^2 + v_y^2 - 3  \, \lambda^2 )\, v_y  $ \\   \hline 
$ 8 $ &  $ r_x $ &  $  [ {35\over12} \, (v_x^2 + v_y^2)^2 - {63\over4} \,\lambda^2\,(v_x^2 + v_y^2) + {101\over6}\,\lambda^4 ] \, v_x  $ \\  
$ 9 $ &  $ r_y $ &  $  [ {35\over12} \, (v_x^2 + v_y^2)^2 - {63\over4} \,\lambda^2\,(v_x^2 + v_y^2) + {101\over6}\,\lambda^4 ] \, v_y  $ \\  
$ 10 $ &  $ h $ &  $   {77\over2} \, (v_x^2 + v_y^2)^2 -  {361\over2} \,\lambda^2\,(v_x^2 + v_y^2) + 140 \,\lambda^4   $ \\  
$ 11 $ &  $ xx_e  $ &  $   [ {17\over12} \, (v_x^2 + v_y^2) -  {65\over12} \,\lambda^2 ] \, (v_x^2 - v_y^2)  $ \\   \hline 
$ 12 $ &  $ h_3 $ &  $   {137\over24} \, (v_x^2 + v_y^2)^3 -  {273\over8} \,\lambda^2 \,(v_x^2 + v_y^2)^2   +  {581\over12} \,\lambda^4 \,(v_x^2 + v_y^2)
                         -12 \,\lambda^6      $ \\  
\hline \end{tabular}}

\caption{Moment polynomials for  D2Q13 thermal Navier Stokes}
\label{polynomes-d2q13-thermique} \end{table}
%%%%%%%%%%%%%%%%%%%%%%%%%%%%%%%%%%%%%%%%%%%%%%%%%%%%%%%%%%%%%%%%%%%%%%%%%%%%%%%%%%%

%%% \newpage 
%%%%%%%%%%%%%%%%%%%%%%%%%%%%%%%%%%%%%%%%%%%%%%%%%%%%%%%%%%%%%%%%%%%%%%%%%%%%%%% 
\monitem D2Q17 for thermal Navier Stokes  

\smallskip

%%%%%%%%%%%%%%%%%%%%%%%%%%%%%%%%%%%%%%%%%%%%%%%%%%%%%%%%%%%%%%%%%%%%%%%%%%%%%%%%%%% moments d2q17 
\begin{table}    [H]  \centering
\centerline{
\begin{tabular}{|l|c|c|} \hline
$ k $ &  $m_k$ &  $p_k (v)$ \\  \hline
$ 0 $ &  $\rho$ &  $ 1 $ \\  
$ 1 $ &  $j_x $ &  $ v_x $ \\  
$ 2 $ &  $j_y $ &  $ v_y $ \\  
$ 3 $ &  $ \varepsilon $ &  $  17 \, (v_x^2 + v_y^2) - 60 \, \lambda^2 $ \\   \hline
$ 4 $ &  $ xx $ &  $  v_x^2 - v_y^2 $ \\  
$ 5 $ &  $ xy $ &  $  v_x \, v_y $ \\  
$ 6 $ &  $ q_x $ &  $  [ 3\, (v_x^2 + v_y^2) - 17 \, \lambda^2 ] \, v_x  $ \\
\hline \end{tabular}}

\caption{Moment polynomials for  D2Q17 thermal Navier Stokes}
\label{d2q17-ns-polynomes} \end{table}
%%%%%%%%%%%%%%%%%%%%%%%%%%%%%%%%%%%%%%%%%%%%%%%%%%%%%%%%%%%%%%%%%%%%%%%%%%%%%%%%%%%

\newpage %%% \bigskip 
%%%%%%%%%%%%%%%%%%%%%%%%%%%%%%%%%%%%%%%%%%%%%%%%%%%%%%%%%%%%%%%%%%%%%%%%%%%%%%% 
\monitem D2V17 for thermal Navier Stokes 

\smallskip

%%%%%%%%%%%%%%%%%%%%%%%%%%%%%%%%%%%%%%%%%%%%%%%%%%%%%%%%%%%%%%%%%%%%%%%%%%%%%%%%%%% moments d2v17 
\begin{table}    [H]  \centering
\centerline{
\begin{tabular}{|l|c|c|} \hline
$ k $ &  $m_k$ &  $p_k (v)$ \\  \hline
$ 0 $ &  $\rho$ &  $ 1 $ \\  
$ 1 $ &  $j_x $ &  $ v_x $ \\  
$ 2 $ &  $j_y $ &  $ v_y $ \\  
$ 3 $ &  $ \varepsilon $ &  $  17 \, (v_x^2 + v_y^2) - 80 \, \lambda^2 $ \\    \hline
$ 4 $ &  $ xx $ &  $  v_x^2 - v_y^2 $ \\  
$ 5 $ &  $ xy $ &  $  v_x \, v_y $ \\  
$ 6 $ &  $ q_x $ &  $  [ 2 \, (v_x^2 + v_y^2) - 15 \, \lambda^2 ] \, v_x  $ \\
\hline \end{tabular}}
\caption{Moment polynomials for  D2V17 thermal Navier Stokes}
\label{d2v17-ns-polynomes} \end{table}

 \newpage 
%%%%%%%%%%%%%%%%%%%%%%%%%%%%%%%%%%%%%%%%%%%%%%%%%%%%%%%%%%%%%%%%%%%%%%%%%%%%%%% 
\monitem D2W17 for thermal Navier Stokes 

\smallskip

%%%%%%%%%%%%%%%%%%%%%%%%%%%%%%%%%%%%%%%%%%%%%%%%%%%%%%%%%%%%%%%%%%%%%%%%%%%%%%%%%%% moments d2w17 
\begin{table}    [H]  \centering
\centerline{
\begin{tabular}{|l|c|c|} \hline
$ k $ &  $m_k$ &  $p_k (v)$ \\  \hline
$ 0 $ &  $\rho$ &  $ 1 $ \\  
$ 1 $ &  $j_x $ &  $ v_x $ \\  
$ 2 $ &  $j_y $ &  $ v_y $ \\  
$ 3 $ &  $ \varepsilon $ &  $  17 \, (v_x^2 + v_y^2) - 52 \, \lambda^2 $ \\  \hline 
$ 4 $ &  $ xx $ &  $  v_x^2 - v_y^2 $ \\  
$ 5 $ &  $ xy $ &  $  v_x \, v_y $ \\  
$ 6 $ &  $ q_x $ &  $  [ 13 \, (v_x^2 + v_y^2) - 15 \, \lambda^2 ] \, v_x  $ \\  
$ 7 $ &  $ q_y $ &  $  [ 13 \, (v_x^2 + v_y^2) - 15 \, \lambda^2 ] \, v_y  $ \\   \hline 
$ 8 $ &  $ r_x $ &  $  [ {57\over4} \, (v_x^2 + v_y^2)^2 - {371\over4} \,(v_x^2 + v_y^2)\,\lambda^2 + {217\over2}\,\lambda^4 ] \, v_x  $ \\  
$ 9 $ &  $ r_y $ &  $  [ {57\over4} \, (v_x^2 + v_y^2)^2 - {371\over4} \,(v_x^2 + v_y^2)\,\lambda^2 + {217\over2}\,\lambda^4 ] \, v_y  $ \\  
$ 10 $ &  $ xy_2  $ &  $ {1\over6} \, v_x \, v_y^2 \, (v_x^2 - v_y^2 ) $ \\   
$ 11 $ &  $ yx_2  $ &  $ {1\over6} \, v_y \, v_x^2 \, (v_x^2 - v_y^2 ) $ \\  
$ 12 $ &  $ h $ &  $    {1\over6} \, (v_x^2 + v_y^2)^2 +   {1543\over2} \, \lambda^2 \,  (v_x^2 + v_y^2)^2 +  684 \, \lambda^4 $ \\
$ 13 $ &  $ xx_e $ &  $   [ {19\over12} \, (v_x^2 + v_y^2) - {91\over12} \,\lambda^2 ] \, (v_x^2 - v_y^2)  $ \\
$ 14 $ &  $ xy_e $ &  $   [ {3\over2} \, (v_x^2 + v_y^2) - 7 \,\lambda^2 ] \, v_x \, v_y    $ \\ \hline 
$ 15 $ &  $ xx_{xy} $ &  $   {1\over6} \, (v_x^2 - v_y^2) \, v_x \, v_y  $ \\  
$ 16 $ &  $ h_3  $ &  $  {313\over20} \, (v_x^2 + v_y^2)^3 -  {2219\over20}  \, (v_x^2 + v_y^2)^2 \,\lambda^2 +  {1733\over10}  \, (v_x^2 + v_y^2) \,\lambda^4
 - 48 \, \lambda^6 $ \\  
  \hline \end{tabular}}

\caption{Moment polynomials for  D2W17 thermal Navier Stokes}
\label{d2w17-ns-polynomes} \end{table}

\newpage 
%%%%%%%%%%%%%%%%%%%%%%%%%%%%%%%%%%%%%%%%%%%%%%%%%%%%%%%%%%%%%%%%%%%%%%%%%%%%%%% 
\monitem D3Q33 for thermal Navier Stokes 

\smallskip

%%%%%%%%%%%%%%%%%%%%%%%%%%%%%%%%%%%%%%%%%%%%%%%%%%%%%%%%%%%%%%%%%%%%%%%%%%%%%%%%%%% moments d3q33 
\begin{table}    [H]  \centering
\centerline{
\begin{tabular}{|l|c|c|} \hline
$ k $ &  $m_k$ &  $p_k (v)$ \\  \hline
$ 0 $ &  $\rho$ &  $ 1 $ \\  
$ 1 $ &  $j_x $ &  $ v_x $ \\  
$ 2 $ &  $j_y $ &  $ v_y $ \\  
$ 3 $ &  $j_z $ &  $ v_z $ \\   
$ 4 $ &  $ \varepsilon $ &  $  11 \, (v_x^2 + v_y^2 + v_z^2)  - 26 \, \lambda^2 $ \\  \hline
$ 5 $ &  $ xx $ &  $  2 \, v_x^2 - v_y^2 - v_z^2  $ \\  
$ 6 $ &  $ ww $ &  $  v_y^2 - v_z^2 $ \\
\hline \end{tabular}}

\caption{Moment polynomials for  D3Q33 thermal Navier Stokes}
\label{d3q33-ns-polynomes}\end{table}
%%%%%%%%%%%%%%%%%%%%%%%%%%%%%%%%%%%%%%%%%%%%%%%%%%%%%%%%%%%%%%%%%%%%%%%%%%%%%%%%%%%

\newpage 
%%%%%%%%%%%%%%%%%%%%%%%%%%%%%%%%%%%%%%%%%%%%%%%%%%%%%%%%%%%%%%%%%%%%%%%%%%%%%%% 
\monitem D3Q27-2 for thermal Navier Stokes 

\smallskip

%%%%%%%%%%%%%%%%%%%%%%%%%%%%%%%%%%%%%%%%%%%%%%%%%%%%%%%%%%%%%%%%%%%%%%%%%%%%%%%%%%% moments d3q27-2-ns 
\begin{table}    [H]  \centering
\centerline{
\begin{tabular}{|l|c|c|} \hline
$ k $ &  $m_k$ &  $p_k (v)$ \\  \hline
$ 0 $ &  $\rho$ &  $ 1 $ \\  
$ 1 $ &  $j_x $ &  $ v_x $ \\  
$ 2 $ &  $j_y $ &  $ v_y $ \\  
$ 3 $ &  $j_z $ &  $ v_z $ \\  
$ 4 $ &  $ \varepsilon $ &  $  3 \, (v_x^2 + v_y^2 + v_z^2)  - 8 \, \lambda^2 $ \\   \hline
$ 5 $ &  $ xx $ &  $  2 \, v_x^2 - v_y^2 - v_z^2  $ \\  
$ 6 $ &  $ ww $ &  $  v_y^2 - v_z^2 $ \\
\hline \end{tabular}}

\caption{Moment polynomials for  D3Q27-2 thermal Navier Stokes}
\label{d3q27-2-ns-polynomes}\end{table}
%%%%%%%%%%%%%%%%%%%%%%%%%%%%%%%%%%%%%%%%%%%%%%%%%%%%%%%%%%%%%%%%%%%%%%%%%%%%%%%%%%%

%
%%%%%%%%%%%%%%%%%%%%%%%%%%%%%%%%%%%%%%%%%%%%%%%%%%%%%%%%%%%%%%%%%%%%%%%%%%%%  merci
\bigskip \bigskip   \noindent {\bf  \large  Acknowledgments }
%%%%%%%%%%%%%%%%%%%%%%%%%%%%%%%%%%%%%%%%%%%%%%%%%%%%%%%%%%%%%%%%%%%%%%%%%%%%%%%%%%%%%%%%%%

\noindent
FD  thanks the Centre National de  la Recherche Scientifique for according a ``Delegation'' 
at the International Research Laboratory 3457   in the ``Centre de Recherches Math\'ematiques''
of the Universit\'e de Montr\'eal during the period February-July 2021.
A part of this contribution was done during this period. 
%
%%%  {\color{blue}
The authors thank also Bruce Boghosian,  Benjamin Graille, Paulo Cesar Philippi, Olivier Pironneau
and Mahdi Tekitek  for fundamental discussions.  %%  }
%
%%%  {\color{red}
The scientific exchanges with the referees greatly improved the matter of this contribution.
We thank them for their rigorous reading!  %% } 
Last but not least, all the formal computations have been implemented with SageMath \cite{sagemath}.

%%%%%%%%%%%%%%%%%%%%%%%%%%%%%%%%%%%%%%%%%%%%%%%%%%%%%%%%%%%%%%%%%%%%%%%%%%%%%%%%%%%%%%%%%%

\newpage
%%%%%%%%%%%%%%%%%%%%%%%%%%%%%%%%%%%%%%%%%%%%%%%%%%%%%%%%%%%%%%%%%%%%%%%%%%%%  references
\bigskip \bigskip      \noindent {\bf  \large  References }
%%%%%%%%%%%%%%%%%%%%%%%%%%%%%%%%%%%%%%%%%%%%%%%%%%%%%%%%%%%%%%%%%%%%%%%%%%%%%%%%%%%%%%%%%%

%%%%%%%%%%%%%%%%%%%%%%%%%%%%%%%%%%       31 decembre 2018       %%%%%%%%%%%%%%%%%%%%%%%%%
\fancyhead[EC]{\sc{Fran\c{c}ois Dubois and Pierre Lallemand}}
\fancyhead[OC]{\sc{Single lattice Boltzmann distribution for Navier Stokes equations}}
%%%%%%%%%%%%%%%%%%%%%%%%%%%%%%%%%%       31 decembre 2018       %%%%%%%%%%%%%%%%%%%%%%%%%
%%%%%%%%%%%%%%%%%%%%%%%%%%%%%%%%%%%%%%%%%%%%%  jolie numerotation des pages
\fancyfoot[C]{\oldstylenums{\thepage}}
%%%%%%%%%%%%%%%%%%%%%%%%%%%%%%%%%%%%%%%%%%%%%  fin jolie numerotation des pages

 %%%%%%%%  \vspace{-.24cm}

%%%%%%%%%%%%%%%%%%%%%%%%%%%%%%%%%%%%%%%%%%%%%%%%%%%%%%%%%%%%%%%%%%%%%%%%%%%%%%%%%%%%%%%%%%%%%%%%%%%%%%%%%%%%%%%%%%%%%%
%%%%%%%%%%%%%%%%%%%%%%%%%%%%%%%%%%%%%%%%%%%%%%%%%%%%%%%%%%%%%%%%%%%%%%%%%%%%%%%%%%%%%%%%%%%%%%%%%%%%%%%%%%%%%%%%%%%%%%
%%%%%%%%%%%%%%%%%%%%%%%%%%%%%%%%%%%%%%%%%%%%%%%%%%%%%%%%%%%%%%%%%%%%%%%%%%%%%%%%%%%%%%%%%%%%%%%%%%%%%%%%%%%%%%%%%%%%%%
%%%%%%%%%%%%%%%%%%%%%%%%%%%%%%%%%%%%%%%%%%%%%%%%%%%%%%%%%%%%%%%%%%%%%%%%%%%%%%%%%%%%%%%%%%%%%%%%%%%%%%%%%%%%%%%%%%%%%%
%%%%%%%%%%%%%%%%%%%%%%%%%%%%%%%%%%%%%%%%%%%%%%%%%%%%%%%%%%%%%%%%%%%%%%%%%%%%%%%%%%%%%%%%%%%%%%%%%%%%%%%%%%%%%%%%%%%%%%
%%%%%%%%%%%%%%%%%%%%%%%%%%%%%%%%%%%%%%%%%%%%%%%%%%%%%%%%%%%%%%%%%%%%%%%%%%%%%%%%%%%%%%%%%%%%%%%%%%%%%%%%%%%%%%%%%%%%%%

\end{document}